\documentclass[10pt,a4paper]{amsart} 
\usepackage{amsfonts,amssymb,amsthm,amsmath,ifpdf,stmaryrd,amsrefs,dsfont,esint,mathrsfs,dsfont,xparse,arydshln}
\usepackage[hidelinks]{hyperref}
\usepackage{enumitem} 
\usepackage[ddmmyyyy]{datetime}

\usepackage{tikz}
\usetikzlibrary{matrix,arrows}
\usepackage{tikz-cd}

\usepackage[overload]{textcase}

\ifpdf
  \usepackage[final,expansion=alltext,protrusion=alltext]{microtype} 
\fi

\usepackage{xargs} 
\usepackage{mathtools} 
\usepackage{ifthen} 

\theoremstyle{plain}

\newtheorem*{Th*}{Theorem}
\newtheorem*{ThA}{Superorbit Theorem}
\newtheorem*{ThB}{Supersymplectic Orbit Form Theorem}

\newtheorem*{Cor*}{Corollary}

\theoremstyle{definition}

\theoremstyle{remark}

\newif\ifTNS 
\TNStrue

\def\printtheoremname#1{\csname#1name\endcsname}
\def\printtheoremnames#1{\csname#1names\endcsname}

\def\thmref#1#2{\printtheoremname{#1}\ifTNS~\fi\ref{#1:#2}}

\def\uc#1#2{\MakeUppercase{#1}{#2}} 

\newcommand{\DefTheorem}[2]{\newenvironmentx{#1}[2][1=\empty,2=\empty]{%
    \ignorespaces%
    \ifx##2\empty%
      \begin{#2}%
    \else%
      \begin{#2}[{\uc##2}]%
    \fi%
    \ifx##1\empty%
      {}%
    \else%
      \label{#1:##1}%
    \fi%
    \ignorespaces}{\end{#2}\ignorespacesafterend}}

\newcommand{\prfof}[2]{\protect{Proof of~\thmref{#1}{#2}}}

\makeatother

\DefTheorem{Th}{theorem}
\DefTheorem{Prop}{proposition}
\DefTheorem{Cor}{corollary}
\DefTheorem{Lem}{lemma}
\DefTheorem{Def}{definition}
\DefTheorem{Rem}{remark}
\DefTheorem{Par}{para}
\DefTheorem{Not}{notation}
\DefTheorem{Exer}{exercise}
\DefTheorem{Ex}{example}
\DefTheorem{Cons}{construction}
\DefTheorem{Scho}{scholium}

\newenvironment{Par*}{\ignorespaces\noindent\ignorespaces}{\ignorespacesafterend}

\numberwithin{equation}{section}

\newcommand\Define[2][\empty]{\ignorespaces%
  \emph{#2}}%

\tikzset{
  commutative diagrams/.cd,
  arrow style=tikz,
  diagrams={>=stealth},
  shift up/.style={
    to path={([yshift=#1]\tikztostart.east) -- ([yshift=#1]\tikztotarget.west) \tikztonodes}},
  shift up left/.style={
    to path={([yshift=#1]\tikztostart.west) -- ([yshift=#1]\tikztotarget.east) \tikztonodes}},
  mathdouble/.style={-,double equal sign distance}
}

\def\ger{\mathfrak}

\newcommand\CategoryTypeface{\mathbf}
\def\cat{\CategoryTypeface}

\newcommand\SheafTypeface{\mathcal}
\def\sh{\SheafTypeface}

\def\DMO{\DeclareMathOperator}

\newcommand\ev{{\bar 0}}
\newcommand\odd{{\bar 1}}

\newcommand{\defi}{\coloneqq}     

\def\diff#1^#2{\ensuremath{\partial_{#1}^{#2}}}

\def\der#1/#2{\ifthenelse {\equal{#1}{}}
              {\ensuremath{\partial_{#2}{#1}}}
              {\ensuremath{\frac{\partial #1}{\partial #2}}}
        }

\def\derf#1/#2{\ifthenelse  {\equal{#1}{}}
              {\ensuremath{\frac{\partial #1}{\partial #2}}}
              {\ensuremath{\partial_{#2}{#1}}}
        }

\newcommand{\Fa}{For all }
\newcommand{\fa}{for all }

\newcommand{\fs}{for some }

\newcommand{\scth}{such that }
\newcommand{\AND}{and}
\newcommand{\OR}{or}
\newcommand\mathtxt[1]{\quad\text{{#1}}\quad}

\newcommand{\nd}{\mathtxt{\AND}}

\newcommand{\Or}{\mathtxt{\OR}}

\newcommand\cf{\emph{cf.}~}

\newcommand\eg{\emph{e.g.}~}
\newcommand\ie{\emph{i.e.}~}

\newcommand\via{\emph{via}~}
\newcommand\loccit{\emph{loc.~cit.}}
\newcommand\opcit{\emph{op.~cit.}}

\newcommand\etc{\emph{etc.}}

\newcommand\vphi{\varphi}

\newcommand\vkappa{\varkappa}
\newcommand\eps{\varepsilon}

\newcommand\ints{\mathbb{Z}}

\newcommand\reals{\mathbb{R}}
\newcommand\cplxs{\mathbb{C}}
\newcommand\knums{\mathbb K}

\def\aff{{\mathbb A}}

\newcommand\sle{\leqslant}
\newcommand\sge{\geqslant}

\DMO\dom{\mathrm{dom}}
\DMO\rk{\mathrm{rk}}
\DMO\Ad{\mathrm{Ad}}
\DMO\ad{\mathrm{ad}}
\DMO\GL{\mathrm{GL}}
\DMO\id{\mathrm{id}}
\DMO\pr{\mathrm{pr}}
\DMO\gr{\mathrm{gr}}
\DMO\sll{\ger{sl}}
\DMO\sdim{\mathrm{sdim}}
\DMO\sgn{\mathrm{sgn}}
\DMO\re{\mathrm{Re}}
\DMO\Gal{\mathrm{Gal}}
\DMO\Ann{\mathrm{Ann}}
\DMO\coker{\mathrm{coker}}
\DMO\im{\mathrm{im}}
\DMO\coim{\mathrm{coim}}
\DMO\Spec{\mathrm{Spec}}
\DMO\codim{\mathrm{codim}}
\DMO\chr{\mathrm{char}}
\DMO\supp{\mathrm{supp}}
\DMO\str{\mathrm{str}}
\DMO\tr{\mathrm{tr}}

\DMO\longrightarrowp{\cat{Top}}
\DMO\Sets{\cat{Sets}}
\DMO\SMan{\cat{SMan}}
\DMO\SRSp{\cat{SRSp}}
\DMO\SSp{\cat{SSp}}
\DMO\SVSp{\cat{SVSp}}
\DMO\SVec{\cat{SVec}}
\DMO\SAlg{\cat{SAlg}}
\DMO\Mod{\cat{Mod}}
\DMO\Op{\cat{Op}}
\DMO\Cov{\cat{Cov}}

\DMO\Ob{\mathrm{Ob}}

\newcommand\ssplfg[2][\empty]{%
  \smash{\SSp
    \ifx#1\empty%
      ^{\mathrm{lfg}}_{#2}
    \else
      ^{#1,\mathrm{lfg}}_{#2}
    \fi
  }
}

\newcommand{\lBr}{[\kern-.35ex[}
\newcommand{\rBr}{]\kern-.35ex]}

\makeatletter
\newcommand\Size[7][1]{
                                 \ifx#20%
                                        \def\r@l{}\def\r@m{}\def\r@r{}%
                                 \else%
                                    \ifx#21%
                                           \def\r@l{\bigl}\def\r@r{\bigr}\def\r@m{\bigm}%
                                    \else%
                                           \ifx#22%
                                                 \def\r@l{\Bigl}\def\r@r{\Bigr}\def\r@m{\Bigm}%
                                            \else%
                                                 \ifx#23%
                                                        \def\r@l{\biggl}\def\r@r{\biggr}\def\r@m{\biggm}%
                                                  \else
                                                        \ifx#24%
                                                        \def\r@l{\Biggl}\def\r@r{\Biggr}\def\r@m{\Biggm}%
                                                        \fi%
                                                  \fi%
                                            \fi%
                                      \fi%
                                 \fi%
                                 \ifx#10%
                                       \def\r@m{}%
                                 \fi%
                                 \r@l#3{#4}\r@m#5{#6}\r@r#7%
}%

\makeatother

\makeatletter
\def\Set@Scallop[#1]#2#3{{#1}\Parens{#2}{#3}}
\newcommand\DeclareScalableOperator[2]{%
  \expandafter\def\csname#1\endcsname{\@ifnextchar[{{#2}\Set@Scallop}{{#2}\Set@Scallop[{}]}}
}
\makeatother

\def\DSO{\DeclareScalableOperator}

\DSO{Presh}{\cat{Presh}}
\DSO{Sh}{\cat{Sh}}

\DSO{ShHom}{\sh Hom} 
\DSO{ShGHom}{\underline{\sh Hom}} 
\DSO{ShEnd}{\sh End} 
\DSO{ShGEnd}{\underline{\sh End}} 
\DSO{ShDer}{\sh Der} 
\DSO{ShGDer}{\underline{\sh Der}} 
\DSO{ShCliff}{\sh C\text{\emph{liff}}} 

\DSO{Ct}{\mathcal{C}} 
\DSO{Cc}{\mathcal{C}_c} 
\DSO{Lp}{\mathbf{L}} 
\DSO{Hom}{\mathrm{Hom}} 
\DSO{End}{\mathrm{End}} 
\DSO{Aut}{\mathrm{Aut}} 
\DSO{Uenv}{\mathfrak{U}} 
\DSO{Cuenv}{\widehat{\mathfrak{U}}} 
\DSO{ABer}{\Abs0{\mathrm{Ber}}}
\DSO{Proj}{\mathrm{Proj}}
\DSO{Ber}{\mathrm{Ber}}
\DSO{Vol}{\mathrm{Vol}}
\DSO{Form}{\Omega}
\DSO{Ind}{\mathrm{Ind}}
\DSO{Coind}{\mathrm{Coind}}
\DSO{Der}{\mathrm{Der}}
\DSO{GDer}{\underline{\mathrm{Der}}}
\DSO{Alg}{\mathrm{Alg}}
\DSO{GHom}{\underline{\mathrm{Hom}}} 
\DSO{GEnd}{\underline{\mathrm{End}}} 
\DSO{Maps}{\mathrm{Maps}}
\DSO{Cliff}{\mathrm{Cliff}} 
\DSO{Jac}{\mathrm{Jac}} 
\DSO{Db}{\mathcal{D}'} 

\newcommand\Set[3]{
                                 \Size{#1}{\{}{#2}{|}{#3}{\}}%
}%
\newcommand\Sdual[3]{
                                 \Size{#1}{\langle}{#2}{|}{#3}{\rangle}%
}%
\newcommand\Dual[3]{
                                 \Size[0]{#1}{\langle}{#2}{,}{#3}{\rangle}%
}%
\newcommand\Parens[2]{
  \Size[0]{#1}{(}{#2}{}{}{)}
}
\newcommand\Bracks[2]{
  \Size[0]{#1}{[}{#2}{}{}{]}
}

\newcommand\Norm[2]{
  \Size[0]{#1}{\lVert}{#2}{}{}{\rVert}
}
\newcommand\Abs[2]{
  \Size[0]{#1}{\lvert}{#2}{}{}{\rvert}
}
\newcommand\Span[2]{
  \Size[0]{#1}{\langle}{#2}{}{}{\rangle}
}
\makeatletter

\makeatletter

\newif\if@smallmat
\newif\if@none
\newif\if@paren
\newif\if@brack
\newif\if@brace
\newif\if@vline
\newenvironment{Matrix}[2][1]
                                 {\ifx#20%
                                        \@smallmattrue%
                                  \else%
                                         \@smallmatfalse
                                  \fi%
                                  \ifx#11%
                                         \@nonefalse\@parentrue\@brackfalse\@bracefalse\@vlinefalse%
                                  \else%
                                       \ifx#12%
                                            \@nonefalse\@parenfalse\@bracktrue\@bracefalse\@vlinefalse%
                                        \else%
                                            \ifx#13%
                                                 \@nonefalse\@parenfalse\@brackfalse\@bracetrue\@vlinefalse%
                                            \else%
                                                 \ifx#14%
                                                       \@nonefalse\@parenfalse\@brackfalse\@bracefalse\@vlinetrue
                                                 \else%
                                                       \ifx#15%
                                                             \@nonefalse\@parenfalse\@brackfalse\@bracefalse\@vlinefalse%
                                                       \else%
                                                             \@nonetrue\@parenfalse\@brackfalse\@bracefalse\@vlinefalse%
                                                       \fi%
                                                 \fi%
                                            \fi%
                                        \fi%
                                   \fi%
                                   \if@smallmat%
                                        \if@none%
                                             \begin{smallmatrix}%
                                        \else%
                                            \if@paren%
                                                  \left(\begin{smallmatrix}%
                                            \else%
                                                  \if@brack%
                                                          \left[\begin{smallmatrix}%
                                                  \else%
                                                          \if@brace%
                                                               \left\{\begin{smallmatrix}%
                                                          \else%
                                                               \if@vline%
                                                                    \left\lvert\begin{smallmatrix}%
                                                                \else%
                                                                    \left\lVert\begin{smallmatrix}%
                                                                \fi%
                                                          \fi%
                                                  \fi%
                                            \fi%
                                        \fi%
                                   \else%
                                        \if@none%
                                             \begin{matrix}%
                                        \else%
                                            \if@paren%
                                                  \begin{pmatrix}%
                                            \else%
                                                  \if@brack%
                                                          \begin{bmatrix}%
                                                  \else%
                                                          \if@brace%
                                                               \begin{Bmatrix}%
                                                          \else%
                                                               \if@vline%
                                                                    \begin{vmatrix}%
                                                                \else%
                                                                    \begin{Vmatrix}%
                                                                \fi%
                                                          \fi%
                                                  \fi%
                                            \fi%
                                        \fi%
                                   \fi}%
                                  {\if@smallmat%
                                        \if@none%
                                             \end{smallmatrix}%
                                        \else%
                                            \if@paren%
                                                  \end{smallmatrix}\right)%
                                            \else%
                                                  \if@brack%
                                                          \end{smallmatrix}\right]%
                                                  \else%
                                                          \if@brace%
                                                               \end{smallmatrix}\right\}%
                                                          \else%
                                                               \if@vline%
                                                                    \end{smallmatrix}\right\rvert%
                                                                \else%
                                                                    \end{smallmatrix}\right\rVert%
                                                                \fi%
                                                          \fi%
                                                  \fi%
                                            \fi%
                                         \fi%
                                   \else%
                                        \if@none%
                                             \end{matrix}%
                                        \else%
                                            \if@paren%
                                                  \end{pmatrix}%
                                            \else%
                                                  \if@brack%
                                                          \end{bmatrix}%
                                                  \else%
                                                          \if@brace%
                                                               \end{Bmatrix}%
                                                          \else%
                                                               \if@vline%
                                                                    \end{vmatrix}%
                                                                \else%
                                                                    \end{Vmatrix}%
                                                                \fi%
                                                          \fi%
                                                  \fi%
                                            \fi%
                                        \fi%
                                   \fi}%
                         
\makeatother

\def\clap#1{\hbox to 0pt{\hss#1\hss}} 

\newif\iffullproof
\fullprooftrue

\def\IndexL#1#2#3{}
\def\IndexG#1#2#3{}
\def\IndexO#1#2#3{}

\begin{document}

\title{Superorbits}

\author[Alldridge]
{Alexander Alldridge}

\address{Universit\"at zu K\"oln\\
Mathematisches Institut\\
Weyertal 86-90\\
50931 K\"oln\\
Germany}
\email{alldridg@math.uni-koeln.de}

\author[Hilgert]
{Joachim Hilgert}

\address{Universit\"at Paderborn\\
Institut f\"ur Mathematik\\
33095 Paderborn\\
Germany}
\email{hilgert@math.upb.de}

\author[Wurzbacher]
{Tilmann Wurzbacher}

\address{%
Institut \'E.~Cartan (IECL)\\
Universit\'e de Lorraine et
C.N.R.S.\\
57045 Metz, France}

\email{tilmann.wurzbacher@univ-lorraine.fr}

\thanks{Research supported by Deutsche Forschungsgemeinschaft (DFG), grant nos.~SFB/TR 12 (all authors), the Heisenberg grant AL 698/3-1 (A.A.), the Leibniz prize to M.~Zirnbauer ZI 513/2-1 (A.A.), SFB TRR 183 (A.A.), and the Institutional Strategy of the University of Cologne within the German Excellence Initiative (all authors)}



\begin{abstract}
  We study actions of Lie supergroups, in particular, the hitherto elusive notion of orbits through odd (or more general) points. Following categorical principles, we derive a conceptual framework for their treatment and therein prove general existence theorems for the isotropy (or stabiliser) supergroups and orbits through general points. In this setting, we show that the coadjoint orbits always admit a (relative) supersymplectic structure of Kirillov--Kostant--Souriau type. Applying a family version of Kirillov's orbit method, we decompose the regular representation of an odd Abelian supergroup into an odd direct integral of characters and construct universal families of representations, parametrised by a supermanifold, for two different super variants of the Heisenberg group.
\end{abstract}

\maketitle

\section{Introduction}

The present formulation of the theory of actions and representations of Lie supergroups does not appropriately address all relevant phenomena: Consider the basic example of the additive Lie supergroup $G$ of an odd-super vector space $\ger g$. The coadjoint action is trivial, so the orbit through the unique point $0\in\ger g^*$ is again a point. Similarly, $G$ has only the trivial irreducible unitary representation. Although this confirms the idea of the orbit method in a narrow sense, there is no hope of decomposing the regular representation of $G$ on $\sh O_G=\bigwedge\ger g^*$ by these means, nor can one reasonably expect thereby to construct representations of $G$ in any generality.

This suggests that it is crucial to broaden the notion of points. Following A.~Grothendieck, a \Define{$T$-valued point} of a space $X$ is a map $x:T\longrightarrow X$. This idea is based on considering an ordinary point as a map $*\longrightarrow X$ where $*$ is a singleton, allowing the parameter space to acquire additional degrees of freedom. The $G$-isotropy (or stabiliser) through $x$ should then be a `group bundle' $G_x\longrightarrow T$, and the orbit a `bundle' $G\cdot x\longrightarrow T$ with a fibrewise $G$-action. 

For any Lie supergroup $G$ with Lie superalgebra $\ger g$ acting on a supermanifold $X$ and any $x:T\longrightarrow X$, we obtain the following.

\begin{ThA}
    The isotropy supergroup $G_x$ exists as a Lie supergroup over $T$ if and only if the orbit morphism is of locally constant rank, which is the case if and only if the $\sh O_T$-module $x^*(\sh A_\ger g)$ is a locally direct summand of $x^*(\sh T_X)$. Here, $\sh A_\ger g$ is the \Define{fundamental distribution} generated by the fundamental vector fields. 

    Moreover, in this case, the orbit $G\cdot x\longrightarrow T\times X$ through $x$ exists as an equivariant local embedding of supermanifolds over $T$.
\end{ThA}

For the special case of orbits through ordinary points, the Superorbit Theorem was first proved by B.~Kostant \cite{kostant} in the setting of Lie--Hopf algebras, by C.P.~Boyer and O.A.~S\'anchez-Valenzuela \cite{bsv} for differentiable Lie supergroups, and by L.~Balduzzi, C.~Carmeli, and G.~Cassinelli \cite{bcc} using a functorial framework and super Harish-Chandra pairs. We recover the case of usual orbits through ordinary points as a special case. 

In the case of the coadjoint action of $G$ on $\ger g^*$ and of a $T$-valued point $f$ of $\ger g^*$, we prove the following result.

\begin{ThB}
  If $G_f$ exists as a Lie supergroup, then the coadjoint orbit $G\cdot f$ admits a canonical supersymplectic structure over $T$.
\end{ThB}

We stress that our point of view allows us to stay within the realm of \emph{even} supersymplectic forms, whereas in previous work \cites{tuyn10a,tuyn10b}, it was necessary to work with inhomogeneous symplectic forms. 

Furthermore, we introduce a general framework of supergroup representations over $T$ to extend Kirillov's method  \cite{kirillov} to orbits through $T$-valued points. As a proof of concept, we apply this to derive a Plancherel formula for the odd Abelian supergroup $\ger g$, presenting its regular representation as an `odd direct integral' of `unitary' characters. In a similar vein, we construct representations for two super versions of the three-dimensional Heisenberg group which arise by assigning suitable parities to the generators in the commutation relation $[x,y]=z$. In this case, we find `universal' parameter spaces $T$ and `universal' representations over $T$. Not surprisingly, these bear a striking similarity to the Schr\"odinger representation.

\medskip\noindent
The idea that irreducible representations should be constructed from orbits on some universal $G$-space is suggested by the general philosophy of geometric quantisation. The case where this works best is that of nilpotent Lie groups, where it was established by A.A.~Kirillov in the form of his orbit method.

The goal of extending this method to Lie supergroups was first addressed by B.~Kostant, in his seminal paper \cite{kostant}. In fact, as he remarks in his note \cite{kostant-harm}: Lie supergroups are ``likely to be [\ldots] useful [objects] only insofar as one can develop a corresponding theory of harmonic analysis''. Similarly, V.~Kac \cite{k77}*{5.5.4} poses the problem of constructing Lie supergroup representations \via the orbit method, in particular infinite-dimensional ones. For nilpotent Lie supergroups through \emph{ordinary} points, it was shown by H.~Salmasian \cite{salmasian} (and further investigated by Neeb--Salmasian \cite{ns11}) that indeed, there is a one-to-one correspondence of coadjoint orbits through ordinary points, \ie through elements of $\ger g_\ev^*$, with irreducible unitary representations in the sense of Varadarajan \emph{et al.} \cites{cctv,cctv-err}.

As remarked at the beginning of this introduction, this does not yet attain the goal of a theory of harmonic analysis for Lie supergroups, even in the Abelian case. These limitations are overcome by considering orbits through $T$-valued points.

A framework for the study of orbits through $T$-valued points was formulated in the category of schemes by D.~Mumford in his influential monograph \cite{mfk}, based on foundational work by A.~Grothendieck and P.~Gabriel. Although these ideas remain fruitful, the algebraic theory cannot be simply transferred to the differentiable category, and indeed the technical obstructions are formidable. At the same time, the differentiable setting is necessary for the envisaged applications: While all Lie groups are real analytic, any non-analytic (complete) vector field gives rise to an action which is not analytic (much less algebraic). Such situations are ubiquitous, particularly in the context of solvable Lie groups and their super generalisations.

A careful study of coadjoint orbits (through regular semi-simple elements) of the orthosymplectic and special linear supergroups in the algebraic category was conducted by R.~Fioresi and M.A.~Lled\'o in Ref.~\cite{fl}. The first one to consider coadjoint orbits through non-even functionals was G.~Tuynman \cites{tuyn10a,tuyn10b} in the form of a case study. His considerations are geared toward a specific example and formulated for DeWitt type supermanifolds. It is not clear whether this can be built into a general procedure and translated to Berezin--Kostant--Leites supermanifolds. Moreover, in his approach, he has to consider inhomogeneous ``symplectic'' forms. 

\medskip\noindent
We conclude the introduction by summarising the paper's contents. We present general categorical notions for the study of actions in Section \ref{sec:cats}. We emphasize the technique of base change known from algebraic geometry. This allows, among other things, to give a general definition of isotropy (or stabiliser) groups at $T$-valued points. In Section \ref{sec:groupoid-quots}, we review categorical quotients in the setting of differentiable and analytic superspaces and suggest a weak notion of geometric quotients. In order to treat quotients by group actions and equivalence relations on an equal footing, and with a view toward future applications, we introduce and employ the language of groupoids and their quotients. In Section \ref{sec:super-quot}, we specialise the discussion to supermanifolds. We prepare our discussion of isotropy supergroups at $T$-valued points by generalising the notion of morphisms of constant rank to relative supermanifolds (over a possibly singular base). We prove a rank theorem in this context (\thmref{Prop}{constant-rank-fibprod}); this is based on a family version of the inverse function theorem presented in Appendix \ref{app:invfun} (\thmref{Th}{invfun-loc}), also valid over a singular base. We investigate when the orbit morphism through a general point has constant rank (\thmref{Th}{action-locconst}) and, as an application, show the representability of isotropy supergroups under general conditions (\thmref{Th}{trans-iso}). This gives the existence of orbits under the same assumptions (\thmref{Th}{orbit}) and also implies that the isotropy supergroups exist only if the orbit morphism has constant rank. This relies on a family version of the closed subgroup theorem that we prove in Appendix \ref{app:closed-subgrp} (\thmref{Th}{imm-lie}). In Section \ref{sec:coad}, we construct the relative Kirillov--Kostant--Souriau form for coadjoint orbits through general points (\thmref{Th}{coadj-sympl}). Finally, in Section \ref{sec:quant}, we define the concept of representations over $T$. We then decompose the left-regular representation $\aff^{0|n}$ as a direct integral of characters and construct representations over appropriate parameter superspaces $T$ for super variants of the Heisenberg group. 

\medskip\noindent
\emph{Acknowledgements.} We gratefully acknowledge the hospitality of the Max-Planck Institute for Mathematics in Bonn, where much of the work on this article was done. We wish to thank Torsten Wedhorn for helpful discussions on module sheaves. 

\section{A categorical framework for group actions}\label{sec:cats}

\subsection{Categorical groups and actions}\label{subs:catgrp}

Groups and actions can be defined quite generally for categories with finite products. In this subsection, we recall the relevant notions and give a number of examples from differents contexts, which will serve to illustrate our further elaborations. 

In what follows, let $\cat C$ be a category with a terminal object $*$. For any $S,T\in\Ob\cat C$, let $\cat C_T^S$ be the category of objects in $\cat C$, which are under $S$ and over $T$. That is, objects and morphisms are given by the commutative diagrams depicted below:
\[
  \begin{tikzcd}[row sep=small]
    S\dar{}&S\dar{}\rar[mathdouble]{}&S\dar{}\\
    X\dar{}&X\rar{}\dar{}&Y\dar{}\\
    T&T\rar[mathdouble]{}&T.
  \end{tikzcd}  
\]
Similarly, we define the categories $\cat C_T$ of objects over $T$ and $\cat C^S$ of objects under $S$. 

We recall the definition of group objects and actions. These concepts are well-known, see \eg Ref.~\cite{maclane}. If $X,S\in\Ob\cat C$, then we write $x\in_SX$ for the statement `$x:S\longrightarrow X$ is a morphism in $\cat C$'. We also say `$x$ is an $S$-valued point of $X$' and denote the set of all these by $X(S)$. This defines the object map of the \Define{point functor} $X(-)$ of $X$. For a morphism $f:X\longrightarrow Y$ in $\cat C$ and $x\in_SX$, we define $f(x)\defi f\circ x$. Applying this procedure to $S$-valued points of $X$ for various $S$ defines the point functor on morphisms. 

\begin{Def}[grp-act][groups and actions]
  A \Define[group!in $\cat C$]{$\cat C$-group} is the data of $G\in\Ob\cat C$, \scth all non-empty finite products $G\times\dotsm\times G$ exist in $\cat C$, together with morphisms 
  \[
    1=1_G:*\longrightarrow G,\quad i:G\longrightarrow G,\quad m:G\times G\longrightarrow G
  \]
  called, respectively, the \Define{unit}, the \Define{inverse}, and the \Define{multiplication} of $G$, which are assumed to satisfy, for any $S\in\Ob\cat C$ and any $r,s,t\in_SG$, the group laws
  \[
    1r=r1=r,\quad rr^{-1}=1=r^{-1}r,\quad (rs)t=r(st),
  \]
  where we denote $st\defi m(s,t)$ and $s^{-1}\defi i(s)$. In particular, $*$ is in a unique fashion a $\cat C$-group, called the \Define[group!trivial]{trivial $\cat C$-group}. Given a $\cat C$-group $G$ with structural morphisms $1$, $i$, and $m$, we define the \Define[group!opposite]{opposite $\cat C$-group} $G^\circ$ to $G$, together with the morphisms $1$ and $i$, and the multiplication $m^\circ:G\times G\longrightarrow G$, where the latter is defined by $m^\circ(s,t)\defi m(t,s)$ \fa $T\in\Ob\cat C$ and $s,t\in_TG$.

  Let $X\in\Ob\cat C$ and assume that the non-empty finite products $Y_1\times\dotsm\times Y_n$ exist in $\cat C$, where $Y_j=G$ or $Y_j=X$ for any $j$. A (left) \Define[group!action (left)]{action} of a $\cat C$-group $G$ in $\cat C$, interchangeably called a (left) \Define[gspace@$G$-space (left)]{$G$-space}, consists of the data of $X$ and a morphism
  \[
    a:G\times X\longrightarrow X,
  \] 
  written $g\cdot x=a(g,x)$, for which we have
  \[
    1\cdot x=x,\quad (rs)\cdot x=r\cdot(s\cdot x)
  \]
  for any $S\in\Ob\cat C$, $x\in_SX$, and $r,s\in_SG$. Slightly abusing terminology, it is sometimes the morphism $a$ that is called an action and the space $X$ that is called a $G$-space. A $G^\circ$-space is called a \Define[gspace@$G$-space!right]{right $G$-space}. An action of $G^\circ$ is called a \Define[action!right]{right action} of $G$.
\end{Def}

\begin{Rem}[grp-dep]
  The data in the definition of a $\cat C$-group are not independent. Given $m$ and $1$ satisfying all above equations not involving $i$, there is at most one morphism $i$ with the above conditions verified. Similarly, $1$ is determined uniquely by $m$. 

  Since the Yoneda embedding preserves limits, a $\cat C$-group is the same thing as an object $G$ of $\cat C$ whose point-functor $G(-)=\Hom[_\cat C]0{-,G}$ is group-valued. Actions can be characterised similarly.
\end{Rem}

\begin{Ex}[sgrp-ex]
  Group objects and their actions are ubiquitous in mathematics. Since our main interest lies in supergeometry, we begin with two examples from this realm. 
  \begin{enumerate}[wide]
    \item\label{item:grpob-ex-vi} The general linear supergroup $\GL(m|n)$ is a complex Lie supergroup (\ie a group object in the category of complex-analytic supermanifolds). Its functor of points is given on objects $T$ by 
    \[
      \GL(m|n)(T)\defi\Set3{
        \begin{Matrix}1
          A&B\\ C&D
        \end{Matrix}
      }{
        \begin{Matrix}[0]1
          A\in\GL(m,\sh O_\ev(T)),B\in\sh O_\odd(T)^{m\times n}\\
          C\in\sh O_\odd(T)^{n\times m},D\in\GL(n,\sh O_\ev(T))
        \end{Matrix}
      }.
    \]
    Here, we let $\sh O_k(T)\defi\Gamma(\sh O_{T,k})$, $k=\ev,\odd$, $\Gamma$ denoting global sections and $\sh O_T$ the structure sheaf of $T$, with graded parts $\sh O_{T,\ev}$ and $\sh O_{T,\odd}$. The group structure is defined by the matrix unit, matrix inversion and multiplication at the level of the point functor. 

    For $X=\aff^{m|n}$, we have 
    \[
      X(T)=\Set3{
        \begin{Matrix}1
          a\\ b
        \end{Matrix}
      }{
        a\in\sh O_\ev(T)^{m\times 1},b\in\sh O_\odd(T)^{n\times 1}
      }.
    \]
    Hence, an action of $\GL(m|n)$ on $X$ is given at the level of the functor of points by the multiplication of matrices with column vectors.

    As another example, consider $X=\mathrm{Gr}_{p|q,m|n}$, the super-Grassmannian of $p|q$-planes in $m|n$-space (where $p\sle m$ and $q\sle n$). For affine $T$, the point functor takes on the form 
    \[
      X(T)=\Set1{
        Z
      }{
        Z\text{ rank $p|q$ direct summand of }\sh O(T)^{m|n}
      }.
    \]
    Again, $\GL(m|n)$ acts by left multiplication of matrices on column vectors. For general $T$ (which need not be affine), the functor of points can be computed in terms of locally direct subsheaves, compare Ref.~\cite{manin}.
    \item\label{item:grpob-ex-vii} In the category $\cat C$ of $(\knums,\Bbbk)$-supermanifolds \cite{ahw-sing}, where $\Bbbk\subseteq\knums$ and both are $\reals$ or $\cplxs$, consider the affine superspace $G\defi\aff^{0|1}$ with the odd coordinate $\tau$. Then $G(T)=\sh O_\odd(T)$, and the addition of odd superfunctions gives $G$ the structure of a supergroup. 

    Let $X$ be a manifold. The total space $\Pi TX$ of the parity reversed tangent bundle of $X$ has the underlying manifold $X$ and the structure sheaf $\sh O_{\Pi TX}=\Omega^\bullet_X$, the sheaf of $\knums$-valued differential forms, with the $\ints/2\ints$ grading induced by the $\ints$-grading. 

    The supermanifold $\Pi TX$ has the point functor 
    \[
      \Pi TX(T)\cong\Hom[_\cat C]0{T\times\aff^{0|1},X}.
    \]
    We denote elements on the left-hand side by $f$ and the corresponding elements on the right-hand side by $\tilde f$.

    We may let $x\in_TG$ act on $f\in_T\Pi TX$ by defining $x\cdot f$ \via 
    \[
      (x\cdot f)^\sim:T\times\aff^{0|1}\longrightarrow X:(t,y)\in_R(T\times\aff^{0|1})\longmapsto \tilde f(t,y+x(t))\in_RX.
    \]
    If $X$ has local coordinates $(x^a)$, then $\Pi TX$ has local coordinates $(x^a,dx^a)$. If $f\in_T\Pi TX$, then in terms of the point functor above, we have 
    \[
      f^\sharp(x^a)=j^\sharp(\tilde f^\sharp(x^a)),\quad f^\sharp(dx^a)=j^\sharp\Parens1{\tfrac\partial{\partial\tau}\tilde f^\sharp(x^a)}.
    \]
    Here, $j:T\longrightarrow T\times\aff^{0|1}$ is the unique morphism over $T$ defined by $j^\sharp(\tau)\defi0$, $\tau$ denoting the standard odd coordinate function on $\aff^{0|1}$.

    From this description, we find that the action of $G$ on $\Pi TX$ is the morphism 
    \[
      a:G\times\Pi TX\longrightarrow\Pi TX,\quad a^\sharp(\omega)=\omega+\tau d\omega.
    \]

    Expanding on this example a little, one may consider the action $\alpha$ of $(\aff^1,+)$ on $\aff^{0|1}$ given by dilation, \ie $\alpha^\sharp(\tau)=e^t\tau$. This defines a semi-direct product supergroup $G'\defi\aff^1\ltimes\aff^{0|1}$, and the action $a$ considered above may be extended to $G'$ by dilating and translating in the $\aff^{0|1}$ argument. 

    In terms of local coordinates, the thus extended action is given by 
    \[
      a^\sharp(\omega)=e^{nt}(\omega+\tau d\omega),
    \]
    for $\omega$ of degree $n$, compare \cite{hkst}*{Lemma 3.4, Proposition 3.9}.
    \item\label{item:grpob-ex-runninggag} Let $G\defi\aff^{0|1}$ with its standard additive structure and $X\defi\aff^{1|1}$. Then $G$ acts on $X$ \via $a:G\times X\longrightarrow X$, defined by 
    \[
      a\Parens0{\gamma,(y,\eta)}\defi(y+\gamma\eta,\eta)
    \]
    for all $R$ and $\gamma\in_RG$, $(y,\eta)\in_RX$. In terms of the standard coordinates $\gamma$ on $G$ and $(y,\eta)$ on $X$, we have
    \[
      a^\sharp(y)=y+\gamma\eta,\quad a^\sharp(\eta)=\eta.
    \]
    \end{enumerate}  
\end{Ex}

\begin{Ex}[grpob-ex]
  Complementing our examples from supergeometry, we give a list of examples for categorical groups and actions from different contexts. 
  \begin{enumerate}[wide]
    \item\label{item:grpob-ex-i} Let $G$ be a $\cat C$-group. Any $X\in\Ob\cat C$ can be endowed with a natural $G$-action, given by taking $a:G\times X\longrightarrow X$ to be the second projection. That is, $g\cdot x\defi x$ \fa $T\in\Ob\cat C$, $g\in_TG$, and $x\in_TX$. This action is called \Define[action!trivial]{trivial}.\index{gspace@$G$-space!trivial} 
    \item\label{item:grpob-ex-ii} Any $\cat C$-group $G$ is both a left and a right $G$-space, by the assignments 
    \[
      g\cdot x\defi gx\Or x\cdot g\defi xg,
    \]
    respectively, \fa $T\in\Ob\cat C$, $g\in_TG$, and $x\in_TX$.
    \item\label{item:grpob-ex-iii} Topological groups and Lie groups, and their actions on topological spaces and smooth manifolds, respectively, are examples of categorical groups and actions. 
    \item Group schemes and their actions on schemes are examples of categorical groups and actions as well, see \citelist{\cite{mfk}*{Definitions 0.2--3} \cite{demazure-gabriel}*{Chapitre II, \S 1.1}}.
    \item\label{item:grpob-ex-iv} A pointed (compactly generated) topological space $(W,w_0)$ is called an \textit{$H$-group}, if it is equipped with based continuous maps $\mu:W\times W \longrightarrow W$, $e:W\longrightarrow W$ with $e(W)=w_0$, and $j:W\longrightarrow W$ such that the following holds:
    \begin{gather*}
      \mu \circ (e,{\id_W})\simeq \mu \circ ({\id}_W,e)\simeq \id_W,\\
      \mu \circ (\mu\times{\id}_W)\simeq \mu \circ ({\id}_W \times \mu),\quad
      \mu \circ ({\id}_W,j)\simeq\mu \circ (j,{\id}_W)\simeq e,
    \end{gather*}
    where $\simeq$ denotes based homotopy equivalence, \cf \cite{agp}*{Section 2.7}. Given a pointed, compactly generated topological space $(X,x_0)$, its based loop space $\Omega X$ is a prime example of an $H$-group.

    In the category $\cat C$ of pointed, compactly generated topological spaces with based homotopy classes of continuous maps as morphisms, an $H$-group together with the homotopy classes of $e$, $j$, and $\mu$ is simply a $\cat C$-group. The basic theorem that the set $[X,W]_*=\Hom[_\cat C]0{X,W}$ of based homotopy classes has a group structure that is natural in the variable $X$ if and only if $W$ is an $H$-group \cite{agp}*{Theorem 2.7.6} is an instance of \thmref{Rem}{grp-dep}.

    If now $(G,1_G)=(W,w_0)$ is an $H$-group and $(X,x_0)$ a pointed topological space, then a pointed continuous map $a:G\times X\longrightarrow X$ is a group action in $\mathbf C$ if and only if $a(1_G,\cdot)$ is pointed homotopy equivalent to ${\id}_X$ and the diagram
    \[
      \begin{tikzcd}[column sep=large]
        G\times G\times X\dar[swap]{\mu\times{\id}_X}\rar{{\id}_G\times a}&G\times X\dar{a}\\
        G\times X\rar{a}&X
      \end{tikzcd}
    \]
    commutes up to a pointed homotopy.

    \item\label{item:grpob-ex-viii} In the theory of integrable systems one encounters the following situation: $(M,\omega)$ is a symplectic manifold of dimension $2n$ and $\rho:M\longrightarrow B$ is a fibration whose fibres are compact, connected Lagrangian submanifolds. Then there is a smooth fibrewise action of $T^*B$ on $M$. In the above language, $T^*B\longrightarrow B$ is a group in the category of smooth manifolds over $B$, and it acts on $X=(M\longrightarrow B)$.

    To see this latter fact, let $m\in M$, $b=\rho(m)$, and $M_b:=\rho^{-1}(b)$. The dual of the differential of $\rho$ is an injective linear map $(T_m\rho)^*:T^*_bB\longrightarrow T_m^*M$ whose image is the annihilator of $T_m(M_b)$. Since $M_b$ is Lagrangian, the musical isomorphism $\omega_m^\flat :T_m^*M\longrightarrow T_mM$ identifies this annihilator space with $T_m(M_b)$. We thus have canonical linear isomorphisms $T_b^*B\longrightarrow T_m(M_b)$ depending smoothly on $m$. Given $v\in T_b^*B$, we obtain a smooth vector field $\hat{v}$ on $M_b$. 

    It is easy to see that these vector fields extend to a commuting family of Hamiltonian vector fields on $M$, and that a linearly independent set of elements of $T_b^*B$ yields vector fields on the fibre $M_b$ that are everywhere independent. Since $M_b$ is compact, we obtain an action of the additive group of $T_b^*B$ whose isotropy is a cocompact lattice $\Lambda_b$ \cite{gs}*{Theorem 44.1}.
  \end{enumerate}
\end{Ex}

\subsection{Isotropies at generalised points}

For many applications of group actions, the notion of isotropy (or stabiliser) groups is essential. In the categorical framework, we can consider isotropy groups through $T$-valued points, by following the general philosophy of base change and specialisation: As we shall see, this allows us to consider $T$-valued points as ordinary points in the category of objects over $T$, leading to a general definition of isotropy groups. 

\begin{Cons}[grp-act-spec][base change of groups and actions]
  Let $G$ be a $\cat C$-group, $X$ a $G$-space and $T\in\Ob\cat C$. We assume that the finite products $T\times Y_1\times\dotsm\times Y_n$ exist in $\cat C$ for any choice of $Y_j=X$ or $Y_j=G$. 

  Consider the category $\cat C_T$. The morphism $\id_T:T\longrightarrow T$ is a terminal object in $\cat C_T$. Non-empty finite products in $\cat C_T$, provided they exist, are fibre products $\times_T$ over $T$ in $\cat C$. Thus, if we denote 
  \[
    G_T\defi T\times G,\quad X_T\defi T\times X,
  \]
  then 
  \[
    (Y_1)_T\times_T\dotsm\times_T(Y_n)_T=T\times Y_1\times\dotsm\times Y_n=(Y_1\times\dotsm\times Y_n)_T
  \]
  exist as finite products in $\cat C_T$. So it makes sense to define on $G_T$ and $X_T$ the structure of a $\cat C_T$-group and a $G_T$-space, respectively. The $\cat C_T$-group structure 
  \[
    1=1_{G_T}:T\longrightarrow G_T,\quad i=i_{G_T}:G_T\longrightarrow G_T,\quad m=m_{G_T}:G_T\times_TG_T\longrightarrow G_T
  \]
  on $G_T$ is defined by the equations
  \[
    1(t)\defi(t,1),\quad (t,g)^{-1}\defi(t,g^{-1}),\quad (t,g)(t,h)\defi(t,gh)
  \]
  for all $g,h\in_RG$ and $t\in_RT$, where we have written all morphisms in $\cat C$ and used the notational conventions from \thmref{Def}{grp-act}.

  Similarly, $X_T$ is a $G_T$-space \via 
  \[
    G_T\times_TX_T\longrightarrow X_T:(t,g)\cdot(t,x)\defi(t,g\cdot x)
  \]
  \fa $g\in_RG$, $x\in_RG$, and $t\in_RT$.
\end{Cons}

As we have seen, groups and actions are easily defined in the full generality of categories with terminal objects. Possibly after base change and specialisation, it will be sufficient to consider isotropy groups only through ordinary points. Their definition at the level of functors presents no difficulty. 

We will define isotropy groups at ordinary points, passing to the general case of $T$-valued points only after base change. This definition will be equivalent to the one given in Ref.~\cite{mfk}*{Definition 0.4} in the case of schemes over some base scheme. 

\begin{Def}[isotropy][isotropy group]
  Let $G$ be a $\cat C$-group and $X$ a $G$-space. We write $X_0\defi X(*)$ and call the elements of this set the \Define{ordinary points} of $X$. Let $x\in X_0$. The \Define{isotropy at $x$} (a.k.a.~the \Define{stabiliser at $x$}) is the functor $G_x:\cat C\longrightarrow\Sets$ whose object map is defined by 
  \[
    G_x(R)\defi\Set1{g\in_RG}{g\cdot x=x},
  \]
  for any $R\in\Ob\cat C$. In other words, $G_x$ is the fibre product defined by the following diagram in the category of set-valued functors on $\cat C$:
  \begin{center}
    \begin{tikzcd}  
      G_x\dar{}\rar{}&G\dar{a_x}\\
      *\rar{x}&X.
    \end{tikzcd}
  \end{center}
  Here, $a_x:G\longrightarrow X$ is the \Define{orbit morphism} defined by 
  \begin{equation}\label{eq:orbit-mor-def}
    a_x(g)\defi g\cdot x
  \end{equation}
  \fa $R\in\Ob\cat C$ and $g\in_RG$. 

  The functor $G_x$ is group-valued. Indeed, let $R\in\Ob\cat C$. By construction, an $R$-valued point $g\in G_x(R)$ is just $g\in_RG$ \scth $g\cdot x=x$. If $g,h\in G_x(R)$, then 
  \[
    (gh)\cdot x=g\cdot(h\cdot x)=g\cdot x=x,
  \]
  so $gh\in G_x(R)$. Taking this as the definition of the group law on $G_x$, we see that the canonical morphism $G_x\longrightarrow G$ preserves this operation. Since $G(R)$ is a group, so is $G_x(R)$, and this proves the assertion. In particular, if $G_x$ is representable and the finite direct products $G_x\times\dotsm\times G_x$ exist, then $G_x$ is a $\cat C$-group. 
\end{Def}

Although the above definition defines the isotropy group only for ordinary points, we may use the procedure of base change from \thmref{Cons}{grp-act-spec} to give a satisfactory definition of the isotropy of an action at a $T$-valued point, as we now proceed to explain in detail.

\begin{Cons}[iso-tpt][\protect{$T$-valued points as ordinary points}]
  Recall the natural bijection
  \begin{equation}\label{eq:comma-iso}
    \Hom[_\cat C]0{A,B}\longrightarrow\Hom[_{\cat C_T}]0{A,B_T}:f\longmapsto(p_A,f),  
  \end{equation}
  valid for any $(p_A:A\longrightarrow T)\in\Ob\cat C_T$ and any $B\in\Ob\cat C$. This allows us to consider any morphism in $\cat C$ from an object over $T$ as a morphism over $T$.

  Applying this to $A=T=*_T$, we obtain in the notation of \thmref{Def}{isotropy}
  \[
    (X_T)_0=\Hom[_{\cat C_T}]0{*_T,X_T}=\Hom[_\cat C]0{T,X}=X(T).
  \]
  Thus, \emph{we may consider any $T$-valued point $x$ of $X$ as an ordinary point of the base change $X_T\in\Ob\cat C_T$ of $X$}. This is one of the main distinguishing traits of our general point of view. 

  Let now $G$ be a $\cat C$-group, $X$ a $G$-space, and $x\in_TX$. By \thmref{Cons}{grp-act-spec}, $G_T$ is a $\cat C_T$-group and $X_T$ is a $G_T$-space. In particular, we obtain an \Define{orbit morphism} $a_x:G_T\longrightarrow X_T$ in $\cat C_T$, from Equation \eqref{eq:orbit-mor-def}. It is the composite
  \begin{center}
    \begin{tikzcd}[column sep=huge]
      T\times G\rar{({\id}_T,x)\times{\id}_G}&T\times X\times G\rar{{\id}_T\times(a\circ\sigma)}&T\times X,
    \end{tikzcd}
  \end{center}
  denoting the action of $G$ on $X$ by $a$, and by $\sigma$ the exchange of factors, \ie
  \begin{equation}\label{eq:orbmap-def}
    a_x(t,g)=\Parens1{t,g\cdot x(t)},\quad\forall t\in_RT,g\in_RG.
  \end{equation}  
\end{Cons}

The objects $T=*_T$, $G_T$, and $X_T$ in the category $\cat C_T$ are promoted to contravariant functors on $\cat C_T$. Similarly, $x$ and $a_x:G_T\longrightarrow X_T$ are promoted to morphisms of functors. We now pose the following definition. 

\begin{Def}[isotrop-def][isotropy functor]
  The \Define{isotropy functor} (a.k.a.~\Define{stabiliser functor}) $G_x\defi(G_T)_x:\cat C_T\longrightarrow\Sets$ is the fibre product defined by the diagram 
  \begin{center}
    \begin{tikzcd}  
      G_x\dar{}\rar{}&G_T\dar{a_x}\\
      T=*_T\rar{x}&X_T
    \end{tikzcd}
  \end{center}
  in the category of (contravariant) set-valued functors on $\cat C_T$. 
\end{Def}

\begin{Rem}
  This coincides with Mumford's definition \cite{mfk}*{Definition 0.4} in the case of $\cat C=\mathop{\mathrm{Sch}}_S$.  
\end{Rem}

Consider now the following diagram in the category $\cat C$:
\begin{center}
  \begin{tikzcd}  
    {}&T\times G\dar{a_x}\\
    T\rar{(\id_T,x)}&T\times X.
  \end{tikzcd}
\end{center}
Its limit in the functor category is the fibre product functor given on $R\in\Ob\cat C$ by
\begin{align*}
  \Parens1{T\times_{T\times X}(T\times G)}(R)
  &=\Set3{(t_1,t_2,g)\in_R(T\times T\times G)}{\begin{aligned}t_1&=t_2\\x(t_1)&=g\cdot x(t_2)\end{aligned}}\\
  &=\Set1{(t,g)\in_R(T\times G)}{g\cdot x(t)=x(t)}.
\end{align*}
If $R$ comes with morphisms $R\longrightarrow T$ and $R\longrightarrow T\times G$ in $\cat C$ completing the fibre product diagram above, then we may consider $R\in\Ob\cat C_T$ \via either of the $T$-projections thus obtained. The above computation then gives 
\[
  G_x(R)=\Parens1{T\times_{T\times X}(T\times G)}(R).
\]
Hence, the representability of the functor $G_x=(G_T)_x$ in $\cat C_T$ is equivalent to the existence of this fibre product in $\cat C$. 

\begin{Ex}[isotropy-runninggag]
  Recall the notation from \thmref{Ex}{sgrp-ex} \eqref{item:grpob-ex-runninggag}. We will investigate the representability of the isotropy functor for different choices of points. To that end, recall the category $\ssplfg{\knums}=\ssplfg[\varpi]{\knums,\Bbbk}$ of locally finitely generated superspaces from Section \ref{sec:groupoid-quots} below and/or Ref.~\cite{ahw-sing}. This category is finitely complete and contains $\SMan_\knums=\SMan_{\knums,\Bbbk}$ as a full subcategory. Finite limits in $\SMan_\knums$, when they exist, are finite limits in $\ssplfg{\knums}$.

  Any point $p\in X_0=X(*)$ gives rise to $p_R\in X(R)$ and we obviously have $\gamma\cdot p_R=p_R$ \fa $\gamma\in_RG$ and all $R\in\Ob\ssplfg{\knums}$. Thus, we find $G_p=G$ as functors, so $G_p$ is represented by the Lie supergroup $G$.

  By contrast, take $T=\aff^{0|1}$ with the odd coordinate $\theta$ and define $x\in_TX$ by 
  \[
    x^\sharp(y)\defi0,\quad x^\sharp(\eta)\defi\theta.  
  \]
  where we might as well take any other number for $x^\sharp(y)$. That is, for any $R\in\ssplfg{\knums}$, we have 
  \[
    x(\theta)=(0,\theta),\quad\forall\theta\in_RT.
  \]
  In this case, the isotropy functor $G_x$ evaluates on any $R\in\ssplfg{T}$ as 
  \[
    G_x(R)=\Set1{(\theta,\gamma)\in_R(T\times G)}{\gamma\theta=0}.
  \]
  Therefore, $G_x$ is represented by the superspace 
  \[
    \Spec\knums[\theta,\gamma]/(\theta\gamma)=\Parens1{*,\knums[\theta,\gamma]/(\theta\gamma)},
  \]
  where $\theta,\gamma$ are odd indeterminates. It lies over $T$ \via the morphism 
  \[
    p:G_x\longrightarrow T,\quad p^\sharp(\theta)\defi\theta.
  \]
  The group multiplication works out to be 
  \[
    m:G_x\times_T G_x\longrightarrow G_x,\quad m^\sharp(\gamma)\defi\gamma_1+\gamma_2,
  \]
  where $\gamma_i\defi p_i^\sharp(\gamma)$. Thus, $G_x$ is a group object in $\ssplfg{T}$ but not given by a Lie supergroup over $T$.
\end{Ex}

\begin{Def}[spec-pt][specialisation of a point]
  Let $\cat C$ be a category, $T_1,T_2,X$ be objects in $\cat C$. Given two points $x_1\in_{T_1}X$ and $x_2\in_{T_2}X$, we say that $x_2$ is a \Define{specialisation} of $x_1$ if for some morphism $\vphi:T_2\longrightarrow T_1$ in $\cat C$, the following diagram commutes: 
  \begin{center}
    \begin{tikzcd}
      T_2\arrow{rr}{\vphi}\arrow[swap]{dr}{x_2}&&T_1\arrow{dl}{x_1}\\
      &X.
    \end{tikzcd}
  \end{center}
\end{Def}

\begin{Prop}[isotropy-basechange]
  Let $G$ be a $\cat C$-group and $X$ a $G$-space. Let $x_1\in_{T_1}X$ and $x_2\in_{T_2}X$ \scth $x_2$ is a specialisation of $x_1$. Then there is a natural isomorphism 
  \[
    T_2\times_{T_1}G_{x_1}=G_{x_2}
  \]
  of $\Sets$-valued contravariant functors on $\cat C_{T_2}$. 

  In particular, if $G_{x_1}$ is representable in $\cat C_{T_1}$, then $G_{x_2}$ is representable in $\cat C_{T_2}$ if and only if the fibre product $T_2\times_{T_1}G_{x_1}$ exists in $\cat C$.
\end{Prop}

\begin{proof}
  By assumption, we have $x_2=x_1\circ\vphi$ for some morphism $\vphi:T_2\longrightarrow T_1$ in $\cat C$. We compute for each $R\in\Ob\cat C$ and $(t,g)\in_RG_{T_2}$ that 
  \[
    g\cdot x_2(t)=g\cdot x_1(\vphi(t)),
  \]
  so that the map $(t,g)\longmapsto(t,\vphi(t),g)$ on $R$-valued points defines a natural bijection 
  \[
    G_{x_2}(R)\longrightarrow\Parens1{T_2\times_{T_1}G_{x_1}}(R).
  \]
  This proves the assertion. 
\end{proof}

\begin{Def}[free-act][free $G$-spaces]
  Let $G$ be a $\cat C$-group and $X$ a $G$-space. Given a $T$-valued point $x\in_TX$, the $G$-space $X$ is called \Define[action!free at $x$]{free at $x$}\index{gspace@{$G$-space}!free at $x$} if $(G_T)_x$ is the trivial group in the category of $\Sets$-valued contravariant functors on $\cat C_T$. It is simply called \Define[action!free]{free}\index{gspace@{$G$-space}!free} if it is free at any $x\in_TX$, for any $T\in\Ob\cat C$.
\end{Def}

As the following corollary to \thmref{Prop}{isotropy-basechange} shows, it is equivalent to require that $X$ be free at the generic point $x=\id_X\in_XX$.

\begin{Cor}
  Let $G$ be a $\cat C$-group and $X$ a $G$-space. Assume that $X$ is free at the generic point $x=\id_X\in_XX$. Then $X$ is free. 
\end{Cor}

\subsection{Quotients and orbits}\label{subs:quot-orb}

In this subsection, we introduce basic facts and terminology relating to quotients and orbits. Although we are mainly interested in quotients by group actions, we shall need a general statement on quotients by equivalence relations for our applications (see \thmref{Prop}{eq-quot}, which is applied in the proof of \thmref{Prop}{orb-quot}). 

In order to be able to treat quotients by group actions and equivalence relations on the same footing, the language of groupoids, introduced to this context by P.~Gabriel \cite{gabriel}*{\S~1}, has proved to be convenient. Moreover, applications in forthcoming work actually rely on this generality. We briefly recall the main definitions and give a number of motivating examples before going into the applications. In what follows, we let $\cat C$ be a category with all finite products.

\begin{Def}[groupoid][groupoids]
  Let $X\in\Ob\cat C$. A \Define{$\cat C$-groupoid on $X$} is a $\Gamma\in\Ob C$, together with morphisms $s,t:\Gamma\longrightarrow X$---called \Define{source}  and \Define{target}---\scth all finite fibre products
  \[
    \Gamma^{(n)}\defi\Gamma\times_X\Gamma\times_X\dotsm\times_X\Gamma=\Gamma\times_{s,X,t}\Gamma\times_{s,X,t}\dotsm\times_{s,X,t}\Gamma  
  \]
  exist, and morphisms 
  \[
    1:X\longrightarrow\Gamma,\quad i:\Gamma\longrightarrow\Gamma,\quad m:\Gamma^{(2)}\longrightarrow\Gamma
  \]
  ---where the first and third are over $X\times X$ (where we consider $X$ as lying over $X\times X$ \via $\Delta_X$ and $\Gamma$ as lying over $X\times X$ \via $(t,s)$) and the second is over the flip $\sigma:X\times X\longrightarrow X\times X$---\scth the following diagrams commute:
  \[
    \begin{tikzcd}
      \Gamma^{(3)}\rar{m\times_X{\id}}\dar[swap]{{\id}\times_X m}&\Gamma^{(2)}\dar{m}\\
      \Gamma^{(2)}\rar[swap]{m}&\Gamma
    \end{tikzcd}
    \begin{tikzcd}
      \Gamma\arrow[mathdouble]{rd}\rar{(1\circ t)\times_X{\id}}\dar[swap]{{\id}\times_X(1\circ s)}&\Gamma^{(2)}\dar{m}\\
      \Gamma^{(2)}\rar[swap]{m}&\Gamma
    \end{tikzcd}
    \begin{tikzcd}
      \Gamma\rar{s}\dar[swap]{(i,{\id})}&X\dar{1}\\
      \Gamma^{(2)}\rar[swap]{m}&\Gamma
    \end{tikzcd}
    \begin{tikzcd}
      \Gamma\rar{({\id},i)}\dar[swap]{t}&\Gamma^{(2)}\dar{m}\\
      X\rar[swap]{1}&\Gamma.
    \end{tikzcd}
  \]
  
  A morphism $\vphi:X\longrightarrow Y$ in $\cat C$ that coequalises $s$ and $t$, \ie
  \[
    \vphi\circ s=\vphi\circ t:\Gamma\longrightarrow Y
  \]
  will be called \Define{$\Gamma$-invariant}. 

  A \Define{subgroupoid} of $\Gamma$ is a monomorphism $j:\Gamma'\longrightarrow\Gamma$ with the induced source and target morphisms, \scth $1$, $i\circ j$, and $m\circ(j\times_Xj)$ factor through $j$.
\end{Def}

\begin{Ex}[ex-groupoid]
  We will need the following three simple examples of groupoids. 
  \begin{enumerate}[wide]
    \item Let $G$ be a $\cat C$-group and $X$ be a $G$-space with action morphism $a$. Then $\Gamma\defi G\times X$ is a $\cat C$-groupoid over $X$, called the \Define{action groupoid} of $a$. Its structural morphisms are 
    \[
      s\defi p_2:\Gamma\longrightarrow X,\quad t\defi a:\Gamma\longrightarrow X,\quad 1\defi(1_G,{\id}_X):X\longrightarrow\Gamma,
    \]
    as well as the inversion $i$ and multiplication $m$ defined by 
    \[
      i(g,x)\defi(g^{-1},g\cdot x),\quad m(g_1,x,g_2)\defi(g_1g_2,x),\quad\forall g_1,g_2\in_TG,x\in_TX,
    \]
    respectively. Here, we identify $\Gamma^{(2)}=G\times X\times G$ \via the morphism induced by ${\id}_\Gamma\times p_1:\Gamma\times\Gamma\longrightarrow G\times X\times G$.
    \item Let $X\in\Ob\cat C$. Then $\Gamma\defi X\times X$ is a $\cat C$-groupoid over $X$, called the \Define{pair groupoid} of $X$. Its structural morphisms are
    \[
      s\defi p_1,t\defi p_2:\Gamma\longrightarrow X,\quad 1\defi\Delta_X:X\longrightarrow\Gamma,
    \]
    as well the inversion $i$ and multiplication $m$ defined by 
    \[
      i(x,y)\defi(y,x),\quad m(x,y,z)\defi(x,z),\quad\forall x,y,z\in_TX,
    \]
    respectively. Here, we identify $\Gamma^{(2)}=X\times X\times X$ \via the morphism induced by ${\id}_\Gamma\times p_2:\Gamma\times\Gamma\longrightarrow X\times X\times X$, 
    \item\label{item:ex-groupoid-iii} Let $X\in\Ob\cat C$. By definition, an \Define{equivalence relation} on $X$ is a subgroupoid $R$ of the pair groupoid of $X$. This definition, in the categorical context, seems to be due to P.~Gabriel \cite{gabriel}*{\S 3 e)}. Almorox \cite{almorox}*{Definition 2.1} was the first to adapt this definition to the case of supermanifolds.
  \end{enumerate}
\end{Ex}

We now recall the notion of categorical quotients \cite{mfk}*{Definition 0.5}. Although Mumford does not use the language of groupoids introduced above, the definition immediately extends to this case. 

\begin{Def}[quot-def][categorical quotients]
  Let $X\in\Ob\cat C$ and $\Gamma$ be a $\cat C$-groupoid on $X$. A morphism $\pi:X\longrightarrow Q$ is called a \Define{categorical quotient} of $X$ by $\Gamma$ if it is universal among $\Gamma$-invariant morphisms. That is, the morphism $\pi$ is $\Gamma$-invariant, and for any $\Gamma$-invariant morphism $f:X\longrightarrow Y$, where $Y\in\Ob\cat C$, there is a unique morphism $\tilde f:Q\longrightarrow Y$ \scth the following diagram commutes:
  \[
    \begin{tikzcd}
      X\arrow{rd}[swap]{f}\rar{\pi}&Q\arrow[dashed]{d}{\tilde f}\\
      &Y.
    \end{tikzcd}
  \]
  By abuse of notation, we also say that $Q$ is a categorical quotient (of $X$ by $\Gamma$).

  We say that $\pi:X\longrightarrow Q$ is a \Define{universal categorical quotient} if for all morphisms $Q'\longrightarrow Q$, the fibre products $X'\defi Q'\times_QX$ and $\Gamma'\defi (Q'\times Q')\times_{Q\times Q}\Gamma$ exist, and $\pi'\defi Q'\times_Q\pi:X'\longrightarrow Q'$ is a categorical quotient of $X'$ by $\Gamma'$. 

  We use the notation $X/\Gamma$ for categorical quotients. In case $\Gamma$ is the action groupoid for the left (respectively, right) action of a $\cat C$-group $G$, we write $G\backslash X$ (respectively, $X/G$) for the categorical quotient (if it exists).
\end{Def}

We now apply these notions to pointed spaces, to arrive at a definition of orbits. At this point, we have to depart from Mumford's definitions \cite{mfk}*{Definition 0.4}, since the notion of scheme-theoretic image does not apply to the setting of $\sh C^\infty$ differentiable supermanifolds that we are primarily interested in. 

For any category $\cat C$ with a terminal object $*$, we define the category $\cat C^*$ of \Define{pointed spaces} to be the category of objects and morphisms under $*$. We denote the objects $*\longrightarrow X$ in this category by $(X,x)$.

\begin{Def}[orb-def][categorical orbits]
  Let $G$ be a $\cat C$-group and $X$ be a $G$-space. Let $x\in_TX$, where $T\in\Ob\cat C$ is arbitrary. Assume that $G_x$ is representable in $\cat C_T$. Being a group object in that category, it is naturally pointed by the unit. 
  Since the unit acts trivially, we have a right $G_x$-action on $G_T$ in $(\cat C_T)^*$. If it exists, the categorical quotient $\pi_x:G_T\longrightarrow G_T/G_x$ in $(\cat C_T)^*$ is called the \Define{categorical orbit} of $G$ through $x$, and denoted by $\pi_x:G_T\longrightarrow G\cdot x$. If the quotient is universal categorical, then we say that the orbit is \Define{universal categorical}.

  The space $X_T$ is pointed by
  \[
    x_T\defi({\id}_T,x):T\longrightarrow X_T,
  \]
  and by definition, $G_x$ acts trivially on $x_T$, so if the categorical orbit exists, there is a unique pointed morphism $\tilde a_x:G\cdot x\longrightarrow X_T$ over $T$ \scth $\tilde a_x\circ\pi_x=a_x$. In order to avoid cluttering our terminology, we also refer to $\tilde a_x$ as the \Define{orbit morphism} of $x$. Also, by definition, the categorical orbit $G\cdot x$ is pointed in $\cat C_T$, so that it comes with a section $T\longrightarrow G\cdot x$ whose composite with $\tilde a_x$ is $x$. We call this section \Define{canonical} and will usually also denote it by $x$.
\end{Def}

We now spell out in detail what the definition given above of an orbit through a $T$-valued point is. Let $G$ be a $\cat C$-group, $X$ a $G$-space in $\cat C$, $T\in\Ob\cat C$, and $x\in_TX$. Assume that $G_x$ is representable in $\cat C_T$. As we have seen above, this means that the fibre product
\[
  G_x=T\times_{T\times X}(T\times G)
\]
exists in $\cat C$. So we have in $\cat C$ a fibre product diagram 
\[
  \begin{tikzcd}  
    G_x\dar{}\rar{}&T\times G\dar{a_x}\\
    T\rar{(\id_T,x)}&T\times X.
  \end{tikzcd}  
\]

Recall that we are working under assumption that finite products exist in $\cat C$. Then $G\cdot x$, provided it exists in $(\cat C_T)^*$, is characterised as follows: For every $G_x$-invariant morphism $f$, which fits into a commutative diagram as depicted on the left-hand side of the display below, there is a unique morphism $\tilde f$ completing the right-hand diagram commutatively:
\[
  \begin{tikzcd}[column sep=scriptsize,row sep=tiny]
    &T\arrow[mathdouble]{dd}\arrow[swap]{dl}{({\id}_T,1)}\arrow{dr}{y}&&&T\arrow[mathdouble]{dd}\arrow[swap]{dl}\arrow{dr}{y}\\
    T\times G\arrow[crossing over,near start]{rr}{f}\arrow[swap]{dr}&&Y\arrow{dl}{p_Y}
    &G\cdot x\arrow[dotted,crossing over,near start]{rr}{\exists!\tilde f}\arrow{dr}&&Y\arrow{dl}{p_Y}\\
    &T&&&T\\
  \end{tikzcd}  
\]

In other words, for any such $T$, the set of pointed morphisms $G\cdot x\longrightarrow Y$ in $\cat C_T$ is in natural bijection to the set of morphisms $f:G_T\longrightarrow Y$, which satisfy the conditions: 
\[
  \left\{
  \begin{aligned}
    f(t,1)&=y(t),\\
    p_Y(f(t,g))&=t,\\
    h\cdot x(t)&=x(t)\ \Longrightarrow\ f(t,gh)=f(t,g)
  \end{aligned}\right.
\]
\fa $R\in\Ob\cat C$, $g,h\in_RG$, and $t\in_RT$. Here, we recall that the equation $h\cdot x(t)=x(t)$ characterises the $R$-valued points $(t,h)$ of $G_x$.

\medskip
Universal categorical orbits carry a natural action.

\begin{Prop}[orbit-action]
  Let $G$ be a $\cat C$-group, and $(X,x)$ a pointed $G$-space in $\cat C$. If the $G$-orbit $G\cdot x$ exists and is universal categorical, then the morphism
  \[
    \pi_x\circ m:G\times G\longrightarrow G\cdot x
  \]
  induces an action of $G$ on $G\cdot x$. It is the unique action of $G$ on $G\cdot x$ for which $\pi_x:G\longrightarrow G\cdot x$ is $G$-equivariant. Moreover, the canonical point $x:*\longrightarrow G\cdot x$ of $G\cdot x$ is invariant under the action of $G_x$.
\end{Prop}

\begin{proof}
  By assumption, $G\cdot x$ is universal categorical, so the base change
  \[
    {\id}\times\pi_x:G\times G\longrightarrow G\times(G\cdot x)
  \]
  along the projection $G\times G\cdot x\longrightarrow G\cdot x$ is a categorical quotient in $\cat C$, for the groupoid
  \[
    \Gamma'\defi G\times\Gamma=G\times G\times G_x
  \]
  derived from $\Gamma=G\times G_x$. In particular, ${\id}\times\pi_x$ is an epimorphism. Applying the base change for a further copy of $G$, we see that so is ${\id}\times{\id}\times\pi_x$.

  Consider the multiplication $m$ of $G$. We have 
  \[
    \pi_x(m(g_1,g_2h))=\pi_x(g_1g_2h)=\pi_x(g_1g_2)=\pi_x(m(g_1,g_2))
  \]
  \fa $R\in\Ob\cat C$, $g_1,g_2\in_R G$, and $h\in_R G_x$. It follows that 
  \[
    (p_1,\pi_x\circ m):G\times G\longrightarrow G\times(G\cdot x)
  \]
  is $\Gamma'$-invariant, and hence, there is a unique morphism
  \[
    a_{G\cdot x}:G\times(G\cdot x)\longrightarrow G\cdot x
  \]
  \scth $a_{G\cdot x}\circ({\id}\times\pi_x)=\pi_x\circ m$. In particular, $\pi_x$ will be $G$-equivariant and $a_{G\cdot x}$ uniquely determined by this requirement as soon as we have established that it indeed is an action. To do so, we compute
  \begin{align*}
    a_{G\cdot x}\circ({\id}\times a_{G\cdot x})\circ({\id}\times{\id}\times\pi_x)&=a_{G\cdot x}\circ({\id}\times(\pi_x\circ m))\\
    &=\pi_x\circ m\circ({\id}\times m)\\
    &=\pi_x\circ m\circ(m\times{\id})\\
    &=a_{G\cdot x}\circ(m\times\pi_x)\\
    &=a_{G\cdot x}\circ(m\times{\id})\circ({\id}\times{\id}\times\pi_x),
  \end{align*}
  which shows that 
  \[
    a_{G\cdot x}\circ({\id}\times a_{G\cdot x})=a_{G\cdot x}\circ(m\times{\id}),
  \]
  since ${\id}\times{\id}\times\pi_x$ is an epimorphism. Similarly, one has 
  \[
    a_{G\cdot x}\circ(1\times{\id})={\id}_{G\cdot x}.
  \]
  Hence, $a_{G\cdot x}$ is an action for which $\pi_x$ is $G$-equivariant. We will denote it by $\cdot$, as for any action. 

  Finally, we verify the claim that $x$ is $G_x$-fixed. By construction, $\pi_x$ is pointed, so that $\pi_x(1)=x$. For $h\in_RG_x$, we compute, by use of the left $G$-equivariance and right $G_x$-invariance of $\pi_x$, that 
  \[
    h\cdot x=h\cdot\pi_x(1)=\pi_x(h\cdot 1)=\pi_x(h)=\pi_x(1)=x.
  \]
  This completes the proof of the proposition. 
\end{proof}

\begin{Ex}[orb-ex][examples of orbits]
  Returning to the groups and actions from \thmref{Ex}{grpob-ex}, we explain the notion of isotropy groups and orbits in these cases. In items \eqref{item:orb-ex-i} and \eqref{item:orb-ex-ii} below, let $\cat C$ denote a category \scth all finite products exist.
  \begin{enumerate}[wide]
    \item\label{item:orb-ex-i} Let $G$ be a $\cat C$-group acting trivially on $X\in \Ob\cat C$. Then \fa $x:T\longrightarrow X$ and $R\in\Ob\cat C$, we have $G_x(R)=G_T(R)$. Thus, the isotropy functor $G_x$ is represented by $G_T=T\times G$. Here, the morphism $\pi_x=p_1:G_T\to T$ is a universal categorical orbit, as can be seen as follows: $\pi_x$ is invariant with respect to the action groupoid $\Gamma$ coming from the right $G_T$-action on $G_T$. Given any $\Gamma$-invariant morphism $f:G_T\to Y$ with $Y$ over $T$, it uniquely factors over $\pi_x$ to $\tilde{f}=f\circ({\id}_T\times 1_G)$. 

    Furthermore, given $Q\to T$, the fibre product $Q\times_TG_T=G_Q$ exists. Moreover, $Q\times_T\pi_x={\id}_Q\times p_1:G_Q\to Q=Q\times_TT$ is a categorical quotient by the above, since $(Q\times Q)\times_{T\times T}\Gamma$ is the action groupoid for the right $G_Q$-action on $G_Q$.
    \item\label{item:orb-ex-ii} Assume given a $\cat C$-group $G$, viewed as a left $G$-space via left multiplication. For $T\in\Ob\cat C$ and $x\in_TG$, we have 
    \[
      G_x(R)=\Set1{(t,g)\in_RG_T}{g\cdot x(t)=x(t)}
      =\Set1{(t,1_G(t))}{t\in_RT}\cong T(R).
    \]
    Thus, $G_x$ is represented by $T$. Defining $\pi_x$ by ${\id}_{G_T}:G_T \longrightarrow Q\defi G_T$, we obtain for any $Y$ and any $G_T$-invariant $f:G_T\to Y$ a unique factorisation $\tilde{f}\defi f$. Thus, $\pi_x:G_T\longrightarrow Q$ is the categorical quotient of $G_T$ with respect to the $G_x$-action. In other words, it is the categorical orbit of $G$ through $x$.

    Furthermore, given $Q'\in\Ob\cat C$ and $Q'\longrightarrow Q$, we have $Q'\times_QG_T=Q'$ and $Q'\times_Q\Gamma=Q'\times_QG_T=Q'$. The projections ${\id}_Q'\times_Qs$ and ${\id}_Q'\times_Qt$ are the identity of $Q'$, so that $Q'\times_Q\pi_x={\id}_{Q'}$ is the categorical quotient of $Q'$ (the space) by $Q'$ (the groupoid). It follows that $\pi_x:G_T\longrightarrow G_T$ is a universal categorical orbit.
    \item Let a continuous or smooth action $a:G\times X\longrightarrow X$, respectively, of a topological group or Lie group on a topological space or a smooth manifold be given. The isotropies at $x\in X_0=X(*)=X$ are represented by the obvious set-theoretic isotropy groups, endowed with the subspace topology coming from the inclusion into $G$. Since these isotropies are closed, they are notably Lie subgroups in the smooth case. 

    Both in the topological and the smooth case, a categorical orbit through such an $x$ is represented by the set of right cosets with respect to the isotropy group $G_x$, with its canonical structure of topological space or smooth manifold, respectively. For the rest of this example, let us focus on the topological case. 

    Then we can consider arbitrary continuous maps $x:T\longrightarrow X$, defined on some topological space $T$ and observe that
    \[
      G_x=\Set1{(t,g)\in G_T}{g\in G_{x(t)}}
    \]
    with the subspace topology from $T\times G$. We may define an equivalence relation $\sim$ on $G_T$ by 
    \[
      (t,g)\sim(t',g')\ :\Longleftrightarrow\ t=t',g\cdot x(t)=g'\cdot x(t).
    \]
    The quotient space $Q\defi X/\sim$ with the canonical map $\pi_x:G_T\longrightarrow Q$ satisfies the universal property of the categorical orbit of $G$ through $x$. 

    If $\pi_x$ is an open map, then it is already an universal categorical orbit. Indeed, in this case, for any $Q'\longrightarrow Q$, the projection $p_1:Q'\times_QG_T\longrightarrow Q'$ is open and in particular a quotient map. The map $\pi_x$ is open in case $T=*$, which is the situation studied classically. In general, however, this fails to be true, as one may see in the following example: Let $G\defi(\reals,+)$, $T\defi\reals$, and $X\defi\reals^2$. Define the action by 
    \[
      g\cdot(t,s)\defi(t,tg+s)
    \]
    and set $x:T\longrightarrow X,x(t)\defi(t,0)$. Then $G_x=(0\times\reals)\cup(\reals^\times\times0)$ and the projection $G_T\longrightarrow G_T/G_x$ is not open, as the saturation of an open set $U\subseteq G_T$ containing $(0,0)$ is $(\reals\times0)\cup U$, which is open only if $\reals\times 0$ is already contained in $U$.

    The smooth case is more subtle, since in general, the isotropy $G_x$ might not exist as a smooth manifold over $T$. In Section \ref{sec:super-quot}, we study  these questions for the category of supermanifolds. \emph{A fortiori}, these apply to ordinary manifolds.
    \item The existence question for isotropies and orbits in the homotopy category of pointed topological spaces leads immediately to subtle questions concerning homotopy pullbacks and homotopy orbits. We do not dwell on these matters here.
    \item\label{item:orb-ex-viii} From the description of the action of $T^*B$ on $M$ in \thmref{Ex}{grpob-ex} \eqref{item:grpob-ex-viii}, it follows immediately that for any $b\in B$, the action of $T^*_bB$ on the fibre $M_b$ is transitive and the orbits are $n$-dimensional real tori. Furthermore, the isotropy is a cocompact lattice $\Lambda_b$ in $T^*_bB$, depending smoothly on $b$, \cf \cite{gs}*{Theorem 44.1}. The union of the $\Lambda_b$ is the total space a smooth $\ints^n$-principal subbundle $\Lambda$ of $T^*B\longrightarrow B$. 

    The traditional description underlines the ensuing action-angle coordinates: Action for the base directions of $B$, angle for the fibre directions (compare the detailed analysis of Duistermaat \cite{d-glo}). In the terminology introduced above, we find that the isotropy of the generic point $x={\id}_X:T=X\to X$ is the subgroup $G_x=M\times_B\Lambda$ of $G_T=M\times_BT^*B$. 

    By our results below (\thmref{Th}{trans-iso} and \thmref{Th}{orbit}), the orbit 
    \[
      G\cdot x=G_T/G_x=(M\times_BT^*B)/(M\times_B\Lambda)  
    \]
    exists as a universal categorical quotient in the category of manifolds over $M$. Moreover (\loccit), it coincides with the image of the orbit morphism $a_x$, which is a surjective submersion. Hence, we have $G\cdot x\cong M\times_BM$ as manifolds over $M$.
  \end{enumerate}
\end{Ex}

\section{Groupoid quotients of superspaces}\label{sec:groupoid-quots}

We now apply the general setup of Section \ref{sec:cats} to the categories of locally finitely generated superspaces and of relative supermanifolds constructed in Ref.~\cite{ahw-sing}. We will start by recalling some basic definitions, referring to this paper for more details.

We fix a field $\knums\in\{\reals,\cplxs\}$. The category $\SSp_\knums$ has as objects pairs $X=(X_0,\sh O_X)$ where $X_0$ is a topological space and $\sh O_X$ is a sheaf of $\knums$-superalgebras with local stalks. Such objects are called \Define{$\knums$-superspaces}. Morphisms $\vphi:X\longrightarrow Y$ are again pairs $(\vphi_0,\vphi^\sharp)$ where this time, $\vphi_0:X_0\longrightarrow Y_0$ is a continuous map and $\vphi^\sharp:\sh O_Y\longrightarrow(\vphi_0)_*\sh O_X$ is a local morphism of $\knums$-superalgebra sheaves. 

If $S$ is a fixed $\knums$-superspace, the category of objects and morphisms in $\SSp_\knums$ over $S$ will be denoted by $\SSp_S$. Objects are denoted by $X/S$ and morphisms by $\vphi:X/S\longrightarrow Y/S$. 

Now we fix a subfield $\Bbbk$ of $\knums$ containing $\reals$ and a `differentiability' class $\varpi\in\{\infty,\omega\}$. Here, $\infty$ means `smooth' and $\omega$ means `analytic' (over $\Bbbk$). We consider model spaces adapted to these data. Namely, let a finite-dimensional super-vector space $V=V_\ev\oplus V_\odd$ over $\Bbbk$ be given, together with a compatible $\knums$-structure on $V_\odd$. Then we may consider on the topological space $V_\ev$ the sheaf $\sh C_{V_\ev}^\varpi$ of $\knums$-valued functions of differentiability class $\varpi$. We set
\[
  \aff(V)\defi\Parens1{V_\ev,\sh C^\varpi_{V_\ev}\otimes_\knums\textstyle\bigwedge(V_\odd)^*}
\]
and call this the \Define{affine superspace} associated with $V$. It depends on the data of $(\knums,\Bbbk,\varpi)$, but we will usually omit them from the notation.

By definition, a \Define{supermanifold over $(\knums,\Bbbk)$ of class $\sh C^\varpi$} is a Hausdorff $\knums$-superspace $X$ which admits a cover by open subspaces which are isomorphic to open subspaces of affine superspaces. We will usually just say that $X$ is a \Define{supermanifold}. The full subcategory of $\SSp_\knums$ comprised of these objects with be denoted by $\SMan_\knums$.

In the literature, the case $\knums=\Bbbk=\reals$ corresponds to (smooth or real-analytic) real supermanifolds \citelist{\cite{leites}*{3.1.2} \cite{ccf}*{Definitions 3.2.1 and 4.2.1}}, and the case $\knums=\Bbbk=\cplxs$ corresponds to (holomorphic) supermanifolds \citelist{\cite{manin}*{Chapter 4, \S 1, Definition 1} \cite{ccf}*{Definition 4.8.1}}. In the case of $\knums=\cplxs$ and $\Bbbk=\reals$, supermanifolds are also known as `\emph{cs} manifolds' \cite{deligne-morgan}*{\S 4.8}. We take this opportunity to replace this unfortunate terminology with a hopefully less confusing one.

In Ref.~\cite{ahw-sing}, we construct a full subcategory $\ssplfg{\knums}=\ssplfg[\varpi]{\knums,\Bbbk}$ of $\SSp_\knums$ that admits finite fibre products and contains $\SMan_\knums$ as a subcategory closed under finite products. Here, `lfg' stands for `locally finitely generated'. For any $S\in\Ob\ssplfg{\knums}$, the category of objects and morphisms over $S$ in $\ssplfg{\knums}$ will be denoted by $\ssplfg{S}$. Given any super-vector space $V$ as above, we define $\aff_S(V)\defi S\times\aff(V)$ (where the product is taken in $\ssplfg{\knums}$). Using these as model spaces, we arrive at a definition of \Define{supermanifolds over $S$}, compare \opcit{} We denote the corresponding full subcategory of $\ssplfg{S}$ by $\SMan_S$. Note that this now makes sense for a wide class of \Define{singular} base spaces $S$ and, moreover, that, contrary to the setting of schemes, it would not be appropriate to instead consider products in $\SSp_S$, as already the embedding $\SMan_\knums\longrightarrow\SSp_\knums$ does not preserve products. For this reason, the use of the intermediate category $\ssplfg{\knums}$ is essential.

\subsection{Geometric versus categorical quotients}

In what follows, fix $S\in\ssplfg{\knums}$, and let $\cat C$ be a full subcategory of $\ssplfg{S}$ admitting finite products. Particular cases are $\cat C=\ssplfg{S}$ and $\cat C=\SMan_S$, by \cite{ahw-sing}*{Corollaries 5.27, 5.42}. Furthermore, let $X\in\Ob\cat C$ and $\Gamma$ be a groupoid over $X$ in $\cat C$.

\begin{Prop}[wgeom-cat]
  The coequaliser $\pi:X\longrightarrow Q$ of $s,t:\Gamma\longrightarrow X$ exists in $\SSp_S$ and is regular in the sense of \cite{ahw-sing}*{Definition 4.12}. If $Q\in\Ob\cat C$, then $Q$ is the categorical quotient of $X$ by $\Gamma$.
\end{Prop}

\begin{proof}
  The existence and regularity of $Q$ is immediate from \cite{ahw-sing}*{Propositions 2.17, 5.5}. By definition, the morphism $\pi:X\longrightarrow Q$ is a coequaliser in $\SSp_S$. But since $\cat C$ is a full subcategory of $\SSp_S$, $\ssplfg{S}$ being full in the latter, $Q$ is the coequaliser of $s,t$ in $\cat C$ if $Q\in\Ob\cat C$, and thus has the properties required by \thmref{Def}{quot-def}.
\end{proof}

\begin{Rem}[colimit-explicit]
  We can describe the colimit $Q$ of $s,t:\Gamma\longrightarrow X$ explicitly. Indeed, by \cite{ahw-sing}*{Remark 2.18}, $\sh O_Q$ is the equaliser in the category $\Sh0{Q_0}$ of sheaves on $Q_0$, defined by the diagram
  \[
    \begin{tikzcd}
      \sh O_Q\rar{\pi^\sharp}&
      \pi_{0*}\sh O_X\arrow[shift up=0.5ex]{r}{s^\sharp}\arrow[shift up=-0.5ex]{r}[swap]{t^\sharp}&
      (\pi_0\circ s_0)_*\sh O_\Gamma.
    \end{tikzcd}
  \]
  Moreover, since the embedding of $\SSp_S$ in $\SSp$ preserves colimits, one may see easily that $Q_0$ is the coequaliser of $s_0,t_0:\Gamma_0\longrightarrow X_0$, \ie the topological quotient space of $X_0$ by the equivalence relation generated by $s_0(\gamma)\sim t_0(\gamma)$.
\end{Rem}

\begin{Ex}[orbit-runninggag]
  Recall the action from \thmref{Ex}{sgrp-ex} \eqref{item:grpob-ex-runninggag} and the $T$-valued point $x$ from \thmref{Ex}{isotropy-runninggag}. Recall that the isotropy supergroup $G_x$ is in this case representable by the group object
  \[
    G_x=\Spec\knums[\theta,\gamma]/(\theta\gamma),\quad p^\sharp(\theta)=\theta,\quad m^\sharp(\gamma)=\gamma_1+\gamma_2,\quad 1^\sharp(\gamma)=0
  \]
  in $\ssplfg{T}$, where $\theta,\gamma$ are odd indeterminates. In particular, it lies in $(\ssplfg{T})^*$.

  Let $\eps$ be an even indeterminate and define
  \[
    Q\defi\Spec\knums[\eps|\theta]/(\eps^2,\theta\eps).
  \]
  We then have morphisms 
  \[
    p_Q:Q\longrightarrow T,\quad p_Q^\sharp(\theta)\defi\theta,\quad q:T\longrightarrow Q,\quad q^\sharp(\eps)\defi0,\quad q^\sharp(\theta)\defi\theta.
  \]
  The morphism 
  \[
    \pi_x:G_T\longrightarrow Q,\quad\pi_x^\sharp(\theta)\defi\theta,\quad\pi_x^\sharp(\eps)\defi\theta\gamma
  \]
  is in the category $(\ssplfg{T})^*$. We claim that $\pi_x:G_T\longrightarrow Q$ is the categorical orbit of $G$ through $x$.

  To establish this claim, let $b:G_T\times_TG_x\longrightarrow G_T$ denote the action by right multiplication of the isotropy, \ie $b^\sharp(\gamma)=\gamma_1+\gamma_2$. We compute 
  \[
    (\pi_x\circ b)^\sharp(\eps)=b^\sharp(\theta\gamma)=\theta(\gamma_1+\gamma_2)=\theta\gamma_1=p_1^\sharp(\theta\gamma)=(\pi_x\circ p_1)^\sharp(\eps)
  \]
  so $\pi_x$ is indeed $G_x$-invariant. If $f$ is a function on $G_T$, then 
  \[
    f=f_0+f_\theta\theta+f_\gamma\gamma+f_{\theta\gamma}\theta\gamma
  \]
  where $f_\alpha\in\knums$ for $\alpha=0,\theta,\gamma,\theta\gamma$. Then 
  \[
    b^\sharp(f)-p_1^\sharp(f)=f_\gamma\gamma_2,
  \]
  so $f$ is $G_x$-invariant if and only $f_\gamma=0$. In this case, 
  \[
    f=\pi_x(\tilde f),\quad\tilde f=f_0+f_\theta\theta+f_{\theta\gamma}\eps,
  \]
  and $\tilde f$ is unique with this property. It is easy to conclude that $\pi_x:G_T\longrightarrow Q$ is the categorical quotient of $G_T$ by $G_x$, and thus the claim follows. Notice that $G\cdot x=Q$ is not a supermanifold over $T$.
\end{Ex}

\begin{Def}[wgeom-quot][weakly geometric quotients]
  The coequaliser $\pi:X\longrightarrow Q$ of $s,t:\Gamma\longrightarrow X$ is called a \Define{weakly geometric quotient} of $X$ by $\Gamma$ if $Q\in\Ob\cat C$.  We say that $\pi:X\longrightarrow Q$ is a \Define{universal weakly geometric quotient} if for all morphisms $Q'\longrightarrow Q$, the fibre products $X'\defi Q'\times_QX$ and $\Gamma'\defi (Q'\times Q')\times_{Q\times Q}\Gamma$ exist in $\cat C$, and $\pi'\defi Q'\times_Q\pi:X'\longrightarrow Q'$ is the weakly geometric quotient of $X'$ by $\Gamma'$. 
\end{Def}

\begin{Rem}
  The terminology is justified as follows: If $G$ is a group scheme acting on a scheme $X$, then a morphism $\pi:X\longrightarrow Q$ is called a \Define{geometric quotient} of $X$ by $G$ if it is the coequaliser of $p_2,a:G\times X\longrightarrow X$ in the category of locally ringed spaces, and in addition, the scheme-theoretic image of $(p_2,a):G\times X\longrightarrow X\times X$ is $X\times_QX$, see \cite{mfk}*{Definition 0.6}. 
\end{Rem}

In terms of the above terminology, we may rephrase \thmref{Prop}{wgeom-cat} as follows. The result is a generalisation of \cite{mfk}*{Proposition 0.1}.

\begin{Cor}[wgeom-quot-is-cat]
  Let the (universal) weakly geometric quotient $Q$ of $X$ by $\Gamma$ exist in $\cat C$. Then $Q$ is the (universal) categorical quotient of $X$ by $\Gamma$ in $\cat C$.
\end{Cor}

\section{Existence of superorbits}\label{sec:super-quot}

In this section, we will derive general sufficient conditions for the existence of isotropies at and orbits through generalised points in the category $\SMan_S$ of supermanifolds over $S$, where in what follows, $S$ will denote some object of $\ssplfg{\knums}$.

The material is organised as follows: In Subsection \ref{subs:const-rank}, we discuss at length the notion of morphisms of constant rank basic for our considerations. In particular, we characterise precisely when the orbit morphism of a generalised point is locally of constant rank. Subsequently, in Subsection \ref{subs:isotropy}, we study the isotropy of a supergroup action at a generalised point. This leads, in Subsection \ref{subs:orbits}, to a characterisation of the existence of orbits through generalised points. 

\subsection{\texorpdfstring{Tangent sheaves of supermanifolds over $S$}{Tangent sheaves of supermanifolds over S}}

We briefly collect some definitions and facts concerning tangent sheaves. These are totally classical if $S$ is itself a supermanifold. 

\begin{Def}[tan-sh][tangent sheaf]
  Let $p_X:X\longrightarrow S$ and $p_Y:Y\longrightarrow S$ be superspaces over $S$ and $\vphi:X/S\longrightarrow Y/S$ a morphism over $S$. Let $U\subseteq X_0$ be open. An $\smash{p_{X,0}^{-1}\sh O_S}$-linear sheaf map 
  \[
    v:\vphi_0^{-1}\sh O_Y|_U\longrightarrow\sh O_X|_U
  \]
  will be called a \Define{vector field along $\vphi$ over $S$} (defined on $U$) if $v=v_\ev+v_\odd$ where
  \[
    v_i(fg)=v_i(f)\vphi^\sharp(g)+(-1)^{i\Abs0f}\vphi^\sharp(f)v_i(g)
  \]
  for all $i=\ev,\odd$ and all homogeneous local sections $f,g$ of $p_{X,0}^{-1}\sh O_Y|_U$.

  The sheaf on $X_0$ whose local sections over $U$ are the vector fields along $\vphi$ over $S$ defined on $U$ will be denoted by $\sh T_{X/S\to Y/S}$ or $\sh T_{\vphi:X/S\to Y/S}$ if we wish to emphasize $\vphi$. It is an $\sh O_X$-module, and will be called the \Define{tangent sheaf along $\vphi$ over $S$}. In particular, we define $\sh T_{X/S}\defi\sh T_{{\id}_X:X/S\to X/S}$ and $\sh T_X\defi\sh T_{X/*}$, the \Define{tangent sheaf of $X$ over $S$} and the \Define{tangent sheaf of $X$}, respectively.
\end{Def}

Let $\tau$ be an even and $\theta$ an odd indeterminate. Whenever $X$ is a $\knums$-superspace, we define 
\[
  X[\tau|\theta]\defi\Parens1{X_0,\sh O_X[\tau|\theta]/(\tau^2,\tau\theta)}.
\]
There is a natural morphism $(\cdot)|_{\tau=\theta=0}:X\longrightarrow X[\tau|\theta]$ whose underlying map is the identity and whose pullback map sends $\tau$ and $\theta$ to zero.

\begin{Lem}[der-mor][superderivations and super-dual numbers]
  Let $X/S$ and $Y/S$ be superspaces over $S$ and $\vphi:X/S\longrightarrow Y/S$ be a morphism over $S$. There is a natural bijection
  \[
    \Set1{\phi\in\Hom[_S]1{X[\tau|\theta],Y}}{\phi|_{\tau=\theta=0}=\vphi}\longrightarrow\Gamma(\sh T_{X/S\longrightarrow Y/S}):\phi\longmapsto v
  \]
  given by the equation
  \begin{equation}\label{eq:der-mor}
    \phi^\sharp(f)\equiv\vphi^\sharp(f)+\tau v_\ev(f)+\theta v_\odd(f)\pod{\tau^2,\tau\theta}
  \end{equation}
  for all local sections $f$ of $\sh O_Y$. Symbolically, we write
  \[
    v_\ev(f)=\frac{\partial\phi^\sharp(f)}{\partial\tau}\nd v_\odd(f)=\frac{\partial\phi^\sharp(f)}{\partial\theta}.
  \]
\end{Lem}

\begin{proof}
  Since $X[\tau|\theta]$ is a thickening of $X$ \cite{ahw-sing}*{Definition 2.10}, the underlying map of $\phi$ is fixed by $\phi_0=\vphi_0$. The assertion follows easily.
\end{proof}

\begin{Def}[inf-flow][infinitesimal flow]
  Let $v\in\Gamma(\sh T_{X/S\longrightarrow Y/S})$. The unique morphism $\phi^v\in\Hom[_S]0{X[\tau|\theta],Y}$, \scth $\smash{\phi^v|_{\tau=\theta=0}=\vphi}$, associated with $v$ \via \thmref{Lem}{der-mor}, is called the \Define{infinitesimal flow} of $v$.
\end{Def}

The infinitesimal flow construction allows us to introduce for each fibre coordinate system a family of fibre coordinate vector fields. 

\begin{Cons}[fib-der][fibre coordinate vector fields]
  Let $S\in\ssplfg{\knums}$ and $X/S$ be in $\SMan_S$ with a global fibre coordinate system $x=(x^a)$.
  
  By \cite{ahw-sing}*{Propositions 5.18, 4.19, Corollary 5.22}, we have $X[\theta|\tau]\in\Ob\ssplfg{\knums}$, and there are unique morphisms $\phi^a\in\Hom[_S]0{X[\tau|\theta],X}$ \scth
  \[
    \phi^{a\sharp}(x^b)=\begin{cases}x^b+\tau\delta_{ab}&\text{for }\Abs0{x^a}=\ev,\\ x^b+\theta\delta_{ab}&\text{for }\Abs0{x^a}=\odd.\end{cases}
  \]
  Evidently, we have $(\phi^a|_{\tau=\theta=0})^\sharp(x^b)=x^b$, and hence $\phi^a|_{\tau=\theta=0}={\id}_X$.
  
  On account of \thmref{Lem}{der-mor}, the morphisms $\phi^a$ are the infinitesimal flows of unique vector fields  over $S$, denoted by $\smash{\frac\partial{\partial x^a}}\in\Gamma(\sh T_{X/S})$. We call these \Define{fibre coordinate vector fields} and simply \Define{coordinate vector fields} in case $S=*$.

  Observe that the meaning of each individual $\smash{\frac\partial{\partial x^a}}$ depends on the entire fibre coordinate system $(x^b)$, and not only on the coordinate $x^a$.
\end{Cons}

As we shall presently see, the coordinate vector fields give systems of generators for the relative tangent bundle. 

\begin{Prop}[tan-coord][coordinate expression of vector fields]
  Let $S$ be in $\ssplfg{\knums}$, $X/S$ be in $\ssplfg{S}$, $Y/S$ be in $\SMan_S$, and $\vphi:X/S\longrightarrow Y/S$ be a morphism over $S$. Let $(y^a)$ be a local fibre coordinate system on an open subset $V\subseteq Y_0$. Let $U\subseteq X_0$ be an open subset, \scth $\vphi_0(U)\subseteq V$, and $v\in\sh T_{X/S\longrightarrow Y/S}(U)$. Then
  \begin{equation}\label{eq:tan-coord}
    v=\sum\nolimits_av(y^a)\,\vphi^\sharp\circ\frac\partial{\partial y^a}.
  \end{equation}
  In particular, we have 
  \[
    \sh T_{X/S\to Y/S}=\vphi^*(\sh T_{Y/S})\defi\sh O_X\otimes_{\vphi_0^{-1}\sh O_Y}\vphi_0^{-1}\sh T_{Y/S},   
  \]
  and this $\sh O_X$-module is locally free, of rank $\rk_x\sh T_{X/S\to Y/S}=\dim_{S,\vphi_0(x)}Y$ for $x\in X_0$.
\end{Prop}

\begin{proof}
  We may assume that $(y^a)$ is a globally defined fibre coordinate system. Define the vector field $v'\in\vphi^*(\sh T_{Y/S})(U)\subseteq\sh T_{X/S\to Y/S}(U)$ by
  \[
    v'\defi\sum\nolimits_av(y^a)\,\vphi^\sharp\circ\frac\partial{\partial y^a}.
  \]  
  Let $\phi$ and $\phi'$ be the infinitesimal flows of $v$ and $v'$, respectively. For any index $a$, we have $v'(y^a)=v(y^a)$, and hence $\phi^\sharp(y^a)=\phi'^\sharp(y^a)$. This implies that $\phi=\phi'$, by reason of \cite{ahw-sing}*{Propositions 5.18, 4.19, Corollary 5.22}. Hence, we have $v'=v$.

  In particular, the vector fields $\smash{\vphi^\sharp\circ\frac\partial{\partial y^a}}$ form a local basis of sections of $\sh T_{X/S\longrightarrow Y/S}$, and this readily implies the remaining assertions.
\end{proof}

\begin{Cor}[tan-free][local freeness of $\sh T_{X/S}$]
  Let $S\in\ssplfg{\knums}$ and $X/S\in\SMan_S$. Then $\sh T_{X/S}$ is locally free, with $\rk_x\sh T_{X/S}=\dim_{S,x}X$, for $x\in X_0.$
\end{Cor}

A special case of the above concerns the relative tangent spaces. 

\begin{Def}[tangent-def][tangent space]
  Let $p=p_X:X\longrightarrow S$ be a superspace over $S$. For any point $x\in X_0$ we let $\ger m_{X,x}$ be the maximal ideal of $\sh O_{X,x}$ and $\vkappa(x)\defi\sh O_{X,x}/\ger m_{X,x}$. Setting $s\defi p_{X,0}(x)$, we define 
  \[
    T_x(X/S)\defi\GDer[_{\sh O_{S,s}}]0{\sh O_{X,x},\vkappa(x)},
  \]
  the $\ints$-span of all homogeneous $v\in\GHom[_{\sh O_{S,s}}]0{\sh O_{X,x},\vkappa(x)}$ \scth
  \begin{equation}\label{eq:tanvector-def}
    v(fg)=v(f)g(x)+(-1)^{\Abs0f\Abs0v}f(x)v(g).   
  \end{equation}
  This is naturally a super-vector space over $\vkappa(x)$, called the \Define{tangent space at $x$ over $S$}. For $S=*$, we also write $T_xX$. The elements are called \Define{tangent vectors} (over $S$).

  As is immediate from the definitions, the tangent space coincides with the tangent sheaf over $S$ along the morphism $(*,\vkappa(x))\longrightarrow X$.
\end{Def}

\begin{Cor}[tansp-dim][dimension of $T_{S,x}X$]
  Let $S\in\ssplfg{\knums}$, $X/S$ be a supermanifold over $S$, and $x\in X_0$. Then $\dim_\knums T_{S,x}X=\dim_{S,x}X$.
\end{Cor}

\begin{Def}[tan-mor][tangent morphism]
  Let $\vphi:X/S\longrightarrow Y/S$ be a morphism of superspaces over $S$. We define the \Define{tangent morphism}
  \[
    \sh T_{\vphi/S}:\sh T_{X/S}\longrightarrow\sh T_{X/S\to Y/S}
  \]
  by setting 
  \[
    \sh T_{\vphi/S}(v)\defi v\circ\vphi^\sharp
  \]
  for any locally defined vector field $v$ over $S$. In view of \thmref{Prop}{tan-coord}, if $Y$ is in $\SMan_S$, then the range of $\smash{\sh T_{\vphi/S}}$ is in $\smash{\vphi^*(\sh T_{Y/S})}$.

  Similarly, we obtain for any $x\in X_0$ a \Define{tangent map}
  \[
    T_x(\vphi/S):T_x(X/S)\longrightarrow T_{\vphi_0(x)}(Y/S)
  \]
  by setting 
  \[
    T_x(\vphi/S)(v)\defi v\circ\vphi^\sharp_x
  \]
  for any tangent vector $v$ over $S$.
\end{Def}

\subsection{Morphisms of constant rank}\label{subs:const-rank}

In order to handle supergroup orbits through $T$-valued points, we will need to understand morphisms of locally constant rank in the setting of relative supermanifolds. Already for ordinary supermanifolds, the notion is somewhat different from the standard one used for manifolds. The correct definition was first given in \cite{leites}*{2.3.8}. 

For our present purposes, it is useful to state this in a more invariant form. We need the following definitions and facts, which are more or less standard.  

\begin{Def}[qcoh][Conditions on module sheaves]
  Let $\sh E$ be a sheaf (of $\ints$-modules) and $I=(I_\ev,I_\odd)$ a graded set, \ie~a pair of sets. We write $\sh E^{(I)}$ for the direct sum $\bigoplus_{I_\ev}\sh E\oplus\bigoplus_{I_\odd}\sh E$ with its natural $\ints/2\ints$-grading. 

  Let $X$ be a superspace and $\sh F$ an $\sh O_X$-module (understood to be graded). We say that $\sh E$ is \Define{locally generated by sections} if any $x\in X_0$ admits an open neighbourhood $U\subseteq X_0$ and a surjective morphism of $\sh O_X|_U$-modules $(\sh O_X|_U)^{(I)}\longrightarrow\sh E|_U$ \fs $I$. If $I$ can be chosen to be finite for any $x$, we say that $\sh E$ is of \Define{finite type}. 
\end{Def}

\begin{Prop}[free-vals]
  Let $X$ be a superspace and $\vphi:\sh E\longrightarrow\sh F$ a morphism of $\sh O_X$-modules, with $\sh E$ of finite type and $\sh F$ finite locally free. For $x\in X_0$, we define
  \[
    \sh E(x)\defi\sh E_x/\ger m_{X,x}\sh E_x.
  \]
  For every $x\in X_0$, the following are equivalent:
  \begin{enumerate}[wide]
    \item\label{item:free-vals-i} The $\vkappa(x)$-linear map defined below is injective: 
    \[
      \vphi(x)\defi\vphi_x\otimes_{\sh O_{X,x}}{\id}_{\vkappa(x)}:\sh E(x)\longrightarrow\sh F(x).
    \]
    \item\label{item:free-vals-ii} For some open neighbourhood $U\subseteq X_0$ of $x$, the morphism $\vphi|_U$ is injective and the $\sh O_X|_U$-module $(\sh F/\im\vphi)|_U$ is locally free. 
    \item\label{item:free-vals-iii} For some open neighbourhood $U\subseteq X_0$ of $x$, $\vphi|_U$ admits a left inverse. 
    \item\label{item:free-vals-iv} There exist an open neighbourhood $U\subseteq X_0$ of $x$ and homogeneous bases of sections for $\sh E|_U$ and $\sh F|_U$, \scth the matrix $M_\vphi$ of $\vphi$ is
    \[
      M_\vphi=
      \begin{Matrix}1
        A&0\\ 0&D
      \end{Matrix},\quad
      A=
      \begin{Matrix}1
        1&0\\ 0&0
      \end{Matrix},\quad
      D=
      \begin{Matrix}1
        1&0\\ 0&0
      \end{Matrix}.
    \]  
  \end{enumerate} 
  The set of all those $x\in X_0$ where this holds is open. Moreover, in this case, $\sh E$ is locally free on an open neighbourhood of $x$.
\end{Prop}

\begin{proof}
  The equivalence of \eqref{item:free-vals-i}--\eqref{item:free-vals-iii} follows from \cite{gro-dieu-ega1new}*{Chapitre 0, Proposition 5.5.4}, and the equivalence with \eqref{item:free-vals-iv} can be seen from its proof. 
\end{proof}

\thmref{Prop}{free-vals} suggests the following definitions.

\begin{Def}[const-rank-def][morphisms of constant rank]
  Let $X$ be a superspace and $\vphi$ a morphism $\sh E\longrightarrow\sh F$ of $\sh O_X$-modules. We say that $\vphi$ is \Define{split} if $\sh F/\im\vphi$ is locally free.  

  Let $f:X/S\longrightarrow Y/S$ be a morphism of superspaces over $S$ and $x\in X_0$. We say that $f$ is of \Define{locally constant rank over $S$ at $x$} if for some open neighbourhood $U$ of $x$, the tangent map 
  \[
    \sh T_{f/S}:\sh T_{X/S}|_U\longrightarrow(f^*\sh T_{Y/S})|_U
  \]
  is a split morphism of $\sh O_X|_U$-modules. We say $f$ is of \Define{locally constant rank over $S$} if it is of locally constant rank over $S$ at $x$ for any $x\in X_0$.
\end{Def}

\begin{Cor}[locconst-char]
  Let $f:X/S\longrightarrow Y/S$ be a morphism over $S$ and $x\in X_0$, where $X/S\in\ssplfg{S}$ and $Y/S$ is a supermanifold over $S$. Then the following are equivalent:
  \begin{enumerate}
    \item\label{item:locconst-char-i} The morphism $f$ has locally constant rank over $S$ at $x$.
    \item\label{item:locconst-char-ii} For every $x'$ in a neighbourhood of $x$, the map 
    \[
      (\im\sh T_{f/S})(x')\longrightarrow(f^*\sh T_{Y/S})(x')=T_{f_0(x')}(Y/S)
    \]
    induced by the inclusion $\im\sh T_{f/S}\longrightarrow f^*(\sh T_{Y/S})$ is injective. 
    \item\label{item:locconst-char-iii} There exist an open neighbourhood $U\subseteq X_0$ of $x$ and homogeneous bases of $\sh T_{X/S}|_U$ and $f^*\sh T_{Y/S}|_U$ \scth the matrix $M$ of $\sh T_{f/S}|_U$ has the form 
    \[
      M=\begin{Matrix}1
        A&0\\ 0&D
      \end{Matrix},\quad
      A=
      \begin{Matrix}1
        1&0\\ 0&0
      \end{Matrix},\quad
      D=
      \begin{Matrix}1
        1&0\\ 0&0
      \end{Matrix}.    
    \]
  \end{enumerate}
\end{Cor}

\begin{proof}
  Locally, $X$ admits an embedding into a supermanifold $Z/S$ over $S$, so that locally, $\sh T_{X/S}$ injects into $\sh T_{X/S\to Z/S}$, which is finite locally free by \thmref{Prop}{tan-coord}. Hence, $\sh T_{X/S}$ is of finite type. By the same proposition, $f^*(\sh T_{Y/S})$ is finite locally free. Therefore, the assumptions of \thmref{Prop}{free-vals} are verified, which proves the assertion.
\end{proof}

With the above definition and corollary, we generalise the rank theorem \cite{leites}*{Theorem 2.3.9, Proposition 3.2.9} in two respects: First, one may consider supermanifolds and morphisms over a general base superspace $S$. Secondly, we show the regularity not only of fibres, but also of the inverse images of subsupermanifolds of the image.

\begin{Prop}[constant-rank-fibprod][rank theorem]
  Let $X/S$ and $Y/S$ be in $\SMan_S$, and let $f:X/S\longrightarrow Y/S$ be a morphism of locally constant rank over $S$. Then the following statements hold true:
  \begin{enumerate}[wide]
    \item\label{item:constant-rank-fibprod-i} For any $x\in X_0$, there is an open subset $U\subseteq X_0$, so that the morphism $f|_U$ factors as $f|_U=j\circ p$. Here, $j:Y'/S\longrightarrow Y/S$ is an injective local embedding of supermanifolds over $S$ and $p:X|_U/S\longrightarrow Y'/S$ is a surjective submersion over $S$. 

    Moreover, we may take $Y'=(Y_0',\sh O_{Y'})$, where $Y_0'\defi f_0(U)$, endowed with the quotient topology with respect to $f_0$, and $\sh O_{Y'}\defi(\sh O_Y/\sh J)|_{Y_0'}$, $\sh J\defi\ker f^\sharp$. The morphism $j$ is given by taking $j_0$ equal to the embedding of $Y_0'$ into $Y_0$, and $j^\sharp$ the quotient map with respect to the ideal $\sh J$. 
    \item\label{item:constant-rank-fibprod-ii} If $f':X'/S\longrightarrow Y/S$ is an embedding of supermanifolds over $S$ with $f'_0(X'_0)\subseteq f_0(X_0)$ and ideal $\sh J'\supseteq\sh J$, then the fibre product $X'\times_YX$ exists as a supermanifold over $S$, and the projection $p_2:X'\times_YX\longrightarrow X$ is an embedding over $S$. We have 
    \begin{equation}\label{eq:constrk-pullback-dim}
      \dim_S(X'\times_YX)=\dim_SX'+\dim_SX-\dim_SY'.
    \end{equation}
  \end{enumerate}

  The supermanifold $Y'/S$ over $S$ constructed in item \eqref{item:constant-rank-fibprod-i} is called the \Define{image of $f|_U$}. For the assertion of item \eqref{item:constant-rank-fibprod-ii} to hold, it is sufficient to assume that $f$ has locally constant rank over $S$ at any $x\in f_0^{-1}(f'_0(x'))$, for any $x'\in X'_0$.
\end{Prop}

\begin{proof}
  The statement of \eqref{item:constant-rank-fibprod-i} is well-known in case $S=*$ \cite{leites}*{Theorem 2.3.9}, in view of \thmref{Cor}{locconst-char}. By \thmref{Th}{invfun-loc}, the inverse function theorem holds over a general base. Thus, by \thmref{Cor}{locconst-char}, the proof of the rank theorem carries over with only incremental changes to the general case. 

  As for \eqref{item:constant-rank-fibprod-ii}, the assumption clearly implies that $f'$ factors through $j$ to an embedding $p':X'/S\longrightarrow Y'/S$ over $S$. Since $p$ is a submersion over $S$, the fibre product $X'\times_{Y'}X$ exists, and has the fibre dimension stated on the right-hand side of \eqref{eq:constrk-pullback-dim}. (See \cite{leites}*{Lemma 3.2.8} for the case of $S=*$, the proof of which applies in general, appealing again to \thmref{Th}{invfun-loc} and its usual corollaries.) 

  Since $j$ is an injective local embedding, it is a monomorphism, and it follows that $X'\times_{Y'}X$ is actually the fibre product of $f'$ and $f$. We have a commutative diagram 
  \[
    \begin{tikzcd}
      X'\times_{Y'}X\dar[swap]{p_1}\rar{p_2}&X\dar[swap]{p}\arrow{rdd}{f}\\
      X'\arrow[swap]{rrd}{f'}\rar{p'}&Y'\arrow{rd}[description]{j}\\
      &&Y
    \end{tikzcd}
  \]
  of morphisms over $S$ \scth the left upper square is a pullback whose lower row is an embedding. In particular, $p_{2,0}$ is injective. The image of $p_{2,0}$ is the locally closed subset $f_0^{-1}(f'_0(X_0'))$ of $X_0$. 

  To show that this map is closed, we shall show that it is proper. Let $K\subseteq X_0$ be a compact subset and $L\defi p_0^{\prime-1}(p_0(K))$, which is a compact subset of $X_0'$. Then $p_{2,0}^{-1}(K)$ is a closed subset of $\smash{(X'\times_YX)_0=X'_0\times_{Y'_0}X_0}$ whose image in $X_0'\times X_0$ is contained in $L\times K$. Thus, $\smash{p_{2,0}^{-1}(K)}$ is compact and $p_{2,0}$ is proper, hence closed by \cite{bourbaki-gt1}*{Chapter I, \S 10, Propositions 1 and 7}. Moreover, $\smash{p_2^\sharp}$ is a surjective sheaf map. Hence, $p_2$ is an embedding. 
\end{proof}

\begin{Rem}[rk-converse]
  From the relative inverse function theorem (\thmref{Th}{invfun-loc}), it is clear that the usual normal form theorems hold for submersions and immersions over $S$. Therefore, the converse of \thmref{Prop}{constant-rank-fibprod} holds: Any morphism $f:X/S\longrightarrow Y/S$ which factors as $f=j\circ p$ where $p$ is a submersion over $S$ and $j$ is an immersion over $S$ has locally constant rank over $S$.
\end{Rem}

\subsection{Isotropies at generalised points}\label{subs:isotropy}

In what follows, fix a Lie supergroup $G$ (\ie a group object in $\SMan_\knums$) and an action $a:G\times X\longrightarrow X$ of $G$ on a supermanifold $X$. Let $T\in\ssplfg{\knums}$ and $x\in_TX$ be a $T$-valued point. We recall from Equation \eqref{eq:orbmap-def} the definition of the orbit morphism through $x$, 
\[
  a_x:G_T/T=(T\times G)/T\longrightarrow X_T/T=(T\times X)/T,
\]
by 
\[
  a_x(t,g)=\Parens1{t,a(g,x(t))}=\Parens1{t,g\cdot x(t)},\quad\forall (t,g)\in_RG_T,
\]
and for any $R\in\ssplfg{\knums}$. When $T=*$ is the singleton space, \ie $x\in X_0$ is an ordinary point, then $a_x:G\longrightarrow X$ is the usual orbit morphism \cite{ccf}*{Definition 8.1.4}.

Let $\ger g$ be the Lie superalgebra of $G$, \ie the set of left-invariant vector fields on $G$. This is a Lie superalgebra over $\knums$. For $v\in\ger g$, let $a_v\in\Gamma(\sh T_X)$ denote the \Define{fundamental vector field} induced by the action. It is characterised by the equality 
\begin{equation}\label{eq:fundvf-def}
  (v\otimes1)\circ a^\sharp=(1\otimes a_v)\circ a^\sharp.
\end{equation}
Let $x\in_TX$ with $T\in\ssplfg{\knums}$. The equation above specialises to 
\begin{equation}\label{eq:fundvf-def-pt}
  \begin{split}
    (1\otimes v)\circ a_x^\sharp&=(p_1,\sigma)^\sharp\circ(1\otimes 1\otimes(x^\sharp\circ a_v))\circ({\id}_T\times a)^\sharp\\
    &=(\Delta_T\times{\id}_G)^\sharp\circ(1\otimes(x^\sharp\circ a_v)\otimes 1)\circ({\id}_T\times(a\circ\sigma))^\sharp    
  \end{split}
\end{equation}
where we denote the flip by $\sigma$ and the diagonal morphism by $\Delta_T$. Let $\sh A_\ger g$ be the \Define{fundamental distribution}, \ie the submodule
\[
  \sh A_\ger g\defi\sh O_X\cdot a_\ger g\subseteq\sh T_X,\quad
  a_\ger g\defi \Set1{a_v}{v\in\ger g}.
\]

We shall need to understand when the orbit morphism $a_x$ for an arbitrary $x\in_TX$ is a morphism of locally constant rank over $T$. The following is a full characterisation. 

\begin{Th}[action-locconst]
  Let $x\in_TX$. The following statements are equivalent:
  \begin{enumerate}[wide]
    \item\label{item:action-locconst-i} The morphism $a_x:X_T\longrightarrow G_T$ has locally constant rank over $T$.
    \item\label{item:action-locconst-ii} The pullback $x^*(\sh A_\ger g)$ is a locally direct summand of the $\sh O_T$-module $x^*(\sh T_X)$.
    \item\label{item:action-locconst-iii} For every $t\in T_0$, the canonical map $\sh A_\ger g(x_0(t))\longrightarrow T_{x_0(t)}X$ is injective.
    \item\label{item:action-locconst-iv} For any $t\in T_0$, there are homogeneous $v_j\in\ger g$ \scth $a_{v_j}(x_0(t))\in T_{x_0(t)}X$ are linearly independent and the $x^\sharp\circ a_{v_j}$ span $x^*(\sh A_\ger g)$ in a neighbourhood of $t$ in $T_0$.
  \end{enumerate}
\end{Th}

In the \emph{proof}, we use the following two lemmas. 

\begin{Lem}[pb-vals]
  Let $f:Y\longrightarrow Z$ be a morphism of superspaces and $\sh E$ an $\sh O_Z$-module. Fix $y\in Y_0$. Then the map $e\longmapsto 1\otimes e:\sh E_{f_0(y)}\longrightarrow(f^*\sh E)_y$ induces an isomorphism 
  \[
    \vkappa_Y(y)\otimes_{\vkappa_Z(f_0(y))}\sh E(f_0(y))\longrightarrow(f^*\sh E)(y)
  \]
  of $\vkappa_Y(y)$-super vector spaces. 
\end{Lem}

\begin{proof}
  Let $z\defi f_0(y)$. Now simply note that $\vkappa_Y(y)$ is an $\sh O_{Z,z}$-module \via the map $f^\sharp_x:\sh O_{Z,z}\longrightarrow\sh O_{Y,y}$. In particular, we have
  \[
    (f^*\sh E)(y)=\vkappa_Y(y)\otimes_{\sh O_{Y,y}}\sh O_{Y,y}\otimes_{\sh O_{Z,z}}\sh E_z=\vkappa_Y(y)\otimes_{\sh O_{Z,z}}\sh E_z=\vkappa_Y(y)\otimes_{\vkappa_Z(z)}\sh E(z).      
  \]
  This proves our claim.
\end{proof}

\begin{Lem}[im-funddist]
  The map 
  \[
    x^*(\sh A_\ger g)\longrightarrow ({\id}_T,1_G)^*\Parens1{\im\sh T_{a_x/T}}:w\longmapsto\Delta_T^\sharp\circ(1\otimes w)
  \]
  is an isomorphism.
\end{Lem}

\begin{proof}
  First, we define a map $\vphi:\sh T_{x:T\to X}\longrightarrow({\id}_T,1_G)^*\Parens0{\sh T_{({\id}_T,x):T\to X_T/T}}$ by
  \[
    \vphi(w)\defi\Delta_T^\sharp\circ(1\otimes w).
  \]
  It admits a left inverse $\psi$, defined by
  \[
    \psi(u)\defi u\circ p_2^\sharp
  \]
  where $p_2:X_T\longrightarrow X$ is the second projection. Indeed,
  \[
    \psi(\vphi(w))=\Delta_T^\sharp\circ (1\otimes w)\circ p_2^\sharp=w.
  \]
  Moreover, we have
  \[
    \begin{split}
      \vphi(x^\sharp\circ a_v)&=\Delta_T^\sharp\circ(1\otimes (x^\sharp\circ a_v))\circ({\id}_T\times a(1_G,\cdot))^\sharp\\
      &=({\id}_T,1_G)^\sharp\circ (\Delta_T\times{\id}_G)^\sharp\circ (1\otimes(x^\sharp\circ a_v)\otimes 1)\circ({\id_T}\times(a\circ\sigma))^\sharp\\
      &=({\id}_T,1_G)^\sharp\circ v\circ a_x^\sharp=({\id}_T,1_G)^\sharp\circ\sh T_{a_x/T}(v)    
    \end{split}
  \]
  by \eqref{eq:fundvf-def-pt}, so $\vphi$ descends to a map $x^*(\sh A_\ger g)\longrightarrow ({\id}_T,1_G)^*\Parens1{\im\sh T_{a_x/T}}$, as claimed. The above computation also shows that 
  \[
    \psi(({\id}_T,1_G)^\sharp\circ\sh T_{a_x/T}(v))=\psi(\vphi(x^\sharp\circ a_v))=x^\sharp\circ a_v,
  \]
  so this map admits a left inverse. But the $({\id}_T,1_G)^\sharp\circ\sh T_{a_x/T}(v)$ for $v\in\ger g$ generate $({\id}_T,1_G)^*(\im\sh T_{a_x/T})$, so $\vphi$ is surjective, and is an isomorphism. 
\end{proof}

\begin{proof}[\prfof{Th}{action-locconst}]
  For every $(t,g)\in T_0\times G_0$, we consider the canonical map 
  \begin{equation}\label{eq:orb-canmap}
    \iota_{(t,g)}\defi\iota_{\sh T_{a_x/T,(t,g)}}:\Parens1{\im\sh T_{a_x/T}}(t,g)\longrightarrow T_{g\cdot x_0(t)}(X_T/T).
  \end{equation}
  By \thmref{Cor}{locconst-char}, the morphism $a_x$ is of locally constant rank over $T$ if and only if \fa $(t,g)$, the map $\iota_{(t,g)}$ is injective. Since $a_x$ is $G$-equivariant, it is equivalent that it be injective for all points of the form $(t,1)$, where $t\in T_0$ is arbitrary. 

  By \thmref{Lem}{im-funddist}, we have $x^*(\sh A_\ger g)\cong({\id}_T,1_G)^*(\im\sh T_{a_x/T})$. Because all residue fields in question are equal to $\knums$, \thmref{Lem}{pb-vals} shows that
  \[
    \begin{split}
      \Parens1{\im\sh T_{a_x/T}}(t,1)&=\Parens1{x^*(\sh A_\ger g)}(t)=\sh A_\ger g(x_0(t)),\\
      \Parens1{a_x^*\Parens1{\sh T_{X_T/T}}}(t,1)&=(\sh T_{X_T/T})(t,x_0(t))=T_{x_0(t)}X=(x^*(\sh T_X))(t).
    \end{split}
  \]
  Thus, by \thmref{Prop}{free-vals}, conditions \eqref{item:action-locconst-i}--\eqref{item:action-locconst-iii} are equivalent.

  If \eqref{item:action-locconst-iv} holds, then $\iota_{(t,1)}$ maps a generating set of $(x^*\sh A_\ger g)(t)=(\im\sh T_{a_x/T})(t,1)$ to a basis of $T_{x_0(t)}X$, so it is injective and \eqref{item:action-locconst-i} holds. Conversely, assume \eqref{item:action-locconst-ii} and \eqref{item:action-locconst-iii}. Thus, we may choose homogeneous $v_j\in\ger g$ \scth $a_{v_j}(x_0(t))\in T_{x_0(t)}X$ are linearly independent and span the image of $(x^*(\sh A_\ger g))(t)$. By assumption, the canonical images of the $x^\sharp\circ a_{v_j}$ in $(x^*(\sh A_\ger g))(t)$ are linearly independent, so that $(x^\sharp\circ a_{v_j})_t$ form a minimal generating set of $(x^*(\sh A_\ger g))_t$ by the Nakayama Lemma. Since this module is free, they form a basis. Since $x^*(\sh A_\ger g)$ is finite locally free, \cite{gro-to}*{4.1.1} shows that the $x^\sharp\circ a_{v_j}$ form a local basis of sections, proving \eqref{item:action-locconst-iv}, and thus, the theorem.
\end{proof}

\begin{Cor}[ordpt]
  Let $T=*$ and $x\in X_0$. Then the orbit morphism $a_x:G\longrightarrow X$ has locally constant rank. 
\end{Cor}

\begin{proof}
  In this case, $x^*\sh E=\sh E(x_0(*))$ for any $\sh O_X$-module $\sh E$, so the condition \eqref{item:action-locconst-ii} of \thmref{Th}{action-locconst} becomes void.
\end{proof}

We now apply these general results to the problem of the representability of the isotropy supergroup functor. To that end, we define for $t\in T_0$:
\[
  \ger g_x(t)\defi\Set1{v\in\ger g}{a_v(x_0(t))=0}.
\]

\begin{Th}[trans-iso]
  Let $x\in_TX$ with $T\in\ssplfg{\knums}$. Assume that $a_x$ has locally constant rank over $T$. Then the functor $G_x:\ssplfg{T}\longrightarrow\Sets$ from \thmref{Def}{isotrop-def} is representable by a supermanifold over $T$ of fibre dimension
  \begin{equation}\label{eq:trans-iso-dim}
    \dim_{T,(t,g)}G_x=\dim_\knums\ger g_x(t).  
  \end{equation}
  The canonical morphism $G_x\longrightarrow G_T$ is a closed embedding.

  Conversely, assume that $G_x$ is representable in $\ssplfg{T}$. Then the canonical morphism $j:G_x\longrightarrow G_T$ is an injective immersion with closed image. If $G_x$ is representable in $\SMan_T$, then $j$ is a closed embedding. 
\end{Th}

\begin{proof}
  By \thmref{Prop}{constant-rank-fibprod}, locally in the domain, the image of $a_x$ exists as a supermanifold over $T$, and has fibre dimension 
  \[
    \dim_{T,x_0(t)}\im a_x=\rk T_{(t,g)}(a_x/T)=\dim\ger g-\dim\ger g_x(t).
  \]
  In view of \thmref{Prop}{constant-rank-fibprod}, it will be sufficient to prove for any superfunction $f$ defined on an open subspace of $X_T$: 
  \[
    a_x^\sharp(f)=0\ \Longrightarrow\ x_T^\sharp(f)=0.
  \]
  But for any supermanifold $R$ and any $t\in_RT$, we have 
  \[
    a_x^\sharp(f)(t,1_{G_T})=f(t,1_{G_T}\cdot x(t))=f(t,x(t))=x_T^\sharp(f)(t),
  \]
  so this statement is manifestly verified. Hence, $G_x$ is representable and the canonical morphism is a closed embedding. The expression for the fibre dimension of $G_x$ follows from Equation \eqref{eq:constrk-pullback-dim}, since $\dim_TT=0$. 

  Conversely, assume the functor $G_x$ is representable in $\ssplfg{T}$. Then $j$ is manifestly a monomorphism, \ie $G_x(R)\longrightarrow G_T(R)$ is injective for any $R\in\ssplfg{\knums}$. Inserting $R=*$, we see that the underlying map is injective with image 
  \[
    \Set1{(t,g)\in T_0\times G_0}{g\cdot x_0(t)=x_0(t)},
  \]
  which is closed. Inserting $R=*[\tau|\theta]$, we see that the tangent map $T_{(t,g)}(j/T)$ is injective for every $(t,g)$, by \thmref{Lem}{der-mor}.

  If $G_x$ is a Lie supergroup, then $j_0$ is a closed topological embedding by \thmref{Th}{imm-lie}, and hence, $j$ is an embedding (as follows from \thmref{Th}{invfun-loc}).
\end{proof}

Specialising \thmref{Th}{trans-iso} by the use of \thmref{Cor}{ordpt}, we recover the case of orbits through an ordinary point first treated by B.~Kostant \cite{kostant} in the setting of Lie--Hopf algebras, by C.P.~Boyer and O.A.~S\'anchez-Valenzuela \cite{bsv} for differentiable Lie supergroups, and by L.~Balduzzi, C.~Carmeli, and G.~Cassinelli \cite{bcc} using a functorial framework and super Harish-Chandra pairs.

\begin{Cor}
  Let $T=*$ and $x\in X_0$. Then $G_x$ is representable by a supermanifold and the canonical morphism $G_x\longrightarrow G$ is a closed embedding. 
\end{Cor}

\subsection{Orbits through generalised points}\label{subs:orbits} 

Having discussed the representability of the isotropy supergroup functor, we pass now to the existence of orbits. In what follows, to avoid heavy notation, we will largely eschew writing $/S$ for morphisms over $S$, instead mostly stating the property of being `over $S$' in words. 

We have the following generalisation of Godement's theorem \cite{ah-berezin}*{Theorem 2.6}, with an essentially unchanged proof. We have added the detail that in this situation, the quotients are universal categorical. 

\begin{Prop}[eq-quot]
  Let $R/S$ be an equivalence relation on $X/S$ in $\SMan_S$, as defined in \thmref{Ex}{ex-groupoid} \eqref{item:ex-groupoid-iii}. Then the following assertions are equivalent:
  \begin{enumerate}[wide]
    \item\label{item:eq-quot-i} The weakly geometric quotient $\pi:X\longrightarrow X/R$ exists in $\SMan_S$ and, as a morphism, is a submersion over $S$.
    \item\label{item:eq-quot-ii} The subsupermanifold $R$ of $X\times_SX$ is closed, and (one of, and hence both of) $s,t:R\longrightarrow X$ are submersions over $S$.
  \end{enumerate}
  If this is the case, then $\pi:X\longrightarrow X/R$ is a universal weakly geometric quotient. The quotient is \Define{effective}, that is, the morphism $(t,s):R\longrightarrow X\times_{X/R}X$ is an isomorphism. Moreover, its fibre dimension is 
  \begin{equation}\label{eq:quot-rel-dim}
    \dim_S(X/R)=2\dim_S X-\dim_SR.
  \end{equation}
\end{Prop}

\begin{proof}
  Apart from that about universal weakly geometric quotients, all statements are proved for $S=*$ in Refs.~\cites{almorox,ah-berezin}. In general, the proof carries over unchanged.

  Let us prove the claim of universality for the weakly geometric quotient. So, let the assumption of item \eqref{item:eq-quot-i} be fulfilled and set $Q\defi X/R$. Then $\pi$ is a submersion over $S$, and hence, $X'\defi Q'\times_QX$ exists in $\SMan_S$ for any $\psi:Q'\longrightarrow Q$, by \cite{ahw-sing}*{Proposition 5.41} and the normal form theorem for submersions over $S$ (which follows from \thmref{Th}{invfun-loc}). By item \eqref{item:eq-quot-ii}, $s$ is also a submersion over $S$. Then so is $\pi\circ s$, and $R'\defi (Q'\times Q')\times_{Q\times Q}R$ exists in $\SMan_S$, where $R$ lies over $Q\times Q$ \via $(\pi\times\pi)\circ (t,s):R\longrightarrow Q\times Q$.

  First, we claim that condition \eqref{item:eq-quot-ii} holds for the equivalence relation $R'/S$ on $X'/S$ in $\SMan_S$. Note that we have a pullback diagram 
  \[
    \begin{tikzcd}
      R'=Q'\times_QR\dar[swap]{s'}\rar{}&R\dar{\pi\circ s}\\
      Q'\rar{\psi}&Q
    \end{tikzcd}
  \]
  Since $\pi\circ s$ is a submersion over $S$, so is $s'$. Next, consider the morphism 
  \[
    R'=(Q'\times Q')\times_{Q\times Q}R\longrightarrow X'\times_SX'=(Q'\times Q')\times_{Q\times Q}(X\times_SX).
  \]
  It is an embedding by \cite{ahw-sing}*{Corollary 5.29}. Thus, the assumption \eqref{item:eq-quot-ii} is verified for $R'$ and $X'$, the weakly geometric quotient $\pi':X'\longrightarrow X'/R'$ exists in $\SMan_S$, and it is a submersion over $S$. It is categorical by \thmref{Cor}{wgeom-quot-is-cat}.

  The morphism $p_1={\id}_{Q'}\times_Q\pi:X'\longrightarrow Q'$ is manifestly $R'$-invariant, so that there is a unique morphism 
  \[
    \vphi:X'/R'\longrightarrow Q',\quad\vphi\circ\pi'={\id}_{Q'}\times_Q\pi.
  \]
  Since so is $p_1$, $\vphi$ is a surjective submersion. 

  To see that it is a local isomorphism, we compute the dimensions of the supermanifolds over $S$ in question. On one hand, we have 
  \[
    \dim_SQ=2\dim_SX-\dim_SR,
  \]
  and on the other, we have
  \begin{align*}
    \dim_SX'/R'&=2\dim_SX'-\dim_S R'\\
    &=2(\dim_{Q'}X'+\dim_SQ')-(\dim_{Q'\times Q'}R'+2\dim_SQ')\\
    &=2\dim_Q'X'-\dim_{Q'\times Q'}R'=2\dim_QX-\dim_{Q\times Q}R\\
    &=2(\dim_QX+\dim_SQ)-(\dim_{Q\times Q}R+2\dim_SQ)\\
    &=2\dim_SX-\dim_SR
  \end{align*}
  Upon invoking the inverse function theorem (\thmref{Th}{invfun-loc}), this proves that $\vphi$ is a local isomorphism over $S$. Finally, we need to show that $\vphi_0$ is injective. 

  To that end, let $q_j'\in Q_0'$, $x_j\in X_0$, \scth $\psi_0(q_j')=\pi_0(x_j')$. Assume that $\vphi_0(\pi_0'(q_1',x_1))=\vphi_0(\pi_0'(q_2',x_2))$, so that $q_1'=q_2'$, because
  \[
    \vphi_0\circ\pi_0'=p_{1,0}:X'_0=Q_0'\times_{Q_0}X_0\longrightarrow Q_0'.  
  \]
  It follows that $\pi_0(x_1)=\psi_0(q_1')=\psi_0(q_2')=\pi_0(x_2)$, so that $(x_1,x_2)\in R_0$, since $\pi$ is an effective quotient. Then $(q_1',q_2',x_1,x_2)\in R'_0$, so that $\pi_0'(q_1',x_1)=\pi_0'(q_2',x_2)$, proving the injectivity. 
\end{proof}

We now wish to apply this proposition to supergroup actions. Thus, fix a Lie supergroup $G$ and a $G$-supermanifold $X$. Let $x\in_TX$, where $T$ is some supermanifold. We assume that $G_x$ is representable in $\SMan_T$ and that the canonical morphism $G_x\longrightarrow G_T$ is an embedding over $T$ (automatically closed).

We define an equivalence relation $R_x$ on $G_T$ by 
\[
  R_x\defi G_T\times_TG_x,\quad i:R_x\longrightarrow G_T\times_TG_T,
\]
where $i$ is given by 
\[
  i(g,g')\defi(g,gg'),\quad\forall (g,g')\in_{T'/T}\Parens1{G_T\times_TG_x}/T,
\]
and for any supermanifold $T'/T$ over $T$. It is straightforward to check that $i$ is an embedding and indeed, that $R_x$ is an equivalence relation. 

\begin{Prop}[orb-quot]
  Let $G_x$ be representable in $\SMan_T$. Then the universal weakly geometric quotient $\pi_x:G_T\longrightarrow G\cdot x$ of $G_T$ by $G_x$ exists in $\SMan_T$. It is an effective quotient and a submersion over $T$. Its fibre dimension is
  \begin{equation}\label{eq:quot-dim}
    \dim_T G\cdot x=\dim G-\dim_TG_x.
  \end{equation}
\end{Prop}

\begin{proof}
  The underlying map of $G_x\longrightarrow G_T$ is injective and a homeomorphism onto its closed image, so it is proper. Therefore, the map underlying the morphism $i:R_x\longrightarrow G_T\times_TG_T$ is closed. The first projection $s$ of $R_x$ is obviously a submersion over $T$. Then \thmref{Prop}{eq-quot} applies, and we reach our conclusion. Equation \eqref{eq:quot-dim} follows from Equation \eqref{eq:quot-rel-dim}, since $\dim_TR_x=\dim G+\dim_TG_x$.
\end{proof}

\begin{Not}
  By abuse of language, the morphism $\tilde a_x:G\cdot x\longrightarrow X_T$ over $T$ induced by $a_x$ will also be called the \Define{orbit morphism}.
\end{Not}

Combining this fact with our previous results, we get the following theorem.

\begin{Th}[orbit]
  Let $x:T\longrightarrow X$. The following are equivalent:
  \begin{enumerate}[wide]
    \item\label{item:orbit-i} The morphism $a_x$ has locally constant rank over $T$.
    \item\label{item:orbit-ii} The isotropy functor $G_x$ is representable in $\SMan_T$.
  \end{enumerate}
  In this case, the canonical morphism $j:G_x\longrightarrow G_T$ is a closed embedding, the weakly geometric and universal categorical quotient $G\cdot x$ exists, $\pi_x:G_T\longrightarrow G\cdot x$ is a surjective submersion over $T$, the fibre dimension of $G\cdot x$ is
  \begin{equation}\label{eq:trans-orb-dim}
    \dim_{T,(t,g\cdot x_0(t))}G\cdot x=\dim G-\dim\ger g_x(t),\quad\forall (t,g)\in (G_T)_0=T_0\times G_0,
  \end{equation}
  and $\tilde a_x$ is an immersion over $T$.

  Moreover, if $U\subseteq X_0$ is open \scth $a_x|_U$ admits an image in the sense of \thmref{Prop}{constant-rank-fibprod}, then so does $\tilde a_x|_{\pi_{x,0}(U)}$, and these images coincide.
\end{Th}

\begin{proof}
  The implication \eqref{item:orbit-i} $\Rightarrow$ \eqref{item:orbit-ii} is the content of \thmref{Th}{trans-iso}.

  Conversely, let $G_x$ be representable in $\SMan_T$. Then $j$ is a closed embedding, by \thmref{Th}{trans-iso}. From \thmref{Prop}{orb-quot}, we conclude that $G\cdot x$ exists and $\pi_x:G_T\longrightarrow G\cdot x$ is a surjective submersion over $T$. Because
  \[
    \ker T_{(t,g)}(\pi_x/T)=\ger g_x(t)=\ker T_{(t,g)}(a_x/T)
  \]
  and $\tilde a_x\circ\pi_x=a_x$, it follows that $\tilde a_x$ is an immersion over $T$. By \thmref{Rem}{rk-converse}, $a_x$ is of locally constant rank over $T$. This shows that \eqref{item:orbit-ii} holds. Equation \eqref{eq:trans-orb-dim} follows from Equation \eqref{eq:quot-dim} and Equation \eqref{eq:trans-iso-dim}. Moreover, since $\pi_{x,0}$ is surjective and $\pi_x^\sharp$ is injective, it follows that the images of $a_x|_U$ and $\tilde a_x|_{\pi_{x,0}(U)}$ are equal whenever one of the two is defined, proving the asserted equivalence. The remaining statements follow from \thmref{Prop}{orb-quot}. 
\end{proof}

\section{Coadjoint superorbits and their super-symplectic forms}\label{sec:coad}

In this section, we construct the Kirillov--Kostant--Souriau form in the setting of coadjoint superorbits through $T$-valued points. For the case of ordinary points, where $T=*$, coadjoint orbits of supergroups were studied by B.~Kostant \cite{kostant}, R.~Fioresi and M.A.~Lled\'o \cite{fl}, and by H.~Salmasian \cite{salmasian}.

By the introduction of the parameter space $T$, it is always possible to work with \emph{even} supersymplectic forms, provided they are considered over $T$. Compare with the work of Tuynman \cites{tuyn10a,tuyn10b}, who is obliged to work with inhomogeneous forms. 

We will follow the notation and conventions of Sections \ref{sec:groupoid-quots}--\ref{sec:super-quot} and Ref.~\cite{ap-sphasym}, only briefly recalling the basic ingredients. Let $G$ be a Lie supergroup---\ie a group object in $\SMan_\knums=\SMan_{\knums,\Bbbk}^\varpi$---with Lie superalgebra $\ger g$. We set $\ger g_\Bbbk\defi\ger g_{\Bbbk,\ev}\oplus\ger g_\odd$, where $\ger g_{\Bbbk,\ev}$ is the Lie algebra of $G_0$. (Note that the latter is a $\Bbbk$-form of $\ger g_\ev$.) The dual $\knums$-super vector space of $\ger g$ will be denoted by $\ger g^*$. Let $\ger g_\Bbbk^*$ be the set of $\knums$-linear functionals $f=f_\ev\oplus f_\odd\in\ger g^*$ \scth $f_\ev(\ger g_{\Bbbk,\ev})\subseteq\Bbbk$. We denote the adjoint action of $G$ on $\aff(\ger g_\Bbbk)$ by $\Ad$. 

The \Define{coadjoint action} is defined by 
\[
  \Dual1{\Ad^*(g)(f)}x\defi\Dual1f{\Ad(g^{-1})(x)},\quad\forall g\in_TG,x\in_T\aff(\ger g_\Bbbk),f\in_T\aff(\ger g_\Bbbk^*),
\]
where $\Dual0\cdot\cdot$ denotes the canonical pairing of $\ger g^*$ and $\ger g$.

\subsection{The super-symplectic Kirillov--Kostant--Souriau form}

Let $T\in\ssplfg{\knums}$ and $f\in_T\aff(\ger g_\Bbbk^*)$ be a $T$-valued point of the dual of the Lie superalgebra $\ger g$. We define an even super-skew symmetric tensor $\Omega_f$,
\[
  \Omega_f:\sh T_{G_T/T}\otimes_{\sh O_{G_T}}\sh T_{G_T/T}\longrightarrow\sh O_{G_T},
\]
by the formula 
\[
  \Omega_f\Parens0{v,w}\defi\Dual1{p_{G_T}^\sharp(f)}{[v,w]},\quad\forall v,w\in(\sh O_{G_T}\otimes\ger g)(U),
\]
where $U\subseteq T_0\times G_0$ is open, $p_{G_T}=p_1:G_T\longrightarrow T$, and we identify $f$ with a section of $\sh O_T\otimes\ger g^*$ \via the natural bijection 
\[
  \Hom1{T,\aff(\ger g_\Bbbk^*)}\longrightarrow\Gamma\Parens1{(\sh O_T\otimes\ger g^*)_{\Bbbk,\ev}},
\]
compare \cite{ahw-sing}*{Corollary 4.26, Proposition 5.18}. The identification is \via
\[
  f^\sharp(x)=\Dual0fx,\quad\forall\,x\in\ger g\subseteq\Gamma(\sh O_{\aff(\ger g_\Bbbk^*)}).
\]

From now on and until the end of this subsection, assume that $G_f$ is representable in $\SMan_T$, so that in particular, $G\cdot f$ exists in $\SMan_T$, by \thmref{Th}{orbit}. 

\begin{Lem}
  The $2$-form $\Omega_f$ descends to a well-defined even super-skew symmetric tensor $\tilde\omega_f$,
  \[
    \tilde\omega_f:\sh T_{G_T/T\to G\cdot f/T}\otimes_{\sh O_{G_T}}\sh T_{G_T/T\to G\cdot f/T}\longrightarrow\sh O_{G_T},
  \]
  by the formula 
  \[
    \tilde\omega_f\Parens1{\sh T_{\pi_f/T}(v),\sh T_{\pi_f/T}(w)}\defi\Dual1{p_{G_T}^\sharp(f)}{[v,w]},\quad\forall v,w\in(\sh O_{G_T}\otimes\ger g)(U),
  \]
  for every open $U\subseteq(G_T)_0$. The $2$-form $\tilde\omega_f$ is non-degenerate. 
\end{Lem}

\begin{proof}
  Let $v\in\sh T_{G_T/T}(U)$ be homogeneous and $x\in\ger g\subseteq\Gamma(\sh O_{\aff(\ger g_\Bbbk^*)})$. Let $(x_j)$ be a homogeneous basis of $\ger g$ and expand $v=\sum_jv^jx_j$.

  Then we compute \fa $R\in\ssplfg{\knums}$ and all $\mu\in_R\aff(\ger g_\Bbbk^*)$ that 
  \begin{align*}
    a_{x_j}(\mu)(x)&=\frac d{d\tau}\Big|_{\tau=0}\Dual1{\Ad^*(\phi^{x_j})(\mu)}x=\frac d{d\tau}\Big|_{\tau=0}\Dual1\mu{\Ad(\phi^{-x_j})(x)}\\
    &=-\Dual1\mu{[x_j,x]}=-\mu(\ad(x_j)(x))=-{\ad^*}(x_j)(\mu)(x).
  \end{align*}
  Here, we let $\Abs0\tau=\Abs0{x_j}$ and follow the conventions of \thmref{Def}{inf-flow}.

  Equation \eqref{eq:fundvf-def-pt} shows that 
  \[
    v\circ a_f^\sharp=\sum_jv^j\cdot(\Delta_T\times{\id}_G)^\sharp\circ(1\otimes(f^\sharp\circ a_{x_j})\otimes1)\circ({\id}_T\times(a\circ\sigma))^\sharp.
  \]
  Therefore, \fa $R$ and all $(t,g)\in_RG_T$, we have 
  \[
    \begin{split}
      (v\circ a_f^\sharp)(x)(t,g)&=\sum_jv^j(t,g)\Dual1{\ad^*(x_j)(f(t))}{\Ad(g^{-1})(x)}\\
      &=\sum_jv^j(t,g)\Dual1{\Ad^*(g)(\ad^*(x_j)(f(t)))}x      
    \end{split}
  \]
  Vector fields are uniquely determined by their action on systems of local fibre coordinates, by \thmref{Prop}{tan-coord}. Moreover, any homogeneous basis of $\ger g$ contained in $\ger g_\Bbbk$ defines a global fibre coordinate system on $\aff_T(\ger g_\Bbbk^*)$. Thus, we have 
  \[
    \sh T_{\pi_f/T}(v)=0\quad\Longleftrightarrow\quad\sum_jv^j(t,g)\Ad^*(g)(\ad^*(x_j)(f(t)))=0\quad\forall\,R,(t,g)\in_RG_T.
  \]

  On the other hand, we may express 
  \[
    \begin{split}
      \Dual1{p_{G_T}^\sharp(f)}{[v,w]}(t,g)&=\sum_jv^j(t,g)\Dual1{\ad^*(x_j)(f(t))}{(t,g)^\sharp\circ w}\\
      &=\sum_jv^j(t,g)\Dual1{\Ad^*(g)(\ad^*(x_j)(f(t)))}{\Ad(g^{-1})((t,g)^\sharp\circ w)}.
    \end{split}
  \]
  This shows immediately that $\tilde\omega_f$ is well-defined. Setting $\check w\defi ({\id}_T\times{\Ad})^\sharp\circ w$, the above computation shows that 
  \[
    \Dual1{p_{G_T}^\sharp(f)}{[v,\check w]}(t,g)=\sum_jv^j(t,g)\Dual1{\Ad^*(g)(\ad^*(x_j)(f(t)))}w.
  \]
  Hence, if $\tilde\omega_f(\sh T_{\pi_f}(v),\sh T_{\pi_f}(\check x_j))=0$ for any $j$, then it follows that $v\circ\pi_f^\sharp=0$, so we see that $\tilde\omega_f$ is non-degenerate. 
\end{proof}

Since $G\cdot f\in\SMan_T$, we have 
\[
  \sh T_{G_T/T\to G\cdot f/T}=\pi_f^*(\sh T_{G\cdot f/T}),
\]
by \thmref{Prop}{tan-coord}, so we may ask whether $\tilde\omega_f$ is induced by some tensor $\omega_f$ on $G\cdot f$. Indeed, this is the case, as we presently show.

The module inverse image and direct image functors $((\pi_f)^*,(\pi_f)_{0*})$ form an adjoint pair, so there is a natural bijection 
\[
	\begin{tikzcd}
		\Hom[_{\sh O_{G\cdot f}}]1{\textstyle\bigwedge^2\sh T_{G\cdot f/T},(\pi_f)_{0*}\sh O_{G_T}}
		\rar{\pi_f^*}
		&\Hom[_{\sh O_{G_T}}]1{\textstyle\bigwedge^2\sh T_{G_T/T\to G\cdot f/T},\sh O_{G_T}}.
	\end{tikzcd}
\]

\begin{Prop}[pf-2form]
  There is a unique even super-skew symmetric tensor
  \[
    \omega_f:\sh T_{G\cdot f/T}\otimes_{\sh O_{G\cdot f}}\sh T_{G\cdot f/T}\longrightarrow\sh O_{G\cdot f}
  \]
  \scth $\pi_f^*(\omega_f)=\tilde\omega_f$.
\end{Prop}

\begin{proof}
	By the above, there is a unique even super-skew symmetric tensor
	\[
	    \omega_f:\sh T_{G\cdot f/T}\otimes_{\sh O_{G\cdot f}}\sh T_{G\cdot f/T}\longrightarrow (\pi_f)_{0*}\sh O_{G_T},
	\]
	\scth $\pi_f^*(\omega_f)=\tilde\omega_f$. We need to show that it takes values in the subsheaf $\sh O_{G\cdot f}$.

  But $G\cdot f=G_T/G_f$ is a weakly geometric quotient by \thmref{Prop}{orb-quot}, so that by \thmref{Rem}{colimit-explicit}, we have 
  \[
    \sh O_{G\cdot f}=\Parens0{(\pi_f)_{0*}\sh O_{G_T}}^{G_f}.
  \]
  It thus remains to prove that $\omega_f$ takes values in the sheaf of invariants.

  To that end, fix a homogeneous basis $(x_j)$ of $\ger g$ contained in $\ger g_\Bbbk$. Take any $v,w\in\sh T_{G\cdot f/T}(U)$, where $U\subseteq (G\cdot f)_0$ is open and define $V\defi(\pi_f)_0^{-1}(U)\subseteq T_0\times G_0$. We may write $\pi_f^\sharp\circ v=\sum_jv^j\,(1\otimes x_j)\circ\pi_f^\sharp)$ \fs $v^j\in\sh O_{G_T}(V)$, $\Abs0{v^j}=\Abs0{x_j}+\Abs0v$, and similarly for $w$. 

  Denote by $(t,g,h)$ the generic point of $G_T|_V\times_TG_f|_V$. We compute for any superfunction $k$ on $G\cdot f$, defined on an open subset of $U$, that
  \[
    (\pi_f^\sharp\circ v)(k)(t,gh)=v(k)\Parens1{(t,gh)\cdot f(t)}=v(k)\Parens1{(t,g)\cdot f(t)}=(\pi_f^\sharp\circ v)(k)(t,g).
  \]
  Here, we are using the fact that $G\cdot f$ is a universal categorical quotient (\thmref{Th}{orbit}), so that, by \thmref{Prop}{orbit-action}, it admits a $G$-action for which $\pi_f$ is equivariant and $f$, considered as a $T$-valued point of $G\cdot f$, is fixed by $G_f$.

  On the other hand, using results from Ref.~\cite{ap-sphasym}, we have 
  \begin{align*}
    \sum_j\Parens1{v^j(x_j\circ\pi_f^\sharp)(k)}(t,gh)
    &=\sum_jv^j(t,gh)\frac d{d\tau}\Big|_{\tau=0}k\Parens1{(t,gh\exp(\tau x_j))\cdot f(t)}\\
    &=\sum_jv^j(t,gh)\frac d{d\tau}\Big|_{\tau=0}k\Parens1{(t,g\exp(\tau\Ad(h)(x_j))h)\cdot f(t)}\\
    &=\sum_jv^j(t,gh)\frac d{d\tau}\Big|_{\tau=0}k\Parens1{(t,g\exp(\tau\Ad(h)(x_j)))\cdot f(t)}\\
    &=\sum_jv^j(t,gh)(\Ad(h)(x_j)\circ\pi_f^\sharp)(k)(t,g).
  \end{align*}
  Combining both computations, we arrive at the equality
  \begin{equation}\label{eq:ad-orb}
    \sum_jv^j(t,gh)(\Ad(h)(x_j)\circ\pi_f^\sharp)=\sum_jv(t,g)(x_j\circ\pi^\sharp_f)
  \end{equation}
  of vector fields over $T$ along the morphism 
  \[
    \pi_f\circ m_{G_T}=\pi_f\circ p_1:G_T\times_TG_f\longrightarrow G\cdot f.
  \] 
  Using Equation \eqref{eq:ad-orb}, we may compute
  \begin{align*}
    \omega_f(v,w)(t,gh)&=\tilde\omega_f(\pi_f^\sharp\circ v,\pi_f^\sharp\circ w)(t,gh)\\
    &=\sum_{jk}(-1)^{\Abs0{x_j}\Abs0{x_k}}(v^jw^k)(t,gh)\Dual1{f(t)}{[x_j,x_k]}\\
    &=\sum_{jk}(-1)^{\Abs0{x_j}\Abs0{x_k}}(v^jw^k)(t,gh)\Dual1{f(t)}{[\Ad(h)(x_j),\Ad(h)(x_k)]}\\
    &=\sum_{jk}(-1)^{\Abs0{x_j}\Abs0{x_k}}(v^jw^k)(t,g)\Dual1{f(t)}{[x_j,x_k]}\\
    &=\tilde\omega_f(\pi_f^\sharp\circ v,\pi_f^\sharp\circ w)(t,g)=\omega_f(v,w)(t,g),
  \end{align*}
  which shows that indeed, $\omega_f(v,w)$ is right $G_f$-invariant, and hence, that $\omega_f$ takes values in the sheaf $\sh O_{G\cdot f}$, as desired.
\end{proof}

We may consider $\omega_f$ as a global section of $\Omega^2_{G\cdot f/T}=\bigwedge^2\Omega^1_{G\cdot f/T}$, \ie a $2$-form over $T$. We show that it is closed.

\begin{Prop}
  The $2$-form $\omega_f$ over $T$ is relatively closed. 
\end{Prop}

\begin{proof}
  The element of $\Gamma(\sh O_{G_T}\otimes\ger g^*)$ corresponding to $f$ is a left $G$-invariant $1$-form (which is, moreover, even and real-valued). We show that it gives a potential for the pullback of $\omega_f$. To that end, we follow ideas of Chevalley--Eilenberg \cite{ce}.

  Let $v,w\in\ger g$. Denote by $d=d_{G_T/T}$ the relative differential. Then 
  \[
    \iota_wd+d\iota_w=\sh L_w,
  \]
  where $\iota_v$, $\Abs0{\iota_v}=\Abs0v$, denotes relative contraction, and $\sh L_v$, $\Abs0{\sh L_v}=\Abs0v$, denotes the relative Lie derivative. We have 
  \begin{align*}
    df(v,w)&=(-1)^{\Abs0v\Abs0w}\iota_w\iota_vdf=(-1)^{\Abs0v\Abs0w}\iota_w(\sh L_vf)\\
    &=-[\sh L_v,\iota_w]f=-\iota_{[v,w]}f=-\Dual1f{[v,w]}=-\tilde\omega_f(\sh T_{\pi_f/T}(v),\sh T_{\pi_f/T}(w)),
  \end{align*}
  since $\iota_wf=\Dual0fw\in\Gamma(\sh O_T)$, so that $d\iota_wf=0=\sh L_v\iota_wf$. Since both sides of the equation are $\sh O_{G_T}$-bilinear, the equation 
  \[
    \tilde\omega_f\Parens1{\sh T_{\pi_f/T}(v),\sh T_{\pi_f/T}(w)}=-df(v,w)
  \]
  holds for any vector fields $v,w$ on $G_T$ over $T$, defined on some open subset. But since $\tilde\omega_f=\pi_f^*(\omega_f)$ by \thmref{Prop}{pf-2form}, we have 
  \[
    \pi_f^\sharp(\omega_f)(v,w)=\tilde\omega_f\Parens1{\sh T_{\pi_f/T}(v),\sh T_{\pi_f/T}(w)}=-df(v,w)
  \]
  for any vector fields $v,w$ on $G_T$ over $T$. Thus, 
  \[
    \pi_f^\sharp(d\omega_f)=d\pi_f^\sharp(\omega_f)=-d^2f=0.
  \]
  Since $\pi_f^\sharp$ is an injective sheaf map, we conclude that $d\omega_f=0$.
\end{proof}

We summarise the above results in the following theorem.

\begin{Th}[coadj-sympl]
  Let $G$ be a Lie supergroup with Lie superalgebra $\ger g$. Let $T\in\ssplfg{\knums}$ and $f:T\longrightarrow\aff(\ger g_\Bbbk^*)$ be \scth $G_f$ is representable and $G_f\longrightarrow G_T$ is an embedding. Then the coadjoint orbit $G\cdot f$ exists, is universal categorical, and with the Kirillov--Kostant--Souriau form $\omega_f$, $G\cdot f$ is a supersymplectic supermanifold over $T$. The assumption is verified if the equivalent conditions in \thmref{Th}{action-locconst} hold.
\end{Th}

\section{Application: Glimpses of the superorbit method}\label{sec:quant}

This section offers an application of our general theory of coadjoint orbits to the geometric construction of representations. By way of example, we show how the formalism can be applied to give certain `universal' $T$-families of representations of certain Lie supergroups, namely, the Abelian supergroup $\aff^{0|n}$ and certain graded variants of the $3$-dimensional Heisenberg group. 

At this point, we will only partially address the issue to which extent unitary structures exist on these families, nor will we make precise in which sense they are universal. We intend to treat these issues in forthcoming work, together with an extension to more general Lie supergroups. 

\subsection{Representations of Lie supergroups over some base}

Fix $T\in\ssplfg{\knums}$. To set the stage both for the general representation theory of supergroups over $T$ and in particular, for the examples to be considered below, we give some very general definitions. 

The functor $\sh O:\Parens1{\ssplfg{T}}^{\mathrm{op}}\longrightarrow\Sets$ defined on objects $U/T$ in $\ssplfg{T}$ by 
\[
  \sh O(U/T)\defi\Gamma(\sh O_{U,\ev})
\]
and on morphisms $f:U'/T\longrightarrow U/T$ in $\ssplfg{T}$ by 
\[
  \sh O(f):\sh O(U/T)\longrightarrow\sh O(U'/T):h\longmapsto f^\sharp(h)
\]
is a ring object in the category $\Bracks1{\Parens1{\ssplfg{T}}^{\mathrm{op}},\Sets}$.

\begin{Def}
  Let $G$ be a supergroup over $T$. A \Define{representation of $G$} is a pair $(\sh H,\pi)$ consisting of:
  \begin{enumerate}[wide]
    \item a $\ints/2\ints$-graded $\sh O$-module object $\sh H:\Parens1{\ssplfg{T}}^{\mathrm{op}}\longrightarrow\Sets$ and 
    \item a morphism $\pi:G\times\sh H\longrightarrow\sh H$, denoted by
    \[
      \pi(g)\psi,\quad \forall\,U\in\ssplfg{\knums},t\in_UT,g\in_tG,\psi\in\sh H(t),
    \]
  \end{enumerate}
  \scth $\pi(g)$ leaves the homogeneous components of $\sh H$ invariant and the following equations are satisfied: 
  \begin{equation}\label{eq:linear-action}
    \begin{gathered}
      \pi(1_G(t))\psi=\psi,\quad\pi(g_1g_2)\psi=\pi(g_1)\Parens1{\pi(g_2)\psi},\\
      \pi(g)(\lambda \psi_1+\psi_2)=\lambda\pi(g)\psi_1+\pi(g)\psi_2,
    \end{gathered}
  \end{equation}
  \fa $t\in_UT$, $g,g_1,g_2\in_tG$, $\psi,\psi_1,\psi_2\in\sh H(t)$, $\lambda\in\sh O(t)$.

  A graded $\sh O$-submodule $\sh H'\subseteq\sh H$ is a \Define{$G$-subrepresentation} if it is $G$-invariant, \ie if $\pi$ descends to a morphism $G\times\sh H'\longrightarrow\sh H'$.
\end{Def}

This concept generalises the existent notions of representations of Lie supergroups in several ways. To make contact with the literature, recall the following construction, which produces a graded $\sh O$-module for any $\sh O_T$-module: Let $\sh H$ be a (graded) $\sh O_T$-module. Define for $t\in_UT$:
\[
  \sh H(t)\defi\Gamma\Parens1{(t^*\sh H)_\ev},\quad\sh H_i(t)\defi\Gamma\Parens1{(t^*\sh H_i)_\ev}
\]
(where $(-)_\ev$ refers to the total grading) and for any commutative diagram 
\[
  \begin{tikzcd}
    U'\arrow{rd}[swap]{t'}\arrow{rr}{f}&&U\arrow{ld}{t}\\ &T,
  \end{tikzcd}
\]
set
\[
  \sh H(f)\defi\Gamma(f^\sharp\otimes1):\sh H(t)\longrightarrow\sh H(t'),
\]
where $\Gamma$ denotes the global sections functor, as usual. The $\sh O$-module structure is given by 
\[
  \sh O(t)\times\sh H(t)\longrightarrow\sh H(t):(h,\psi)\longmapsto h\cdot\psi,
\]
where $\cdot$ is the module structure on global sections. 

In particular, for $T=*$, any super-vector space $V$ over $\knums$ defines such an $\sh O$-module. Assume that $V$ is finite-dimensional and we are given a linear action $\pi:G\times\sh H\longrightarrow\sh H$ where $\sh H=\aff^\knums(V)$ is the functor given on objects $U$ by $\sh H(U)=\Gamma\Parens1{(\sh O_U\otimes V)_\ev}$ and linear actions are defined by the identities in Equation \eqref{eq:linear-action}. Then we may define a representation $(\sh H,\pi)$ by 
\[
  \pi(g)\psi\defi\pi(g,\psi),
\]
\fa $U\in\ssplfg{\knums}$, $g\in_UG$, and $\psi\in\sh H(U)=\Gamma\Parens1{(\sh O_U\otimes V)_\ev}$.

If $G$ is a Lie supergroup, then a linear action is the same thing as a representation of the associated supergroup pair $(G_0,\ger g)$, compare \cite{a-supercw}*{Proposition 1.5}. For the affine algebraic case, compare also \cite{ccf}*{Definition 11.7.2}.

\begin{Ex}[reg-rep][the left-regular representation]
  Let $G$ be a supergroup over $T$. The \Define{left-regular representation} $(\sh H_G,\lambda_G)$ of $G$ is defined by taking 
  \[
     \sh H_G(t)\defi\Gamma\Parens1{\sh O_{U\times_TG,\ev}}
  \] 
  \fa $t\in_UT$,
  \[
    \sh H_G(f)\defi (f\times_T{\id}_G)^\sharp:\sh H_G(t)\longrightarrow\sh H_G(t')
  \]
  \fa $f:t'\longrightarrow t$, and 
  \[
    \lambda_G(g)\psi\defi\Parens1{({\id}_U\times_Tm_G)\circ(({\id}_U,g^{-1})\times_T{\id}_G)}^\sharp(\psi)=\psi(-,g^{-1}(-))
  \]
  \fa $t\in_UT$, $g\in_tG$, $\psi\in\sh H_G(t)=\Gamma\Parens1{\sh O_{U\times_TG,\ev}}$. Here, $g^{-1}$ is $i_G(g)$, as usual.
\end{Ex}

Let $G$ be a supergroup over $T$. By definition, the \Define{Lie superalgebra} of $G$ is the $\sh O_T$-submodule $\ger g$ of the direct image sheaf $p_{0*}(\sh T_{G/T})$ of \Define{left-invariant vector fields} on $G$, defined by 
\[
  \ger g(U)\defi\Set1{v\in\sh T_{G/T}(p_0^{-1}(U))}{m_G^\sharp\circ v=(1\otimes v)\circ m_G^\sharp}
\]
for any open set $U\subseteq T_0$, endowed with the usual bracket of vector fields. We may consider $\ger g$ as a functor, as explained above. Then the \Define{derived representation} $L$ of $(\sh H_G,\lambda)$ is the morphism $d\lambda_G:\ger g\times\sh H_G\longrightarrow\sh H_G$ defined by $d\lambda_G(v)\defi -L_v$ and 
\[
  L_v\psi\defi\Parens1{(1_{U\times_TG}\times_T{\id}_G)^\sharp\circ (v\otimes1)\circ ({\id}_U\times_Tm_G)^\sharp}(\psi)
\]
\fa $U\in\ssplfg{\knums}$, $t\in_UT$, $v\in\ger g(t)=\Gamma((t^*\ger g)_\ev)$, and $\psi\in\sh H_G(t)=\Gamma\Parens1{\sh O_{U\times_TG,\ev}}$. Here, we have 
\[
  1_{U\times_TG}\defi ({\id}_U,1_G(t)):U\longrightarrow U\times_TG.
\]
Similarly, we define $R:\ger g\times\sh H_G\longrightarrow\sh H_G$ by 
\[
  R_v\psi\defi\Parens1{(1_{U\times_TG}\times_T{\id}_G)^\sharp\circ v_{13}\circ ({\id}_U\times_Tm_G)^\sharp}(\psi)
\]
\fa $t\in_UT$, $v\in\ger g(t)=\Gamma((t^*\ger g)_\ev)$, and $\psi\in\sh H_G(t)=\Gamma\Parens1{\sh O_{U\times_TG,\ev}}$. Here, we define 
\[
  v_{13}\defi((1\,2)^{-1}\times_T{\id}_G)^\sharp\circ (1\otimes v)\circ((1\,2)\times_T{\id}_G)^\sharp,
\]
where $(1\,2):G\times_TU\longrightarrow U\times_TG$ is the flip.

\bigskip\noindent
Let us now indicate how to apply these ideas to transplant the orbit method into the world of Lie supergroups. Let $G$ be a Lie supergroup and $f\in_T\aff(\ger g_\Bbbk^*)$ a $T$-valued point of the dual of the Lie superalgebra $\ger g$ (see the introduction to Section \ref{sec:coad}). The Lie superalgebra of $G_T\defi T\times G$ will be denoted by $\ger g_T$ and equals $\sh O_T\otimes\ger g$, as is easy to see. The representations that appear in the superorbit method are all instances of the following simple construction. 

\begin{Prop}[pol-sect]
  Let $\ger h\subseteq\ger g_T$ be an $\sh O_T$-submodule. Then the graded $\sh O$-submodule $\sh H_f^\ger h\subseteq\sh H_{G_T}$ defined by 
  \[
    \sh H_f^\ger h(t)\defi\Set1{\psi\in\sh H_{G_T}(t)}{\forall v\in\ger h(t)=\Gamma((t^*\ger h)_\ev):R_v\psi=-i\Dual0{f(t)}v\psi}
  \]
  \fa $U\in\ssplfg{\knums}$, $t\in_UT$, is a $G_T$-subrepresentation of $(\sh H_{G_T},\lambda_{G_T})$.
\end{Prop}

\begin{proof}
  The action $R$ is $\sh O$-linear and commutes with $L$.
\end{proof}

In the generality they are defined here, the representations $(\smash{\sh H_f^\ger h,\lambda_{G_T}})$ are not interesting. The relevant case is when $\ger h$ is an $\sh O_T$-Lie subsuperalgebra of $\ger g_T$ that is $\Omega_f$-isotropic, \ie $\Omega_f(v,v')=0$ for all local sections $v,v'$ of $\ger h$. If $\ger h$ is maximally isotropic, then one thinks of $\sh H^\ger h_f$ as the `space of $\ger h$-polarised sections of the canonical line bundle on $G\cdot f$' and, following Kirillov's orbit philosophy, expects suitable `completions' thereof to be irreducible. 

In the classical case of a Lie group over $T=*$, this is indeed true, and under certain assumptions, \eg when $G$ is nilpotent, one thus obtains all irreducible unitary representations \cite{rais}. For nilpotent Lie supergroups, it is known by the work of H.~Salmasian and K.-H.~Neeb \cites{salmasian,ns11} that irreducible unitary representations (in the sense of Ref.~\cites{cctv,cctv-err}) are parametrised by coadjoint orbits through ordinary points of $\ger g^*$. The known constructions of these representations are however somewhat roundabout.

As we show below, by way of example, for the Clifford supergroup of dimension $1|2$, they are realised as certain $\sh H_f^\ger h$. Moreover, we show that our approach, for general $T$, also allows for a Plancherel decomposition of the regular representation, at least for the simplest case of the Abelian Lie supergroup $\aff^{0|n}$, where coadjoint orbits through ordinary points are totally insufficient. 

We believe that these examples are mere inklings of a vastly more general picture covering the representation theory and harmonic analysis of nilpotent and possibly more general Lie supergroups.

\subsection{\texorpdfstring{The Plancherel formula for $\aff^{0|n}$}{The Plancherel formula for the superaffine space}}

Let $G=\aff^{0|n}$ be the additive supergroup of the super-vector space $\ger g=\knums^{0|n}$. The coadjoint action $\Ad^*$ of $G$ is trivial. If we let $T\defi\aff(\ger g^*)$ and consider the generic point $f={\id}_T\in_T\aff(\ger g^*)$, then the following diagram commutes:
\[
  \begin{tikzcd}
    G_T\arrow{rr}{a_f}\arrow{rd}[swap]{p_1}&&\aff(\ger g^*)_T\\
    &T\arrow{ru}[swap]{\Delta_T}
  \end{tikzcd}
\]
Thus, $a_f$ factors as the composition of an embedding with a surjective submersion and thus has constant rank by \thmref{Rem}{rk-converse}. Alternatively, observe that the fundamental distribution $\sh A_\ger g=0$, so that the criterion \eqref{item:action-locconst-ii} of \thmref{Th}{action-locconst} is verified. Moreover, the above factorisation coincides with the standard one into $\pi_f$ and $\tilde a_f$, \ie $G_f=G_T$, $G\cdot f=T$, $\pi_f=p_1$, and $\tilde a_f=\Delta_T$. The Kirillov--Kostant--Souriau form $\omega_f$ is zero. 

The general philosophy of `geometric quantisation' or `Kirillov's orbit method' demands the choice of a \Define{polarising} (\ie maximally isotropic) subalgebra $\ger h\subseteq\ger g_T$.  Since $\Omega_f=0$, we must have $\ger h=\ger g_T$. The corresponding $G_T$-subrepresentation $\sh H=\sh H_f^\ger h$ of $(\sh H_{G_T},\lambda_{G_T})$ is given \fa $U\in\ssplfg{\knums}$, $t\in_UT$, by 
\[
  \sh H(t)\defi\Set1{\psi}{\psi\in\Gamma\Parens1{(\sh O_{G_U})_\ev},\forall v\in\ger g:R_{1\otimes v}\psi=-i\Dual0{f(t)}v\psi}.
\]
This is the functor of a free $\sh O_T$-module of rank $1|0$, since it has the basis of sections 
\[
  \psi_0\defi e^{-i\sum_j\theta_j\xi^j}\in\sh H({\id}_T)=\Gamma(\sh O_{G_T,\ev}).
\]
Here, $\theta_1,\dotsc,\theta_n$ is some arbitrary basis of $\ger g$ and $\xi^1,\dotsc,\xi^n$ is the dual basis of $\ger g^*$, considered as coordinate superfunctions on $T=\aff(\ger g^*)$ and $G=\aff(\ger g)$, respectively.

The representation of $G_T$ on $\sh H$ is determined by its action on the special vector
\[
  \psi_t\defi\sh H_{G_T}(t)(\psi_0)=e^{-i\sum_jt_j\xi^j},\quad t_j\defi t^\sharp(\theta_j).
\]
With $\pi$ denoting the restriction of $\lambda_{G_T}$ to $\sh H$ and $g^j\defi g^\sharp(\xi^j)$, it is given by
\[
  \begin{split}
    \pi(g)\psi_t&=(({\id}_U,g^{-1})\times_T{\id}_G)^\sharp({\id}_U\times_Tm_G)^\sharp(e^{-i\sum_jt_j\xi^j})\\
    &=(({\id}_U,g^{-1})\times_T{\id}_G)^\sharp(e^{-i\sum_jt_j(\xi^j_1+\xi^j_2)})\\
    &=e^{-i\sum_jt_j(-g^j+\xi^j)}=e^{i\Dual0tg}\psi_t,
  \end{split}
\]
that is, it is a character, as was to be expected. 

We have the following `abstract' Fourier inversion formula.

\begin{Prop}[odd-fi]
  For any superfunction $f$ on $G$, we have 
  \[
    \int_TD(\theta)\,\str\pi(f)=(-1)^{n(n+1)/2}i^nf_0(0),
  \]   
  where $\pi(f)$ is defined by 
  \[
    \pi(f)\defi\int_GD(\xi)\,f\pi,
  \]
  and the integrals are Berezin integrals. 
\end{Prop} 

For the Berezin integral, see \citelist{\cite{leites}*{Chapter 2, \S 4} \cite{manin}*{Chapter 4, \S 6} \cite{deligne-morgan}*{\S 3.9}}.

\begin{proof}
  Since $\pi$ is a character, the operator $\pi(f)$ is a function:
  \[
    \pi(f)=\int_G D(\xi)\,fe^{i\sum_j\theta_j\xi^j}\in\Gamma(\sh O_T)
  \]
  Therefore, $\str\pi(f)$ is that same function. (Incidentally, this may be viewed as a baby version of Kirillov's character formula.) The assertion now follows from the Euclidean Fourier inversion formula \cite{as-sbos}*{Proposition C.17}.
\end{proof}

We obtain the following Plancherel formula. 

\begin{Cor}
  \Fa superfunctions $f$ and $g$ on $G$, we have 
  \[
    \int_TD(\theta)\,\str(\pi(f)^\dagger\pi(g))=(-1)^{n(n+1)/2}i^n\int_GD(\xi)\,\overline fg.
  \]
  Here, $(-)^\dagger$ is the super-adjoint with respect to the $\sh O_T$-inner product on $\sh H$ normalised by $\Sdual0{\psi_0}{\psi_0}=1$ and $\overline{\id}$ is the antilinear antiautomorphism of $\sh O_G$ defined by $\overline{\xi^j}=\xi^j$.
\end{Cor}

\begin{proof}
  Using the methods of Ref.~\cite{ahl-cliff}, one sees that $\pi(f)^\dagger\pi(g)=\pi(f^**g)$, where $*$ is the convolution product on $G$ and $f^*=i^\sharp(\overline f)$ where $i$ is the inversion of $G$. Since $\delta=\xi^1\dotsc\xi^n$ is the Dirac delta on $G$, the formula follows from \thmref{Prop}{odd-fi}.
\end{proof}

\begin{Rem}
  Thus, by judiciously applying the orbit method to $T$-valued points, we have obtained a decomposition of the left regular representation of $G$ into an `odd direct integral' of `unitary' characters. By contrast, a direct sum decomposition of the function space $\Gamma(\sh O_G)$ into irreducible unitary $G$-representations is impossible, since the only such representation is the trivial one!  
\end{Rem}

\subsection{The orbit method for Heisenberg type supergroups}

Let us consider the Lie superalgebra $\ger g$ over $\knums$ spanned by homogeneous vectors $x,y,z$ satisfying the unique non-zero relation
\[
  [x,y]=z.
\]

When $x,y,z$ are even, $\ger g$ is the classical Heisenberg algebra of dimension $3|0$. When $x,y$ are odd, $z$ must be even. The central element $z$ spans a copy of $\knums$, so $\ger g$ is a unital Lie algebra in the sense of Ref.~\cite{andler-sahi}, and its unital enveloping algebra $\Uenv0{\ger g}/(1-z)$ is the Clifford algebra $\mathop{\mathrm{Cliff}}(2,\knums)$. (NB: We will use a different normalisation below.) For this reason, $\ger g$ is called the Clifford--Lie superalgebra, and its representation theory was studied \eg in Refs.~\cites{salmasian,ahl-cliff}. The construction of the representations used there is \emph{ad hoc}. Below, we show how they arise in a natural fashion.

A third possibility, which does not seem to have been considered before, is that $x,y$ are of distinct parity (but see Ref.~\cite{frydryszak}). In this case, $z$ is odd. As we show below, besides characters, there exists a family of representations (which happen to be finite-dimensional) parametrised by $T=\aff^{0|1}$, which bear a striking resemblance to the Schr\"odinger representation of the Heisenberg group. 

\subsubsection{Parity-independent computations}

A number of computations concerning the Lie superalgebra $\ger g$ of Heisenberg type introduced above are somewhat independent of the parity of its elements. We begin with the coadjoint representation of $\ger g$. Let $x^*,y^*,z^*$ be the basis dual to $x,y,z$. In terms of this basis, we have
\[
  \ad^*(x)=
  \begin{Matrix}1
    0&0&0\\ 0&0&-(-1)^{\Abs0x\Abs0z}\\ 0&0&0
  \end{Matrix},\quad
  \ad^*(y)=
  \begin{Matrix}1
    0&0&(-1)^{\Abs0y}\\ 0&0&0\\ 0&0&0
  \end{Matrix},\quad
  \ad^*(z)=0.
\]

Recall the definitions given at the beginning of Section \ref{sec:super-quot}. We will consider the field $\Bbbk=\reals$, since we are mainly interested in super versions of real Lie groups. A Lie supergroup $G$ (\ie a group object in the category of supermanifolds over $(\knums,\reals)$ of class $\sh C^\varpi$) with Lie superalgebra $\ger g$ is uniquely determined by the choice of a real Lie group $G_0$ whose Lie algebra is a real form $\ger g_{\reals,\ev}$ of $\ger g_\ev$, compare Ref.~\cite{ap-sphasym}.

We fix $\ger g_\reals\defi\ger g_{\reals,\ev}\oplus\ger g_\odd$ by setting $\ger g_{\reals,\ev}\defi\ger g_\ev\cap\Span0{x,y,z}_\reals$. Let $G$ be the connected and simply connected Lie supergroup whose Lie superalgebra is $\ger g$ and whose Lie group has Lie algebra $\ger g_{\reals,\ev}$. Unless $\ger g$ is purely even, $G_0$ is the additive group of $\reals$. With these conventions, $\ad^*(v)$ is the fundamental vector field corresponding to $v\in\ger g$ under the coadjoint action $\Ad^*$ of $G$.

Let $T\in\ssplfg{\knums}$ be arbitary and $f=\alpha x^*+\beta y^*+\gamma z^*\in_T\aff(\ger g_\reals^*)$. Observe that
\[
  \ad^*(ax+by+cz)(f)=-a\gamma y^*+(-1)^{\Abs0y(1+\Abs0z)}b\gamma x^*
\]
for $v=ax+by+cz\in\ger g$, where $a,b,c\in\knums$. Thus, if $u=a'x+b'y+c'z\in\ger g\subseteq\Gamma(\sh O_{\aff(\ger g_\Bbbk^*)})$ with arbitary $a',b',c'\in\knums$, then 
\[
  \Parens1{f^\sharp\circ\ad^*(v)}(u)=-ab'\gamma+(-1)^{\Abs0y(1+\Abs0z)}ba'\gamma.
\]
\thmref{Prop}{tan-coord} gives 
\begin{equation}\label{eq:heis-isotropy}
  f^\sharp\circ\ad^*(v)=
  \begin{cases}
    \displaystyle-\gamma\,f^\sharp\circ\frac\partial{\partial y}, &\text{ if }v=x,\\
    \displaystyle(-1)^{\Abs0y(1+\Abs0z)}\gamma\, f^\sharp\circ\frac\partial{\partial x}, &\text{ if }v=y,\\
    0, &\text{ if }v=z,
  \end{cases}
\end{equation}
where we use $x,y,z$ as a coordinate system on $\aff(\ger g_\Bbbk^*)$. Let $t\in T_0$. The image of $f^\sharp\circ\ad^*(v)$ in $T_{f_0(t)}\aff(\ger g_\Bbbk^*)=(f^*\sh T_{\aff(\ger g_\Bbbk^*)})(t)$ is 
\[
  (f^\sharp\circ\ad^*(v))(t)=\ad^*(v)(f_0(t))=
  \begin{cases}
    \displaystyle-\gamma(t)\,\Parens2{\frac\partial{\partial y}}(f_0(t)), &\text{ if }v=x,\\
    \displaystyle(-1)^{\Abs0y(1+\Abs0z)}\gamma(t)\,\Parens2{\frac\partial{\partial x}}(f_0(t)), &\text{ if }v=y.
  \end{cases}
\]
These are zero if $\gamma(t)=0$ and linearly independent otherwise. In the latter case, condition \eqref{item:action-locconst-iii} of \thmref{Th}{action-locconst} is verified. In the former case, the images of $(f^\sharp\circ\ad^*(v))$, $v=x,y$, in $(f^*\sh A_\ger g)(t)$ are zero if and only if $\gamma_t\in\gamma_t\ger m_{T,t}$. 

For simplicity, let $T\in\SMan_\knums$ and $(\tau,\theta)$ be a local coordinate system at $t$ \scth $\tau^j(t)=0$ \fa $j$. Assume that $\gamma_t=\gamma_th_t$ \fs $h_t\in\sh O_{T,t}$, but $\gamma_t\neq0$. Then in the expansion $\gamma=\sum_J\gamma_J\theta^J$ there is some minimal $I$ \scth $\gamma_I(t)\neq0$. It follows that $\gamma_I(t)=\gamma_I(t)h_0(t)$, so that $h\notin\ger m_{T,t}$. 

Thus, applying \thmref{Th}{action-locconst}, we have proved that for $T\in\SMan_\knums$, $a_f$ has locally constant rank over $T$ if and only if 
\[
  \forall t\in T_0:\Parens1{\gamma(t)=0\ \Longrightarrow\ \gamma_t=0}.
\]
If $T_0$ is connected, then this condition is equivalent to: $\gamma\in\Gamma(\sh O_T^\times)$ or $\gamma=0$. The orbit exists if the orbit map $a_f$ attached to $f$ has locally constant rank, by \thmref{Prop}{orb-quot}. 

To compute the coadjoint action, we realise $G$ in matrix form and $\ger g$ as left-invariant vector fields on $G$. For any $R\in\ssplfg{\knums}$, consider $3\times 3$ matrices with entries in $\sh O_R$. We fix the parity on the matrices by decreeing that the rows and columns of nos.~$1,2,3$ have parities depending on those of $x,y,z$ according to Table~\ref{tab:paritydist}.

\begin{table}[h]
  \caption{Parity distribution for the supergroups of Heisenberg type\label{tab:paritydist}}
  \begin{tabular}{|c|c|c||c|c|c|}
    \hline
    $\Abs0x$&$\Abs0y$&$\Abs0z$&1&2&3\\
    \hline
    \hline
    $\mathstrut^{\displaystyle\mathstrut}\ev$&$\ev$&$\ev$&$\ev$&$\ev$&$\ev$\\
    $\odd$&$\odd$&$\ev$&$\odd$&$\ev$&$\odd$\\
    $\ev$&$\odd$&$\odd$&$\ev$&$\ev$&$\odd$\\
    $\odd$&$\ev$&$\odd$&$\odd$&$\ev$&$\ev$\\
    \hline
  \end{tabular}
\end{table}

Then matrices of the form 
\[
  \begin{Matrix}1
    1&a'&c'\\ 0&1&b'\\ 0&0&1    
  \end{Matrix}
\]
are even if and only if $\Abs0{a'}=\Abs0x$, $\Abs0{b'}=\Abs0y$, and $\Abs0{c'}=\Abs0z$. Let $G'(R)$ be the set of these matrices where in addition $\{a',b',c'\}\subseteq\Gamma(\sh O_{R,\reals})$. Clearly, by defining the group multiplication by the multiplication of matrices, $G'$ is the point functor of a Lie supergroup. As we shall show presently, it is isomorphic to $G$. Since $G'_0=G_0$ is the additive group of $\reals$, unless $G$ is purely even, it will be sufficient to show that the Lie superalgebra of left-invariant vector fields on $G'$ is precisely $\ger g$.

Let $(a,b,c)$ be the coordinate system on $G$ defined on points by 
\[
  h%
  \begin{Matrix}1
    1&a'&c'\\ 0&1&b'\\ 0&0&1
  \end{Matrix}
  \defi
  \begin{cases}
    (-1)^{\Abs0x}a'&h=a,\\
    (-1)^{\Abs0y}b'&h=b,\\
    (-1)^{\Abs0z}c'&h=c.
  \end{cases}
\]
Note that this sign convention is natural in the following sense: Consider the supermanifold $G'$ as the affine superspace of strictly upper triangular matrices. Then $a,b,c$ are the linear superfunctions which constitute the dual basis to the standard basis $(E_{12},E_{13},E_{23})$ of elementary matrices.

Let $\frac\partial{\partial a},\frac\partial{\partial b},\frac\partial{\partial c}$ be the coordinate vector fields given by the coordinate system $(a,b,c)$. Let $R_x,R_y,R_z$ be the left-invariant vector fields on $G'$ determined by
\[
  R_x(1_{G'})=\frac\partial{\partial a}(1_{G'}),\quad
  R_y(1_{G'})=\frac\partial{\partial b}(1_{G'}),\quad
  R_z\defi[R_x,R_y],
\]
where write $R_x(1_{G'})$ for $1_{G'}^\sharp\circ R_x$, \emph{etc.}

We now proceed to compute these explicitly. Let $\phi^x:*[\tau_x]\longrightarrow G'$ be the infinitesimal flow of $R_x(1_{G'})$, where $\Abs0{\tau_x}=\Abs0x$. (Compare \thmref{Def}{inf-flow}.) For any function $h$ on $G'$, we have 
\[
  \frac\partial{\partial\tau_x}\Big|_{{\tau_x}=0}
  h%
  \begin{Matrix}1
    1&(-1)^{\Abs0x}\tau_x&0\\ 0&1&0\\ 0&0&1    
  \end{Matrix}
  =
  \Parens2{\frac\partial{\partial a}h}(1_{G'}),
\]
as one sees by inserting the coordinates $h=a,b,c$. Thus, we have 
\[
  (\phi^x)^\sharp(h)=h
  \begin{Matrix}1
    1&(-1)^{\Abs0x}\tau_x&0\\ 0&1&0\\ 0&0&1    
  \end{Matrix}\!.
\]
Similarly, we obtain
\[
  (\phi^y)^\sharp(h)=h%
  \begin{Matrix}1
    1&0&0\\ 0&1&(-1)^{\Abs0y}\tau_y\\ 0&0&1    
  \end{Matrix}
\]
for the infinitesimal flow $\phi^y$ of $R_y(1_{G'})$. 

We compute
\[
  \begin{split}
    (R_xh)%
    \begin{Matrix}1
      1&a'&c'\\ 0&1&b'\\ 0&0&1
    \end{Matrix}&=
    \frac\partial{\partial\tau_x}\Big|_{\tau_x=0}h%
    \begin{Matrix}1
      1&a'+(-1)^{\Abs0x}\tau_x&c'\\ 0&1&b'\\ 0&0&1
    \end{Matrix}\!,\\
    (R_yh)%
    \begin{Matrix}1
      1&a'&c'\\ 0&1&b'\\ 0&0&1
    \end{Matrix}&=
    \frac\partial{\partial\tau_y}\Big|_{\tau_y=0}h    
    \begin{Matrix}1
      1&a'&(-1)^{\Abs0y}a'\tau_y+c'\\ 0&1&(-1)^{\Abs0y}\tau_y+b'\\ 0&0&1
    \end{Matrix}\!,
  \end{split}
\]
by again inserting the coordinates for $h$. We obtain
\begin{equation}\label{eq:rx-ry}
  R_x=\frac\partial{\partial a},\quad
  R_y=\frac\partial{\partial b}+(-1)^{\Abs0x\Abs0y}a\frac\partial{\partial c}.
\end{equation}
Here, we have used the parity identity $\Abs0x+\Abs0y+\Abs0z=\ev$. From these expressions, we see immediately that 
\begin{equation}\label{eq:rz}
  R_z=[R_x,R_y]=(-1)^{\Abs0x\Abs0y}\Bracks3{\frac\partial{\partial a},a\frac\partial{\partial c}}=(-1)^{\Abs0x\Abs0y}\frac\partial{\partial c},
\end{equation}
and that this is the only non-zero bracket between the vector fields $R_x,R_y,R_z$. The sign $(-1)^{\Abs0x\Abs0y}$ that appears in the case of $\Abs0x=\Abs0y=\odd$ is an artefact of the parity distribution which is non-standard in that case. 

Since $R_x,R_y,R_z$ are linearly independent, they span the Lie superalgebra of $G'$, and it follows that $G\cong G'$. In what follows, we will identify these two supergroups. Moreover, we will identify $x,y,z$ with $R_x,R_y,R_z$, respectively. 

For further use below, we note that the right-invariant vector fields $L_x,L_y,L_z$, defined by 
\[
  L_v\defi -i_G^\sharp\circ R_v\circ i_G^\sharp,\quad v=x,y,z,
\]
take on the form
\begin{equation}\label{eq:lx-ly-lz}
  L_x=\frac\partial{\partial a}+b\frac\partial{\partial c},\quad
  L_y=\frac\partial{\partial b},\quad
  L_z=(-1)^{\Abs0x\Abs0y}\frac\partial{\partial c}.
\end{equation}
One immediately checks the bracket relation $[L_x,L_y]=-L_z$.

We now calculate the adjoint action of $G$ in terms of the matrix presentation. Let $U\in\ssplfg{\knums}$ and $(g,v)\in_UG\times\aff^\knums(\ger g)$ (\cf Ref.~\cite{ap-sphasym} for the notation), where we write
\[
  g=
  \begin{Matrix}1
    1&a'&c'\\ 0&1&b'\\ 0&0&1
  \end{Matrix},\quad
  v=\xi x(1_G)+\eta y(1_G)+\zeta z(1_G)\in\Gamma((1_G(g)^*\sh T_G)_\ev).
\]
According to the definition of $a$, $b$, and $c$, the generic point ${\id}_G\in_GG$ is 
\[
  {\id}_G=
  \begin{Matrix}1
    1&(-1)^{\Abs0x}a&(-1)^{\Abs0z}c\\ 0&1&(-1)^{\Abs0y}b\\ 0&0&1
  \end{Matrix}\!.
\]
Denote the diagonal morphism of $U$ by $\Delta_U$. We compute, for any function $h$ on $G$:
\begin{align*}
  \Ad(g)(v)(h)&=\Delta_U^\sharp(1\otimes v\otimes 1)h(g\,{\id}_G\,g^{-1})\\
  &=
  \Delta_U^\sharp(1\otimes v\otimes 1)\,h\!%
  \begin{Matrix}1
    1&(-1)^{\Abs0x}a&(-1)^{\Abs0z}c+(-1)^{\Abs0y}a'b-(-1)^{\Abs0x}ab'\\
    0&1&(-1)^{\Abs0y}b\\
    0&0&1
  \end{Matrix}\!.
\end{align*}
To evaluate this further, we insert $a,b,c$ for $h$. For $h=a,b$, Equation \eqref{eq:rx-ry} tells us that we get $\xi$ and $\eta$, respectively. For $h=c$, we get, upon applying Equation \eqref{eq:rz}:
\[
  (-1)^{\Abs0x\Abs0y}\zeta+(-1)^{\Abs0x(\Abs0y+\odd)}\eta a'-(-1)^{\Abs0y}\xi b'.
\]
Thus, identifying $x$ with $x(1_G)$, \emph{etc.}, and writing $v$ in columns, we find 
\[
  \Ad\!
  \begin{Matrix}1
    1&a'&c'\\ 0&1&b'\\ 0&0&1
  \end{Matrix}\!
  \begin{Matrix}1
    \xi x\\ \eta y\\ \zeta z   
  \end{Matrix}
  =
  \begin{Matrix}1
    \xi x\\ \eta y\\ (\zeta+(-1)^{\Abs0x}\eta a'-(-1)^{(\Abs0x+\odd)\Abs0y}\xi b')z
  \end{Matrix}\!.
\]
One may verify the correctness of this result by rederiving the bracket relation
\begin{align*}
  [x,y]&=\frac\partial{\partial\tau_y}\Big|_{\tau_y=0}(-1)^{\Abs0x\Abs0y}[x,\tau_yy]\\
  &=\frac{\partial^2}{\partial\tau_y\partial\tau_x}\Big|_{\tau_x=\tau_y=0}(-1)^{\Abs0x\Abs0y}\Ad
  \begin{Matrix}1
    1&(-1)^{\Abs0x}\tau_x&0\\ 0&1&0\\ 0&0&1
  \end{Matrix}
  (\tau_yy)\\
  &=\frac{\partial^2}{\partial\tau_y\partial\tau_x}\Big|_{\tau_x=\tau_y=0}(-1)^{\Abs0x\Abs0y}(-1)^{\Abs0x\Abs0y+\Abs0x}\tau_y((-1)^{\Abs0x}\tau_x)=z.
\end{align*}
It is now straightforward if somewhat tedious to derive 
\begin{equation}\label{eq:Adstar-explicit}
  \Ad^*\!%
  \begin{Matrix}1
    1&a'&c'\\ 0&1&b'\\ 0&0&1
  \end{Matrix}\!%
  \begin{Matrix}1
    \xi^*x^*\\
    \eta^*y^*\\
    \zeta^*z^*
  \end{Matrix}
  =
  \begin{Matrix}1
    (\xi^*+(-1)^{\Abs0y(\Abs0x+\odd)}b'\zeta^*)x^*\\
    (\eta^*-(-1)^{\Abs0x}a'\zeta^*)y^*\\
    \zeta^*z^*
  \end{Matrix}  
\end{equation}
for any 
\[
  (g,v^*)\in_UG\times\aff^\knums(\ger g^*),\quad
  g=
  \begin{Matrix}1
    1&a'&c'\\ 0&1&b'\\ 0&0&1
  \end{Matrix},
  \quad 
  v^*=
  \begin{Matrix}1
    \xi^*x^*\\ \eta^*y^*\\ \zeta^*z^*
  \end{Matrix}\!.
\]
As for the adjoint action, we make a sanity check: 
\begin{align*}
  \ad^*(x)(z^*)&=\frac{\partial^2}{\partial\tau_z\partial\tau_x}\Big|_{\tau_x=\tau_z=0}(-1)^{\Abs0x\Abs0z}\Ad^*%
  \begin{Matrix}1
    1&(-1)^{\Abs0x}\tau_x&0\\ 0&1&0\\ 0&0&1
  \end{Matrix}(\tau_zz^*)\\
  &=\frac{\partial^2}{\partial\tau_z\partial\tau_x}\Big|_{\tau_x=\tau_z=0}(-1)^{\Abs0x\Abs0z+\Abs0x}(-(-1)^{\Abs0x}\tau_x)\tau_zz^*=-(-1)^{\Abs0x\Abs0z}z^*,
\end{align*}
which is in agreement with our previous computations. 

Let us return to our $T$-valued point $f$ in the case where $\alpha=\beta=0$, \ie we have $f=\gamma z^*\in_T\aff(\ger g^*_\reals)$. Then 
\begin{equation}\label{eq:heisen-isotropy-group}
  (t,g)\in_UG_f\ \Longleftrightarrow\ a't^\sharp(\gamma)=b't^\sharp(\gamma)=0,\quad
  g=\begin{Matrix}11&a'&c'\\ 0&1&b'\\ 0&0&1\end{Matrix}\!.
\end{equation}
Moreover, the orbit map $a_f:G_T\longrightarrow\aff(\ger g_\reals^*)$ takes the form 
\[
  (a_f)^\sharp(x)=\gamma b,\quad
  (a_f)^\sharp(y)=-(-1)^{\Abs0x\Abs0z}\gamma a,\quad
  (a_f)^\sharp(z)=\gamma,
\]
in terms of coordinates $a,b,c$ on $G$ and the (linear) coordinates $x,y,z$ on $\aff(\ger g_\reals^*)$, given by
\[
  h%
  \begin{Matrix}1
    \xi^*x^*\\ \eta^*y^*\\ \zeta^*z^*
  \end{Matrix}
  =
  \begin{cases}
      (-1)^{\Abs0x}\xi^*x(x^*)=\xi^*&h=x,\\
      (-1)^{\Abs0y}\eta^*y(y^*)=\eta^*&h=y,\\
      (-1)^{\Abs0z}\zeta^*z(z^*)=\zeta^*&h=z.
    \end{cases}  
\]

We will now analyse this further, separately in the two cases in which $G$ is not a Lie group (\ie when at least one of $x,y,z$ is odd).

\subsubsection{\texorpdfstring{The Clifford supergroup of dimension $1|2$}{The Clifford supergroup of dimension 1|2}}

Assume that $\Abs0x=\Abs0y=\odd$. In this case, $G$ is called the \Define{Clifford supergroup}. This case has been given a definitive treatment by Neeb and Salmasian \cites{ns11,salmasian}, see also Ref.~\cite{ahl-cliff} for the related harmonic analysis. Our emphasis here will be to put it in the general context the orbit method. Moreover, we shall obtain the full family of Clifford modules for any non-trivial central character in one sweep. 

We will take $T\defi\aff^1\setminus0$ and $\gamma\defi u$, the standard coordinate function on $\aff^1$, so $f=\gamma z^*:T\longrightarrow\aff(\ger g_\Bbbk^*)$. Since $\gamma$ is invertible, $a_f$ has locally constant rank over $T$, and in particular, $G_f$ is a Lie supergroup over $T$.

It is completely determined by its underlying Lie group $(G_f)_0$ over $T_0$ and its Lie superalgebra $\ger g_f$ (over $\sh O_T$), defined by 
\[
  \ger g_f(U)\defi\Set2{v=\textstyle\sum_jv^je_j\in\sh O_T(U)\otimes_\knums\ger g}{\textstyle\sum_jv^j(1_G^\sharp\circ e_j\circ a_f^\sharp)=\sum_jv^jf^\sharp\circ a_{e_j}=0},
\]
for any open set $U\subseteq T_0$. In view of Equation \eqref{eq:heis-isotropy}, we have $\ger g_f=\sh O_Tz$. For the superspace $U=*$, the condition in Equation \eqref{eq:heisen-isotropy-group} is void. We conclude that the point functor of $G_f$ is given by 
\[
  G_f(U)=\Set2{\Parens2{t,\begin{Matrix}0
    1&0&c'\\ 0&1&0\\ 0&0&1
  \end{Matrix}}}{t,c'\in\Gamma(\sh O_{U,\reals,\ev})},
\]
\fa $U\in\ssplfg{\knums}$, so that $G_f\cong\aff^1_T$ with the standard addition of $\aff^1$ as multiplication over $T$.

The orbit $G\cdot f=G_T/G_f$ is $T\times\aff^{0|2}$ with fibre coordinates $a,b$. The local embedding $\tilde a_f:G\cdot f\longrightarrow\aff(\ger g_\Bbbk^*)_T$ over $T$ is given by 
\[
  (\tilde a_f)^\sharp(x)=\gamma b,\quad
  (\tilde a_f)^\sharp(y)=-\gamma a,\quad
  (\tilde a_f)^\sharp(z)=\gamma.
\]

Again following the general philosophy of geometric quantisation or Kirillov's orbit method, we choose a polarising subalgebra. To avoid reality problems, we consider the case of $\knums=\cplxs$. In the real case, we would have to complexify anyway.

A polarising subalgebra corresponds here to the preimage $\ger h$ in $\ger g_T=\sh O_T\otimes\ger g$ of a locally direct submodule of $\ger g_T/\ger g_f$ which is maximally totally isotropic with respect to the supersymplectic form induced by $\omega_f$. We will consider the case of
\[
  \ger h\defi\Span0{x,z}_{\sh O_T}.
\]
The image in $\ger g_T/\ger g_f$ is indeed maximally totally isotropic. 

The space of $\ger h$-polarised sections of the canonical line bundle on $G\cdot f$ is the $\sh O$-submodule $\sh H_f^\ger h$ of $\sh H_{G_T}$ introduced in \thmref{Prop}{pol-sect}. It is given by 
\[
  \sh H(t)\defi\Set1{\psi}{\psi\in\Gamma(\sh O_{G_U,\ev}),R_x\psi=0,R_z\psi=-it^\sharp(\gamma)\psi},
\]
for $U\in\ssplfg{\cplxs}$, $t\in_UT$. By Equations \eqref{eq:rx-ry} and \eqref{eq:rz}, this amounts to 
\[
  \psi=\vphi e^{it^\sharp(\gamma)c}
\]
where $\vphi\in\Gamma(\sh O_{U\times\aff^{0|1},\ev})$, and we consider $b$ as fibre coordinate on $(U\times\aff^{0|1})/U$. Thus, $\psi$ admits an expansion in the powers $b^0,b^1$ of $b$, with coefficients in functions on $U$. Thus, $\sh H$ is the functor of the free $\sh O_T$-module of rank $1|1=\dim\Gamma(\sh O_{\aff^{0|1}})$. We denote the corresponding $\sh O_T$-module by the same letter. 

Denoting the restriction of $\lambda_{G_T}$ to $\sh H$ by $\pi$, we compute for $g=\begin{Matrix}01&a'&c'\\ 0&1&b'\\ 0&0&1\end{Matrix}\in_UG$:
\[
  \begin{split}
    \pi(g)\psi&=(({\id}_U,g^{-1})\times m_G)^\sharp\Parens1{\vphi(b_1+b_2)e^{it^\sharp(\gamma)(c_1+c_2+a_1b_2)}}\\
    &=\vphi(b-b')e^{it^\sharp(\gamma)(-c'+a'b'+c-a'b)}=e^{it^\sharp(\gamma)(-a'b+a'b'-c')}\psi(b-b')
  \end{split}
\]
Formally deriving this expression, we readily obtain the infinitesimal action
\[
  d\pi(x)=-ibt^\sharp(\gamma),\quad d\pi(y)=-\frac\partial{\partial b},\quad d\pi(z)=it^\sharp(\gamma).
\]
Since the supercommutator of $\pi(x)$ and $\pi(y)$ is an anticommutator, we recognise this as the `fermionic Fock space' or `spinor module' of the $\sh O_T$-Clifford algebra $\mathop{\mathrm{Cliff}}(2,\sh O_T)\defi(\sh O_T\otimes\Uenv0{\ger g})/(z-i\gamma\cdot 1)$. That is, we have a trivial bundle of `spinor' modules $\cplxs^{1|1}$ over the base space $\reals^\times$, where the central character on the fibre at $t\in\reals^\times$ is $i\gamma(t)=it$. (The fibres are unital algebra representations of $\mathop{\mathrm{Cliff}}(2,\cplxs)$.)

\subsubsection{\texorpdfstring{The odd Heisenberg supergroup of dimension $1|2$}{The odd Heisenberg supergroup of dimension 1|2}}

Assume now that $\Abs0x=\ev$, $\Abs0y=\Abs0z=\odd$. In this case, we call $G$ the \Define{odd Heisenberg supergroup}, since it is a central extension of the Abelian Lie supergroup $\aff^{1|1}$ with respect to a $2$-cocycle corresponding to an odd supersymplectic form. 

We will take $T\defi\aff^{0|1}$ and $\gamma\defi\theta$, the standard coordinate function on $\aff^{0|1}$, so $f=\theta z^*:T\longrightarrow\aff(\ger g_\Bbbk^*)$. In this case, Equation \eqref{eq:heisen-isotropy-group} gives 
\[
  G_f=(\reals,\sh O_{G_f}),\quad\sh O_{G_f}\defi\sh O_{\aff^1}[b,c,\theta]/(a\theta,b\theta),
\]
where $b,c$ are odd, $a$ is the standard coordinate function on $\aff^1$, and the embedding $j:G_f\longrightarrow G_T$ is the obvious one. Clearly, $G_f$ is not a supermanifold over $T$.

To determine the orbit, let $h$ be a function on $G_T$ and expand 
\[
  h=h_0+h_bb+h_cc+h_\theta\theta+h_{bc}bc+h_{b\theta}b\theta+h_{c\theta}c\theta+h_{bc\theta}bc\theta
\]
where $h_I$ are functions on $\aff^1$. The multiplication $m$ of $G_T$ is given by 
\[
  m^\sharp(a)=a_1+a_2,\quad m^\sharp(b)=b_1+b_2,\quad m^\sharp(c)=c_1+c_2+a_1b_2,
\]
where we write $a_i\defi p_i^\sharp(a)$, \etc{} Thus, writing $m'\defi m\circ({\id}_{G_T}\times_Tj)$, we find that 
\[
  \begin{split}
    m^{\prime\sharp}(h)=&\;m^{\prime\sharp}(h_0)+m^{\prime\sharp}(h_b)(b_1+b_2)+m^{\prime\sharp}(h_c)(c_1+c_2+a_1b_2)+m^{\prime\sharp}(h_\theta)\theta\\
    &+m^{\prime\sharp}(h_{bc})(b_1c_1+b_1c_2+a_1b_1b_2-c_1b_2+b_2c_2)+m^{\prime\sharp}(h_{b\theta})b_1\theta\\
    &+m^{\prime\sharp}(h_{c\theta})(c_1\theta+c_2\theta)+m^{\prime\sharp}(h_{bc\theta})(b_1c_1\theta+b_1c_2\theta).
  \end{split}
\]
Since $p_1^\sharp(h)$ contains only $b_1,c_1$, if $h$ is invariant, then all summands involving $b_2$ or $c_2$ have to vanish. Moreover, on $G_T\times_TG_f$, we have 
\[
  m^{\prime\sharp}(h_I)\theta=h_I(a_1+a_2)\theta=h_I(a_1)\theta+h_I'(a_1)a_2\theta=h_I(a_1)\theta=p_1^\sharp(h_I)\theta,
\]
so the invariance condition is verified automatically for the $\theta$ and $b\theta$ components. Therefore, $h$ is left $G_f$-invariant if and only if 
\[
  \begin{cases}
    m^{\prime\sharp}(h_I)=p_1^\sharp(h_I), &\text{ for }I=0,\\
    h_I=0,&\text{ for }I=b,c,bc,c\theta,bc\theta.
  \end{cases}
\]
In other words, $h$ is of the form 
\[
  h=h_0+h_\theta\theta+h_{b\theta}b\theta
\]
where $h_0$ is constant and $h_\theta$, $h_{b\theta}$ are arbitrary. It follows that the colimit in $\SSp_T$ of $m,p_1:G_T\times_TG_f\longrightarrow G_T$ is given by 
\[
  Q\defi\Parens0{*,\sh O_Q},\quad\sh O_Q=\Set1{f\in\Gamma(\sh O_{\aff^1})[\eps|\theta]/(\eps^2,\eps\theta)}{f_0\in\knums},
\]
together with the morphism $\pi_f:G_T\longrightarrow Q$ determined by 
\[
  \pi_f^\sharp(a)=a,\quad\pi_f^\sharp(\eps)=b\theta,\quad\pi_f^\sharp(\theta)=\theta,
\]
see \thmref{Rem}{colimit-explicit}. By \thmref{Prop}{wgeom-cat}, $Q$ is a regular superspace in the sense of \cite{ahw-sing}*{Definition 4.12}, but it is not locally finitely generated, because it is not a subspace of $Y_q\defi(*,\knums\lBr a\rBr[\theta^1,\dotsc,\theta^q])$ for any $q$. (If $Q$ were locally finitely generated, then it would have to be a subspace of some $Y_q$ \cite{ahw-sing}*{Example 3.50}.) 

Nonetheless, we have the subrepresentation $\sh H_f^\ger h$ of $\sh H_{G_T}$ from \thmref{Prop}{pol-sect} for polarising subalgebras $\ger h\subseteq\ger g_T$. We choose
\[
  \ger h\defi\Span0{x,z}_{\sh O_T}.
\]
Once again, $\sh H=\sh H_f^\ger h$ is given by 
\[
  \sh H(t)\defi\Set1{\psi}{\psi\in\Gamma(\sh O_{G_U,\ev}),R_x\psi=0,R_z\psi=-it^\sharp(\gamma)\psi}
\]
for $U\in\ssplfg{\knums}$, $t\in_UT$. We see that the condition on $\psi$ amounts to 
\[
  \psi=\vphi e^{it^\sharp(\gamma)c}=\vphi(1+it^\sharp(\gamma)c),
\]
where $\vphi\in\Gamma(\sh O_{U\times\aff^{0|1},\ev})$ admits a finite expansion in $b$ with coefficients in functions on $U$. Again, this corresponds to the $\sh O_T$-module $\sh O_T\otimes\smash{\cplxs^{1|1}}$. The restriction $\pi$ of $\lambda_{G_T}$ to $\sh H$ is given by the same formula as before:
\[
  \pi(g)\psi=e^{it^\sharp(\gamma)(-a'b+a'b'-c')}\psi(b-b'),\quad\forall g=\begin{Matrix}01&a'&c'\\ 0&1&b'\\ 0&0&1\end{Matrix}.
\]
Formally deriving this expression, one obtains the following infinitesimal action:
\[
  d\pi(x)=-it^\sharp(\gamma)b,\quad d\pi(y)=-\frac\partial{\partial b},\quad d\pi(z)=it^\sharp(\gamma).
\]
Since the supercommutator of $d\pi(x)$ and $d\pi(y)$ is an ordinary commutator, this is a parity reversed Schr\"odinger representation, parametrised by $T=\aff^{0|1}$. 

If instead we consider the polarising subalgebra $\ger h=\Span0{y,z}_{\sh O_T}$, then the dimension of the representation $\sh H_f^\ger h$ changes drastically, although the action is formally very similar. (Essentially, $a$ and $b$ exchange their roles.) This has also been observed by Tuynman \cite{tuyn10b} in his setting and seems to reflect the fact that in this case, the orbit is not a supermanifold. 

\appendix

\section{The inverse function theorem over a singular base}\label{app:invfun}

In this appendix, we prove a relative version of the inverse function theorem, valid for an arbitrary base $S\in\ssplfg{\knums}$. This was used heavily in Section \ref{sec:super-quot}. The case where $S$ is a supermanifold is a corollary of the well-known inverse function theorem for supermanifolds \cite{leites}*{Theorem 2.3.1}. However, the proof in that case does not apply without change to such cases as $S=\Spec\knums\llbracket T\rrbracket$, which is covered by our argument. 

\begin{Th}[invfun-loc][inverse function theorem]
  Let $X/S$ and $Y/S$ be in $\SMan_S$ and $\vphi:X/S\longrightarrow Y/S$ a morphism over $S$. For any $x\in X_0$, the following are equivalent:
  \begin{enumerate}
    \item\label{item:invfun-loc-i} There is an open neighbourhood $U\subseteq X_0$ of $x$ so that $V\defi\vphi_0(U)\subseteq Y_0$ is open, and $\vphi:X|_U\defi(U,\sh O_X|_U)\longrightarrow Y|_V$ is an isomorphism.
    \item\label{item:invfun-loc-ii} The germ $(\sh T_{\vphi/S})_x:\sh T_{X/S,x}\longrightarrow\sh T_{Y/S,\vphi_0(x)}$ is invertible.
    \item\label{item:invfun-loc-iii} The map $T_{S,x}\vphi:T_{S,x}X\longrightarrow T_{S,\vphi_0(x)}Y$ is invertible.
  \end{enumerate}
\end{Th}

\begin{proof}
  The only non-trivial implication is \eqref{item:invfun-loc-iii} $\Longrightarrow$ \eqref{item:invfun-loc-i}. The question is local, so that we may assume that there are globally defined fibre coordinates $(x^a)=(u^j,\xi^k)$ on $X$ and $(y^a)=(v^j,\eta^k)$ on $Y$. Consider the ideal $\sh I_X\subseteq\sh O_X$ that is the tidy closure of that generated by the $\xi^k$ and $\smash{p_{X,0}^{-1}(\sh I_{S_0})}$, where $p_X:X\longrightarrow S$ is the structural projection and $\smash{\sh I_{S_0}}$ is the ideal of the reduction of $S$ \cite{ahw-sing}*{Construction 3.9}. (Here, the notion of \Define{tidy closure} of an ideal is introduced in \cite{ahw-sing}*{Definition 3.40}, where it was called \Define{tidying}.) The similarly defined ideal of $\sh O_Y$ shall be denoted by $\sh I_Y$.

  Let $j_X:X^{(0)}\longrightarrow X$ and $j_Y:Y^{(0)}\longrightarrow Y$ be the thickenings \cite{ahw-sing}*{Definition 2.10} defined by $\sh I_X$ and $\sh I_Y$, respectively. Let $\smash{X^{(n)}}$ and $\smash{Y^{(n)}}$, respectively, be the tidyings of the $k$th normal neighbourhoods of $j_X$ and $j_Y$ (see \cite{ahw-sing}*{Proposition 3.52}, where the notation is different). That is, $\smash{X^{(n)}}=(X_0,\smash{\sh O_X/\sh I_X^{(n+1)}})$ where $\smash{\sh I_X^{(n+1)}}$ is the tidy closure of $\smash{(\sh I_X)^{n+1}}$. There are natural tidy embeddings $\smash{j^{(n)}_X}:\smash{X^{(n)}}\longrightarrow X$ and $\smash{j^{(n+1,n)}_X}:\smash{X^{(n)}}\longrightarrow\smash{X^{(n+1)}}$ \scth the $\smash{j_X^{(n+1,n)}}$ form an inductive system. 

  Since the morphism $\vphi$ is over $S$, we have $\vphi^\sharp(\sh I_Y)\subseteq\sh I_X$ \cite{ahw-sing}*{Proposition 3.47, Corollary 3.49}. Therefore, we obtain commutative diagrams
  \[
    \begin{tikzcd}[column sep=large]
      X\rar{\vphi}&Y\\
      X^{(n+1)}\uar{j_X^{(n+1)}}\rar{\vphi^{(n+1)}}&Y^{(n+1)}\uar[swap]{j_Y^{(n+1)}}\\
      X^{(n)}\uar{j_X^{(n+1,n)}}\rar{\vphi^{(n)}}&Y^{(n)}\uar[swap]{j_Y^{(n+1,n)}}
    \end{tikzcd}
  \]
  Notice that by \cite{ahw-sing}*{Corollary 5.30}, $X^{(0)}$ is the reduction $X_0$ of $X$, so it is reduced and a supermanifold over $S_0$ of even fibre dimension. 

  Assume for the moment that we can prove the theorem for such spaces. For a while, we will proceed similar to the standard proof of the inverse function theorem for supermanifolds \cite{leites}*{proof of Theorem 2.3.1}. Namely, it is easily seen that possibly after shrinking $X$ and $Y$, $\smash{\vphi^{(0)}}$ is an isomorphism over $S_0$. Again shrinking $X$ and $Y$ as necessary, we may assume that there are functions $v'^j$ on $Y$ \scth 
  \[
    j_Y^{(0)\sharp}(v'^j)=\vphi^{(0){-1}\sharp}(j_X^{(0)\sharp}(u^j)).
  \]
  Here, we abbreviate $\smash{((\vphi^{(0)})^{-1})^\sharp}$ by $\smash{\vphi^{(0){-1}\sharp}}$. Moreover, define $\smash{A^{k\ell}\defi\frac{\partial\vphi^\sharp(\eta^k)}{\partial\xi^k}}$. This matrix is invertible, so consider its inverse $A_{k\ell}$. After shrinking $X$ and $Y$ further, there are functions $A'_{k\ell}$ on $Y$ \scth
  \[
    j_Y^{(0)\sharp}(A'_{k\ell})=\vphi^{(0){-1}\sharp}(j_X^{(0)\sharp}(A_{k\ell})).
  \]
  We let $\smash{\eta'^k\defi\sum_\ell A'_{k\ell}\eta^\ell}$. Since $j_Y$ is a thickening, the values of $v'^j$ are determined by those of $\smash{j_Y^{(0)\sharp}(v'^j)}$, and hence, the mapping condition for the functions $y'^a\defi(v'^j,\eta'^j)$ is verified. Thus, there is by \cite{ahw-sing}*{Corollary 5.36} a unique morphism $\psi:Y\longrightarrow X$ over $S$ \scth 
  \[
    \psi^\sharp(x^a)=y'^a.
  \]
  Notice that 
  \[
    j_X^{(0)\sharp}\Parens1{\vphi^\sharp(\psi^\sharp(u^j))}=\vphi^{(0)\sharp}(\vphi^{(0)-1\sharp}(j_X^{(0)\sharp}(u^j)))=j_X^{(0)\sharp}(u^j),
  \]
  so that 
  \[
    (\psi\circ\vphi)^\sharp(u^j)\equiv u^j\pod{\sh I_X}.
  \]
  But since both sides of the equation are even, the equivalence is modulo $\sh I_{X,\ev}$, which is the tidy closure of the ideal generated by $\smash{p_{X,0}^{-1}(\sh I_{S_0,\ev})}$ and the $\xi^k\xi^\ell$. 
  
  We argue similarly for the $A_{k\ell}$: 
  \[
    j_X^{(0)\sharp}(\vphi^\sharp(A'_{k\ell}))=\vphi^{(0)\sharp}\Parens1{\vphi^{(0)-1\sharp}(j_X^{(0)\sharp}(A_{k\ell}))},
  \]
  so that $\vphi^\sharp(A_{k\ell}')\equiv A_{k\ell}$ modulo $\sh I_{X,\ev}$. Since $A_{km}$ is even and $\frac\partial{\partial\xi^\ell}$ is a vector field over $S$, we find that 
  \[
    \frac{\partial A_{km}}{\partial\xi^\ell}\equiv0\pod{\sh I_X}.
  \]
  Therefore we have modulo $\sh I_X$:
  \[
    \frac{\partial\vphi^\sharp(\psi^\sharp(\xi^k))}{\partial\xi^\ell}
    =\sum_m\frac{\partial\vphi^\sharp(A'_{km}\eta^m)}{\partial\xi^\ell}
    \equiv\sum_mA_{km}\frac{\partial\vphi^\sharp(\eta^m)}{\partial\xi^\ell}=\delta_{k\ell},
  \]
  where we use the simple fact that a vector field on a tidy superspace that leaves an ideal invariant also leaves its tidy closure invariant. Thus,
  \[
    (\psi\circ\vphi)^\sharp(\xi^k)\equiv\xi^k\pod{\sh J},
  \]
  where $\sh J$ is the tidy closure of the ideal generated by $\smash{p_{X,0}^{-1}(\sh I_{S_0,\odd})}$ and the $\xi^k\xi^\ell\xi^m$.

  This implies that $\smash{(\psi\circ\vphi)^\sharp}={\id}+\delta$ where $\delta$ annihilates $\smash{p_{X,0}^{-1}(\sh O_S)}$ (because $\vphi$ and $\psi$ are over $S$), $\delta(\sh O_X)\subseteq\sh I_X$ and $\delta(\sh I_X)\subseteq\sh I_X^2$. The identity
  \[
    \delta(fg)=\delta(f)g+f\delta(g)+\delta(f)\delta(g)
  \]
  shows by induction that $\delta(\sh I_X^k)\subseteq\sh I_X^{k+1}$. At this point, we cannot conclude as for the case that $S$ is a supermanifold (compare \cite{leites}*{proof of Theorem 2.3.1}), since $\sh I_X$ may not be nilpotent. However, we can continue as follows. 

  Since the morphism $\psi$ is over $S$, it induces morphisms $\psi^{(n)}:Y^{(n)}\longrightarrow X^{(n)}$ \scth $\smash{(\psi^{(n)}\circ\vphi^{(n)})^\sharp={\id}+\delta^{(n)}}$ for some sheaf endomorphisms $\delta^{(n)}$ which map $\smash{\sh I_X^k/\sh I_X^{(n+1)}}$ to $\smash{\sh I_X^{k+1}/\sh I_X^{(n+1)}}$. Thus, $\smash{(\delta^{(n)})^{n+1}}=0$, and it follows that $\smash{(\psi^{(n)}\circ\vphi^{(n)})^\sharp}$ is invertible. Thus, $\smash{\vphi^{(n)}}$ admits a left inverse $\smash{\phi^{(n)}}$ (say). By construction, we see that
  \[
    \phi^{(n+1)}\circ j^{(n+1,n)}_Y=j_X^{(n+1,n)}\circ\phi^{(n)}.
  \]

  Now, since $X$ is formally Noetherian \cite{ahw-sing}*{Lemmas 3.36, 3.37}, it follows from \cite{ahw-sing}*{Proposition 3.52} that $X=\smash{\varinjlim_nX^{(n)}}$ in the category $\ssplfg{\knums}$. Thus, $\vphi$ admits the left inverse $\phi\defi\smash{\varinjlim_n\phi^{(n)}}$. Applying the above procedure to $\phi$, it follows that $\phi$ admits a left inverse, too. But it also has a right inverse, namely, $\vphi$, so it is invertible. Hence, $\vphi$ is invertible.

  \medskip\noindent
  It remains to prove the theorem in the case where $S=S_0$ is reduced and the fibre dimension of $X$ is purely even. Possibly shrinking $S$, there is an embedding $i:S\longrightarrow S'=\aff^r$. Let $X'\defi S'\times\aff^p$ and similarly for $Y$. Define $i_X:X\longrightarrow X'$ to be the unique morphism over $i$ \scth $\smash{i_X^\sharp(x'^a)}=x^a$, where $(x'^a)$ are the standard fibre coordinates on $X'$ over $S'$. Similarly, define $i_Y:Y\longrightarrow Y'$. Then $i_X$ and $i_Y$ are embeddings by \cite{ahw-sing}*{Corollary 5.29}.

  Define $\vphi^a\defi\smash{\vphi^\sharp(y^a)}$. Possibly shrinking $X$, $X'$, $Y$, and $Y'$, we may assume that there are functions $\vphi'^a$ on $X$ \scth $\vphi^a=\smash{i_X^\sharp(\vphi'^a)}$. Since $X$ is reduced, $\smash{i_X^\sharp}$ is post-composition with $(i_X)_0$. Thus, by taking real parts in the case of $(\knums,\Bbbk)=(\cplxs,\reals)$, we may assume that the functions $\vphi'^a$ are $\Bbbk$-valued. Therefore, there is by \cite{ahw-sing}*{Corollary 5.36} a unique morphism $\vphi':X'\longrightarrow Y'$ over $S'$ \scth $\smash{\vphi'^\sharp(y'^a)}=\vphi'^a$. Then 
  \[
    \vphi'\circ i_X=i_Y\circ\vphi,
  \]
  so that by \thmref{Lem}{tanrk-lowersemi} below, we may assume that $\vphi'$ satisfies the assumption of \eqref{item:invfun-loc-iii}. But $X'$ and $Y'$ are ordinary manifolds, so the local invertibility of $\vphi'$ follows. Thus, $\vphi^\sharp(y^a)=i_X^\sharp(\vphi'^\sharp(y'^a))$ is a system of fibre coordinates. This proves the assertion.
\end{proof}

In the proof of the \thmref{Th}{invfun-loc}, we have used the following easy lemma.

\begin{Lem}[tanrk-lowersemi]
  Let $\vphi:X/S\to Y/S$ be a morphism of supermanifolds over $S$. For any pair $m|n$ of non-negative integers, the set 
  \[
    \Set1{x\in X_0}{\rk T_{S,x}\vphi\sge m|n}
  \]
  is open. Here, we write $p|q\sge m|n$ if and only if $p\sge m$ and $q\sge n$.
\end{Lem}

\begin{proof}
  In local fibre coordinates, $T_{S,x}\vphi$ is represented by the Jacobian matrix $\Jac[_S]0\vphi(x)$, which is a continuous function of $x$. Since the rank of the upper or lower diagonal block of a block matrix is a lower semicontinuous function and the finite intersection of open sets remains open, the assertion follows. 
\end{proof}

\section{Immersions of closed Lie groups over some base}\label{app:closed-subgrp}

Let $T\in\ssplfg{\knums}$ be reduced. The aim of this appendix is to prove the following.

\begin{Th}[imm-lie]
  Let $j:H\longrightarrow G$ be a morphism of Lie groups over $T$ which is an injective immersion and has closed image. Then $j$ is a closed embedding. 
\end{Th}

We let $\ger g\defi T\times_GT(G/T)$, the restriction of the fibrewise tangent bundle of $G$ to $T$, be the Lie algebra of $G$. It is a vector bundle over $T$ and admits a bracket. Similarly, we define $\ger h$ and consider the differential $dj:\ger h\longrightarrow\ger g$ induced by $T(j/T)$. It is an injective vector bundle morphism and therefore a closed embedding. We define $\exp_G:\ger g\longrightarrow G$ by 
\[
  \exp_G(x)\defi\exp_{G_{p_\ger g(x)}}(x),
\]
where $p_\ger g$ is the vector bundle projection of $\ger g$ and we write $G_t$ for the fibre of $G$ over $t\in T$. By the smooth dependence of the solutions of ODE on Cauchy data, $\exp_G$ is a morphism of manifolds over $T$. It is a local isomorphism of manifolds over $T$ by the inverse function theorem (\thmref{Th}{invfun-loc}). Similarly, we may define $\exp_H$. It follows that $j\circ\exp_H={\exp_G}\circ{dj}$, since this is true fibrewise.

Consider the set
\[
  \ger h'\defi\Set1{x\in\ger g}{\exp_G(\reals x)\subseteq j(H)}.
\]
From the fibrewise statement (which is classical), it follows that the fibres of $\ger h'$ are Lie subalgebras of the fibres of $\ger g$, and moreover, that $dj(\ger h_t)=\ger h'_t$ for any $t\in T$. Thus, $\ger h'$ identifies with the image of $dj$ and is therefore a vector subbundle of $\ger g$.

Fix $t\in T$. We may choose open neighbourhoods $U\subseteq\ger g$ of $0_t$ and $V_0\subseteq G$ of $1_t$ \scth $\exp_G:U\longrightarrow V_0$ is a homeomorphism. Since this holds fibrewise, it follows that 
\[
  \exp_G(U\cap\ger h')\subseteq V_0\cap H.
\]

Fix a vector bundle metric $\Dual0\cdot\cdot$ on $\ger g$ (this exists after possibly restricting to a paracompact neighbourhood of $t$ in $T$) and $E\defi(\ger h')^\perp\subseteq\ger g$ a vector subbundle complementary to $\ger h'$. We write $\Norm0{\cdot}=\sqrt{\Dual0\cdot\cdot}$. The proof of the following two lemmas is identical to the classical case \cite{hn}*{Section 9.2.3}.

\begin{Lem}[l1]
  Let $x_k\in U$, $x_k\neq0$, $\exp_G(x_k)\in j(H)$, converge to a point in the zero section of $\ger g$. Any accumulation point of $\Norm0{x_k}^{-1}x_k$ lies in $\ger h'$.
\end{Lem}

\begin{Lem}[l2]
  There is some open neighbourhood $U''\subseteq U\cap E$ of $0_t$ \scth we have $\exp_G(U'')\cap j(H)\subseteq T$.
\end{Lem}

\begin{Lem}[l3]
  Possibly after shrinking the neighbourhood $U''$, there exist open neighbourhoods $U'\subseteq U\cap\ger h'$ of $0_t$ and $V'\subseteq G$ of $1_t$ \scth 
  \[
    \phi:U'\times U''\longrightarrow V':(u',u'')\longmapsto (\exp_G(u'))(\exp_G(u''))
  \]
  is a diffeomorphism over $T$. 
\end{Lem}

\begin{proof}
  The classical case shows that fibrewise, $\phi$ fulfills the assumptions of the inverse function theorem (\thmref{Th}{invfun-loc}).
\end{proof}

We now come to the proof of the theorem.

\begin{proof}[\prfof{Th}{imm-lie}]
  We claim that $\exp_G(U')=V'\cap j(H)$. The inclusion $\subseteq$ is clear, since $U'\subseteq U$ by construction. Conversely, let $g=\phi(u',u'')\in j(H)$, where $(u',u'')\in U'\times U''$. Then 
  \[
    \exp_G(U'')\ni\exp_G(u'')=(\exp_G(u'))^{-1}g\in j(H),
  \]
  so that $g=1$ and $u''=0$ by \thmref{Lem}{l2}. Thus, $g=\exp_G(u')\in\exp_G(U')$, proving the claim. 

  Let $U\defi\exp_H(dj^{-1}(U'))$. This is a neighbourhood of $1_t$ in $H$, and after shrinking $U'$, we may assume that it is open. The claim implies 
  \[
    j(U)=\exp_G(U')=V'\cap j(H),
  \]
  so that $U$ carries the initial topology with respect to $j$. The assertion follows. 
\end{proof}

\end{document}